\documentclass[aap,preprint]{imsart}

\RequirePackage{amsmath,amsthm,amsfonts,amssymb}
\RequirePackage[numbers]{natbib}
\RequirePackage{graphicx}%% uncomment this for including figures

\usepackage{xparse,xifthen} % Argument passing for advanced macro definitions
\usepackage{cleveref} % Allows refs to environments without typing out the evironment type in the text
\usepackage{bm} % bold Greek letters
\usepackage{xcolor}

\usepackage{xcolor}
\usepackage{float}

\startlocaldefs
\theoremstyle{plain}
\newtheorem{theorem}{Theorem}
\newtheorem{corollary}[theorem]{Corollary}

\newtheorem{lemma}[theorem]{Lemma}
\newtheorem{proposition}[theorem]{Proposition}
\theoremstyle{remark}
\newtheorem*{remark}{Remark}

\newcommand{\T}{\top}
\newcommand{\vem}{\varepsilon^{\rm m}}
\newcommand{\ves}{\varepsilon^{\rm s}}
\newcommand{\veq}{\varepsilon^{\rm q}}
\newcommand{\veqd}{\varepsilon^{\rm qd}}
\newcommand{\hvem}{\hat{\varepsilon}^{\rm m}}
\newcommand{\hveq}{\hat{\varepsilon}^{\rm q}}
\newcommand{\cveqd}{\check\varepsilon^{\rm qd}}
\newcommand{\bo}{\mathbf{1}}
\newcommand{\maps}{\colon} % use in f:X->Y for correct spacing
\newcommand{\comp}{^c}
\newcommand{\im}{\mathrm{i} }
\DeclareMathOperator{\Ai}{Ai}
\newcommand{\bv}{\mathbf{v}}

\newcommand{\adot}{\,\cdot\,} % Function argument dot

\newcommand{\mR}{\mathbb{R}}
\newcommand{\mP}{\mathbb{P}}
\newcommand{\tr}{\operatorname{tr}}
\newcommand{\bxi}{\bm{\xi}}

\newcommand{\str}{{(\infty,1)}}

\renewcommand{\Re}{\operatorname{Re}}
\renewcommand{\Im}{\operatorname{Im}}

\newcommand{\R}{\mathbb{R}}
\newcommand{\C}{\mathbb{C}}
\newcommand{\Z}{\mathbb{Z}}
\newcommand{\E}{\mathbf{E}}
\def\Pr{\mathbf{P}}

\newcommand{\rmrightarrow}[1]{\xrightarrow{\mathrm{#1}}}
\newcommand{\drightarrow}{\rmrightarrow{d}}

\newcommand{\diag}{\mathrm{diag}}
\newcommand{\Var}{\mathrm{Var}}

\newcommand{\Ind}{\boldsymbol{1}}
\def\det{\mathrm{det}}
\DeclareMathOperator{\TV}{TV}
\newcommand{\psc}{p_{SC}}
\newcommand{\diff}{\mathop{}\!\mathrm{d}}

\makeatletter

\def\cstostr#1{%
  \expandafter\@gobble\detokenize\expandafter{\string#1}}

\NewDocumentCommand{\@resizedelimiter}{mmO{}m}{%
  \ifthenelse{\equal{\noexpand#3}{\noexpand\extend}}
    {\left#1#4\right#2}
    {\csname\cstostr#3l\endcsname#1 #4 \csname\cstostr#3r\endcsname#2}
}
\newcommand{\abs}{\@resizedelimiter{\lvert}{\rvert}}
\newcommand{\floor}{\@resizedelimiter{\lfloor}{\rfloor}}
\newcommand{\ceil}{\@resizedelimiter{\lceil}{\rceil}}

\NewDocumentCommand{\norm}{omo}{%
  \@resizedelimiter{\lVert}{\rVert}[#1]{#2}\IfNoValueF{#3}{_{#3}}
}

\makeatother
\endlocaldefs
\begin{document}

\begin{frontmatter}
\title{An edge CLT for the log determinant of Wigner ensembles}
%\title{A sample article title with some additional note\thanksref{T1}}
\runtitle{Edge CLT for log determinant}
%\thankstext{T1}{A sample of additional note to the title.}

\begin{aug}
%%%%%%%%%%%%%%%%%%%%%%%%%%%%%%%%%%%%%%%%%%%%%%%
%% Only one address is permitted per author. %%
%% Only division, organization and e-mail is %%
%% included in the address.                  %%
%% Additional information can be included in %%
\author[A]{\fnms{Iain M.} \snm{Johnstone}\ead[label=e1]{imj@stanford.edu}},
\author[B]{\fnms{Yegor} \snm{Klochkov}\ead[label=e2,mark]{yk376@cam.ac.uk}},
\author[B]{\fnms{Alexei} \snm{Onatski}\ead[label=e3,mark]{ao319@cam.ac.uk}}
\and
\author[A]{\fnms{Damian} \snm{Pavlyshyn}\ead[label=e4,mark]{damianp@stanford.edu}}
%%%%%%%%%%%%%%%%%%%%%%%%%%%%%%%%%%%%%%%%%%%%%%
%% Addresses                                %%
%%%%%%%%%%%%%%%%%%%%%%%%%%%%%%%%%%%%%%%%%%%%%%
\address[A]{Department of Statistics,
  Stanford University,
  \printead{e1,e4}}

\address[B]{Faculty of Economics,
  University of Cambridge,
  \printead{e2,e3}}
\end{aug}

\begin{abstract}
We derive a Central Limit Theorem (CLT) for $\log \left\vert\det \left( W_{N}-E_{N}\right)\right\vert,$ where $%
W_{N}$ is a Wigner matrix, 
%from the Gaussian Unitary or Gaussian Orthogonal Ensemble (GUE and GOE), 
and $E_{N}$ is local to the edge of the semi-circle law. Precisely,
$E_N=2+N^{-2/3}\sigma_N$ with $\sigma_N$ being either a
constant (possibly negative), or a sequence of positive real numbers,
slowly diverging to infinity so that $\sigma_N \ll \log^{2} N$. 
%For slowly growing $\sigma_N$, our proofs hold for general Gaussian $\beta$-ensembles. 
We also extend our CLT to cover spiked Wigner matrices. Our interest in the CLT is motivated by its applications to statistical testing in critically spiked models and to the fluctuations of the free energy in the spherical Sherrington-Kirkpatrick model of statistical physics.

%GUE and GOE.  
\end{abstract}

\begin{keyword}[class=MSC]
\kwd[Primary ]{60}%Probability theory and stoch proc
\kwd{62}%Statistics
\kwd{82}%Statistical mechanics, structure of matter
\kwd[; secondary ]{60B20}%Random matrices (probabilistic aspects)
\kwd{60F05}%central limit and other weak theorems
\kwd{60K35}%Interacting random processes, statistical mechanics type models, percolation theory. This classification was used in a 2020 Baik-Lee paper "Free energy of bipartite spherical Sherrington–Kirkpatrick model"
\kwd{62H15}%hypothesis testing in multivariate analysis
\kwd{82B44}%Disordered systems (random Ising models, random Schrodinger operators etc.) in equilibrium statistical mechanics. Again, this was used in the 2020 Baik-Lee.
\end{keyword}


\end{frontmatter}
%%%%%%%%%%%%%%%%%%%%%%%%%%%%%%%%%%%%%%%%%%%%%%
%% Please use \tableofcontents for articles %%
%% with 50 pages and more                   %%
%%%%%%%%%%%%%%%%%%%%%%%%%%%%%%%%%%%%%%%%%%%%%%

%\newpage

%%%%%%%%%%%%%%%%%%%%%%%%%%%%%%%%%%%%%%%%%%%%%%
%%%% Main text entry area:

%\newpage

%\texttt{From Editorial policy ``Most papers should contain an
%  Introduction which presents a discussion of the context and
%  importance of the issues they address and a clear, non-technical
%  description of the main results. The Introduction should be
%  accessible to a wide range of readers.'' }

%\fix{A proposal:
%\begin{enumerate}
%\item State result for real Wigner
%\item Motivations:  a) likelihood ratio, b) Baik-Lee transition
%\item Outline of approach: a) Gaussian, b) Wigner
%\item Related work: put our contribution in context
%\item Organization of paper
%\end{enumerate}
%}

%In the introduction below, \textcolor{cyan}{the text in cyan is old text that I suggest cutting}, and replacing by \textcolor{magenta}{the new text in magenta}.[AO]

\section{Introduction}
\label{sec:introduction}

Let $W_N = (\xi_{ij}/\sqrt{N})$ be an $N \times N$ real or complex
Wigner matrix; in particular $W_N$ is Hermitian and for $i \geq j$ the
entries $\xi_{ij}$ are independent with mean zero,
the variances $\E |\xi_{ij}|^2 = 1$ for $i \neq j$ and are bounded for
$i=j$.
Our conditions, fully specified in \cref{sec:lindeberg-swapping}, imply
that the empirical distribution of the eigenvalues
$\lambda_1 \geq \cdots \geq \lambda_N$ of $W_N$ converges to the
Wigner semi-circle law $\rho_{\rm sc}$ on $[-2,2]$ and that the
largest eigenvalue 
$\lambda_1$ converges almost surely to the right edge $2$,
see, for example, \cite{anderson2010introduction}.
%e.g. AGZ.

The logarithmic linear statistic
$$\mathcal{L}_N = \sum_1^N f(\lambda_j)
=\sum_1^N \log (E-\lambda_j)$$
arises in several applications;
we focus below in particular on statistical testing in
`spiked' models and on the fluctuation behavior of the free energy in
the spherical Sherrington-Kirkpatrick (SSK) model of statistical physics.
Suppose initially that $E > 2$ is fixed.
In this case $\mathcal{L}_N - N \int f \diff \rho_{\rm sc}$
is asymptotically Gaussian with
%\textcolor{cyan}{variance $\sigma^2(E) = \cdots $} 
finite variance that depends on the first four moments of the entries of $W_N$.
Since $f(z) = \log(E-z)$ is analytic in a neighborhood of the
semi-circle support, this follows from  general CLTs for linear
statistics, e.g. \cite{Yao2005}.

This paper concerns Gaussian behavior near, at, or just inside the
edge:
\begin{equation}
  \label{eq:edge}
  E = E_N = 2 + \sigma_N N^{-2/3},  \qquad -\gamma \leq \sigma_N \ll
  \log^2 N
\end{equation}
for some fixed $\gamma > 0$.
Our particular motivations, detailed below, lie in certain transition
zones in the spiked statistical and SSK models.
Here $E$ is sufficiently close to the edge that the functions
$f_N(z) = \log|z-E|$ do not appear to be covered even by recent
mesoscopic CLTs (e.g. \cite{LandonSosoeAAP,LiSchnelliXu}).

The basic identity
\begin{equation*}
  L_N = \sum_{j=1}^N \log |\lambda_j - E|
      = \log | \det(W_N-E)|
\end{equation*}
casts the linear statistic (now with the absolute value under the logarithm) as a log determinant, i.e. in terms of the
characteristic polynomial of $W_N$. The latter is the subject of a substantial
literature, partly reviewed in Section \ref{sec:related-work}.
In particular, as pioneered by Tao and
Vu \cite{tao2012central} for $E=0$,
for Gaussian ensembles $W_N$ drawn from GUE or GOE,
one can use the Trotter equivalence to cast the
matrix in tridiagonal Jacobi form and analyze the recurrence satisfied
by the principal minors.  Lindeberg swapping is used to extend to
Wigner matrices with four matching moments.

In this paper we
% assemble and develop tools to
carry out this program
at the edge \eqref{eq:edge}, to arrive at the following result.
%Introduce centering and scaling constants

\begin{theorem}   \label{th:main}
  Let $W_N$ be a Wigner matrix whose off-diagonal moments match GUE
  ($\alpha = 1$) or  GOE ($\alpha = 2$)
  to third order. For edge values $E = E_N$ satisfying
  \eqref{eq:edge}, we have
    \begin{equation}
  \label{eq:zW-clt1}
   \left(\log\abs{\det(W_N - E)} - \mu_N\right) / \tau_N  \stackrel{\mathrm{d}}{\rightarrow} \mathcal{N}(0,1),
\end{equation}
with
\begin{equation}
  \label{eq:centscal1}
   \mu_N  = \tfrac{1}{2} N + \sigma_N N^{1/3} - \tfrac{2}{3}
   (\sigma_N \vee 0)^{3/2} - \tfrac{1}{6}(\alpha-1) \log N,
   \qquad  \tau_N  = \sqrt{\tfrac{\alpha}{3} \log N}.
    % \tilde{\tau}_N
    %  = \sqrt{ \alpha \rho(\theta_N^{-2})},
    %  \qquad \rho(x) = \log \tfrac{1}{2}[1+(1-x)^{-1/2}].
\end{equation}
\end{theorem}

\subsection{Two motivating applications}
\label{sec:two-motiv-appl}

Although superficially unrelated, both applications involve the
spherical integral
\begin{equation}
  \label{eq:ZalphaN}
  Z_{\alpha,N}(\beta,M)
    =  \int_{S_\alpha^{N-1}} \exp \{ (\beta N/\alpha)  u^* M u \}
    (\diff u),
\end{equation}
where $(\diff u)$ denotes normalized uniform measure on the unit sphere
$S_\alpha^{N-1} = \{ x:  \| x \| = 1 \}$ in $\C^N$ for $\alpha = 1$, or
$\R^N$ for $\alpha = 2$, while $M$ is Hermitian resp. symmetric, and
$\beta > 0$.

\medskip
\textit{Testing critical spiked models.} \
Principal Components Analysis (PCA) seeks low-dimensional summaries of
high-dimensional data. In certain cases, such as genomics e.g.
\cite{Patterson2006},
it can be reasonable to approximate the covariance matrix as
$\Sigma = \sigma^2 I + F$, where $F$ has small rank.
A perennial applied question is to determine this rank, at least
approximately. The simplest version is to test for the presence of a
rank one component. Thus we assume
\begin{equation}
  \label{eq:PCA}   \tag{PCA}
  X_{N \times n} \text{ has i.i.d. columns } \sim N(0,\Sigma),
  \qquad  \Sigma = \sigma^2 I + h \mathbf{v v}^*,
  \qquad  M = X X^*/n.
\end{equation}
The largest eigenvalue $\lambda_1(M)$ is a natural test statistic, but
its utility is limited by a phase transition first exhibited for
complex data in \eqref{eq:PCA} by
Baik, Ben Arous and P\'ech\'e \cite{BaiBenPech}
in the setting of proportional
asymptotics $N/n \to y > 0$.
Below the critical value, $h < \sqrt{y} $, the largest eigenvalue,
after centering and scaling at rate $N^{-2/3}$, has a limiting
Tracy-Widom distribution, and so carries \textit{no information} about
$h$.

Onatski, Moreira and Hallin \cite{Ona2013}
showed that testing below the critical value was still possible,
using a likelihood ratio test of $H_0: h=0$ versus $H_A: h=\beta$.
The asymptotic behavior of this test depends on a logarithmic linear
statistic $\mathcal{L}_N$ for $E = E(\beta)$ located outside the edge of the
Mar\v{c}enko-Pastur bulk. It is also noted in \cite{Ona2013}
that the likelihood ratio
is exponentially small for supercritical alternatives $h >
\sqrt{y}$, but left open the behavior for alternatives $\beta$
near the critical point.

It is commonly noted that the spectra of Wishart matrices $X X^*$
exhibit behavior analogous to that of simpler symmetric Wigner
matrices. For us, the analog of \eqref{eq:PCA} specifies that
\begin{equation}
  \label{eq:SMD}   \tag{SMD}
    M = h \mathbf{v v}^* + Z/\sqrt{N},
\end{equation}
where $Z$ is in general an $N \times N$ Wigner matrix, real symmetric
or complex Hermitian.
In the special case that $Z$ is drawn from GOE resp GUE, the term
deformed G(O/U)E is used.
The BBP transition occurs at threshold $h = 1$ in these models
\cite{Peche2006a,maid07}.
%Pec06, Mai07.
The likelihood ratio against $H_A: h = \beta < 1$ was studied, along
with other spiked models in \cite{johnstone2020testing}, again in
terms of a logarithmic 
linear statistic $\mathcal{L}_N$, now with $E = \beta + 1/\beta > 2$.
Again behavior for $\beta$ near 1 was left open.

One reason for the close parallel of results for PCA and SMD is
that the joint density of the eigenvalues $\Lambda = (\lambda_i)_1^N$
of $M$ has the same form in both cases.
The joint density is found by integrating over the orthogonal group
corresponding to the eigenvectors; for a rank one spike the integral
reduces to one over $S_\alpha^{N-1}$.
The argument goes back at least to \cite{james1954}, see also
\cite[Suppl p. 6]{johnstone2020testing} and
\cite[p. 104]{muir82}.
For SMD the result is
\begin{equation*}
  p(\Lambda,h) = c(\Lambda) d(h) Z_{\alpha,N}(h,\Lambda).
\end{equation*}
The main term $Z_{\alpha,N}(h,\Lambda)$ is given by
\eqref{eq:ZalphaN}, while
$d(h) = \exp(-(N/\alpha)(h^2/2))$ and $c(\Lambda)$ though
explicit is not needed as it disappears on taking ratios for distinct
values of $h$.

For PCA, we have $d(h) = (1+h)^{-n/\alpha}$ and in the $Z_{\alpha,N}$
term, $h$ is replaced by $nh/(N(1+h)$. In view of the foregoing
remarks, we will henceforth focus on SMD.

\bigskip
\textit{SSK model.}  In the spherical version of the
Sherrington-Kirkpatrick model  studied by Kosterlitz, Thouless
and Jones \cite{Kosterlitz1976},
% KTJ76,
the vector of spins $\sigma \in \mR^N$ is constrained
to lie on the sphere $\| \sigma \|^2 = N$.
The Hamiltonian is given by $H_N(\sigma) = \sum_{i < j} M_{ij}^{\rm SSK} \sigma_i
\sigma_j$, where the couplings $M_{ij}^{\rm SSK}$ between distinct
spins are random, and in \cite{Kosterlitz1976} are independent $N^{-1/2}
\mathcal{N}(J / N^{1/2},1)$ variates.
The partition function is then $Z_N = \int e^{\beta H_N(\sigma)} \diff
\omega_N(\sigma)$ with $\diff \omega_N$ being normalized uniform
measure on the sphere $\| \sigma \|^2 = N$.
If $M^{\rm SSK}$ 
% $= (M_{ij}^{\rm SSK})$
is the corresponding symmetric
matrix, then on rescaling to the unit sphere, we have
\begin{equation*}
  Z_N = \int \exp \{ (\beta N/2) \, u^* M^{\rm SSK} u \} (\diff u)
      = Z_{2,N}(\beta, M^{\rm SSK}). 
\end{equation*}
Since the integrals depend only on the eigenvalues of $M$, this is
exactly the integral occurring in the rank one spiked GOE model
\eqref{eq:SMD}, with the sole difference that $M^{\rm SSK}$ has
vanishing diagonal.
 
Kosterlitz et. al.  evaluated the first order limiting behavior of the
free energy, 
finding a phase diagram for $(J, 1/\beta) \in \mR_+^2$ with three regions:
ferromagnetic for $J>1, \beta J >1$, 
paramagnetic for $ \beta<1, \beta J < 1$ and 
spin glass for $J < 1, \beta > 1$.
In particular, we record that
\begin{equation*}
  F_N = \frac{1}{N} \log Z_N
      \to F(\beta) =
      \begin{cases}
        \frac{1}{4} \beta^2   & \beta < 1, J < 1 \\
        \beta - \frac{1}{2} \log \beta - \frac{3}{4}   & \beta > 1, J
        < 1
      \end{cases}.
\end{equation*}

Baik and Lee \cite{baik2016fluctuations, Baik2017} studied the
second-order fluctuations of the free 
energy, making more general assumptions of Wigner type on the
distributions of the couplings $M_{ij}$.
When $\Var M_{ii} = 0$, they refer to $H_N(\sigma)$ as the spherical
SK Hamiltonian with ferromagnetic Curie-Weiss interaction.
We refer to \cite{baik2016fluctuations, Baik2017} for fuller bibliographic discussion of work around
the SSK model.
The main results of \cite{baik2016fluctuations, Baik2017} show that
$N^\gamma (F_N - F(\beta)) 
\stackrel{\rm d}{\to} 
\xi$, where in the three phases respectively
$(\gamma, \xi) = (\frac{1}{2}, \text{Gaussian}),
(1, \text{Gaussian}),$ and
$(\frac{2}{3}, \text{Tracy-Widom})$, suppressing details of the centering
and scaling of $\xi$.

The transition regions between the three phases are studied in
\cite{baik2016fluctuations, Baik2017} and \cite{BaikLeeWu}.
% BL16, 17 and BLW18.  
Two transitions are settled but the spin glass to
paramagnetic transition is left open.
We emphasize that the open case is exactly the transition relevant to studying
the likelihood ratio statistic for testing against near critical
alternatives! 
By equating variances for
$\beta$ above and below $1$, \cite{baik2016fluctuations} conjecture
that the relevant scale 
has $\beta = 1+b N^{-1/3} \sqrt{ \log N}$ for $b \in \mR$.

In a companion paper \cite{johnstone2021}, we apply our Theorem
\ref{th:main} to  verify the scaling conjectured by Baik and Lee: on this
scale, after
centering and scaling $F_{N} - F(\beta)$ converges in law to a
$b$-dependent linear combination of \textit{independent} Tracy-Widom and Gaussian
components. In turn this implies conclusions for the null distribution
of log-likelihood ratio tests of $H_0: h=J \in [0,1)$ versus 
critically spiked alternatives $H_A: h=\beta$.

The loglinear statistic $L_N$ and Theorem \ref{th:main} are basic for
this result. Briefly, 
% -- JKOP has the detail --
the standard first step
casts the spherical integral as a single contour integral
\begin{equation*}
  Z_{\alpha,N}(\beta,\Lambda)
    = C_{\alpha,N} \int_\mathcal{K} e^{(N/\alpha) G_\beta(z)} \diff z,
\end{equation*}
where $G_\beta$ involves the loglinear statistic $\mathcal{L}_N$ 
\begin{equation*}
  G_\beta(z) = (1 + bN^{-1/3}\sqrt{\log N})z - N^{-1} \sum_1^N \log(z-\lambda_j).
\end{equation*}
The contour $\mathcal{K}$ passes to the right of all eigenvalues
$\lambda_j$, and is chosen to allow Laplace approximation of the
integral. 
For $b < 0$, the vertical contour through
$\hat{\gamma}_b = 2 + b^2 N^{-2/3} \log N$ suffices, and the main
approximating term involves $L_N(\hat{\gamma}_b)$.
For $b > 0$ a keyhole contour around $\lambda_1$ and for $b=0$ a
contour of steepest descent both yield a leading approximation term
involving
\begin{equation*}
  -L_N(2) + (\beta-1)N(\lambda_1-2).
\end{equation*}
In each case Theorem \ref{th:main} along with further analysis of the
derivatives of $G_\beta$ lead to the transition theorem.
In addition, the tridiagonal Jacobi method
and Lindeberg swapping tools 
developed here is adapted to
show the asymptotic independence of $L_N(2)$ and $\lambda_1$.

\subsection{Outline of approach}
\label{sec:outline-approach}

The analysis begins with $W_N$ drawn from a Gaussian ensemble: GUE or
GOE. By a unitary/orthogonal transformation
\cite{trotter1984eigenvalue} the eigenvalues of $\sqrt{N} W_N$ are the
same as those of
% Their derivation is based on the equivalence
% \cite{trotter1984eigenvalue} of the joint distribution of the
% eigenvalues of $M_{N}$ to that of  
\begin{equation}
  \sqrt{N} \widehat{W}_{N} = 
      \begin{pmatrix}
        a_1 ~&~ b_1 ~&~ ~&~ ~&~ \\
        b_1 ~&~ a_2 ~&~ b_2 ~&~   ~&~ \\
        ~&~ b_{2} ~&~ \ddots ~&~ \ddots  ~&~\\
%        \hline
         ~&~~&~ \ddots ~&~ \ddots ~&~ b_{N-1} \\
         ~&~~&~~&~ b_{N-1} ~&~ a_{N}
    \end{pmatrix}
%   =\begin{pmatrix}
% a_{1} & b_{1} &  &  \\ 
% b_{1} & \ddots  & \ddots  &  \\ 
% & \ddots  & a_{N-1} & b_{N-1} \\ 
% &  & b_{N-1} & a_{N}%
% \end{pmatrix}%
,  \label{tridiag form earlier}
\end{equation}
where $a_{i}\sim \mathcal{N}\left(0,\alpha\right)$ and $b_{i}^{2}\sim \chi
^{2}(2i/\alpha)/(2/\alpha)$, $i=1,...,N$, are jointly independent, and
$\alpha=1$ for GUE and $\alpha=2$ for GOE. Here, by definition, \(
\chi^{2}(d) \) has the density
\( c_d^{-1}x^{d/2-1} e^{-x/2} \mathbf{1}_{x > 0} \) for $d > 0$
and $c_d = 2^{d/2} \Gamma(d/2)$.

%\textcolor{cyan}{Let $D_i$ denote the determinant of the $i$-th minor
%of $\sqrt{N}(\widehat{W}_N-E)$.
%Using the cofactor expansion yields the recurrence
%\begin{equation}\label{determinant_recurrence_nonzero}
%D_{i}=(a_{i}-E\sqrt{N})D_{i-1}-\left(
%  i-1+\sqrt{i-1}c_{i-1}\right) D_{i-2}, 
%\end{equation}
%for $i \geq 1$ with the conventions $D_0=1, D_{-1}=0$.
%Here $c_{i}=\left( b_{i}^{2}-i\right) /\sqrt{i},$ so that $\E c_{i}=0$
%and $\Var\left( c_{i}\right) =\alpha.$
%A deterministic version of such a recursion ($a_{i}=c_{i-1}=0$) has
%an explosive characteristic root $\rho_{i}^{+}$ with 
%\[
%2\rho_{i}^{\pm}=-\Big(E\sqrt{N} \pm \sqrt{E^2N-4(i-1)}\Big).
%\]
%% which are not eliminated by the scaling
%% % $E_{i}=D_{i}/\sqrt{i!}$
%% $D_{i}/\sqrt{i!}$ because $\theta_N \neq 0$.
%We therefore adopt the normalization
%\begin{equation}
%    \label{new scaling new}
%    E_i=D_i/\prod_{j=1}^{i}\vert\rho_j^{+}\vert.
%\end{equation}
%The characteristic roots of the dynamic equation describing
%$E_{i}$ can be approximated by the pair
%$\rho^{\pm}_i/\left\vert\rho^{+}_i\right\vert$. 
%% For $i>\theta^2_{N}N$, this is a
%% complex-conjugated pair located on the unit circle in $\mathbb{C}$. In
%% contrast,
%For $i\leq E^2N/4$, these characteristic roots are
%real, the first one equals minus one, and the other decreases from zero to minus one as $i$ goes from $1$ to $E^2N/4$.
%}
%
Let $D_i$ denote the determinant of the $i$-th minor
of $\widehat{W}_N-E$.
Using the cofactor expansion yields the recurrence
\begin{equation}\label{determinant_recurrence_nonzero}
  D_{i}=
  \Big(\frac{a_{i}}{\sqrt{N}}-E\Big)D_{i-1}-\frac{
  i-1+\sqrt{i-1}c_{i-1}}{N} D_{i-2}, 
% D_{i}=(a_{i}/\sqrt{N}-E)D_{i-1}-\frac{
%   i-1+\sqrt{i-1}c_{i-1}}{N} D_{i-2}, 
\end{equation}
for $i \geq 1$ with the conventions $D_0=1, D_{-1}=0$.
Here $c_{i}=\left( b_{i}^{2}-i\right) /\sqrt{i},$ so that $\E c_{i}=0$
and $\Var\left( c_{i}\right) =\alpha.$
A deterministic version of such a recursion ($a_{i}=c_{i-1}=0$) has
an explosive characteristic root $\rho_{i}^{+}$ with 
\[
2\rho_{i}^{\pm}=-\Big(E \pm \sqrt{E^2-\frac{4(i-1)}{N}}\Big).
\]
We therefore adopt the normalization
\begin{equation}
    \label{new scaling new}
    M_i=D_i/\prod_{j=1}^{i}\vert\rho_j^{+}\vert.
\end{equation}
The characteristic roots of the dynamic equation describing
$M_{i}$ can be approximated by the pair
$\rho^{\pm}_i/\left\vert\rho^{+}_i\right\vert$. 
% For $i>\theta^2_{N}N$, this is a
% complex-conjugated pair located on the unit circle in $\mathbb{C}$. In
% contrast,
For $i\leq E^2N/4$, these characteristic roots are
real, the first one equals minus one, and the other decreases from zero to minus one as $i$ goes from $1$ to $E^2N/4$.

Qualitatively, for most $i$, $M_{i}$ and $M_{i-1}$ have opposite signs and similar magnitudes, so that $R_i=M_{i}/M_{i-1}+1$
remains close to zero. However, for $i$ approaching $N$, $R_i$ starts
to develop more excited dynamics.

It may be of interest to note that $M_{i}/M_{i-1}$ can be interpreted
as normalized Sturm ratio sequence of matrix $\widehat{W}_N$. Sturm ratios
play a useful role in the analysis of large random matrices (see
e.g. \cite{albrecht2009sturm} or Section 1.9.3 in
\cite{forrester2010loggases}).  

In Section \ref{section away from edge}, we show that as long as the
local parameter, $\sigma_N$, of the singularity is slowly diverging to
infinity so that $\sigma_N \gg (\log\log N)^2$, the dynamics of $R_i$,
$i=1,...,N$, can be well approximated by a linear one. Then we use
this linear approximation  to obtain a CLT for the sums of the
logarithms of the normalized Sturm ratios. This leads to Theorem
\ref{th:main} with $\sigma_N \gg (\log\log N)^2$. In fact, for such
$\sigma_N$, our proof remains valid for matrices $\widehat{W}_N$ from
general Gaussian $\bm{\beta}$-ensembles (with $\bm{\beta}=1/\alpha \in
(0,\infty)$).  

To extend the theorem to slower growing and constant $\sigma_N$,
Section \ref{section up to edge} derives simple asymptotic formulae
for the Stieltjes transform of the empirical spectral distribution of
%$\hat{M}_N/\sqrt{N}$ 
$\widehat{W}_N$ and its derivative at the edge of the support
$[-2,2]$. These formulae and the Taylor expansion of the logarithm
describe the asymptotic behavior of the log statistics at the edge
with $\sigma_N\leq(\log\log N)^{3}$ in terms of that of the statistic
with $\sigma_N\gg(\log\log N)^{2}$. Thus, we obtain Theorem \ref{th:main} in its generality. 

Our proof can be easily extended to cover spiked GUE and GOE
matrices. See \cref{prop:zero-diag spiked} and the remark that follows that proposition.
\vspace{1mm}
%Theorem \ref{theorem spiked close} of Section \ref{section
%  spiked} generalizes Theorem \ref{th:main} for spiked GUE and GOE.}
 
  \textit{Extension to Wigner case}.
  Proving that a Wigner matrix $W'_N$ satisfies a certain property as
long as a matrix $W_N$ from G(O/U)E satisfies this property is often
based on the Lindeberg swapping process, where elements of $W_N$ are
replaced by the elements of $W'_N$ one by one without losing the
property. 
Typically, one needs to show that any individual swap does not change
the expectation $\mathbf{E}Q(M)$ of some smooth function $Q(\cdot)$ of
the matrix $M$ participating in the swapping process too much.

Although our initial interest is in the asymptotic normality of the log-determinant, we will
eventually need
to use Lindeberg swapping for several functionals which depend on the
Stieltjes transform evaluated at $z=E+\im \eta$ for $E$ near the edge and $\eta$ distant
at least $N^{-2/3-\delta}$ from the real axis -- here the gross
$N^{-2/3}$ scale is that appropriate for working at the edge of the
spectrum. 
We outline the swapping approach for the log-determinant example
but with the general class of ``Stieltjes edge functionals'' in mind. 

We adopt the method of \cite{tao2012central}, with modifications to work at the
edge, and under weakened assumptions, as described below.
We call a quantity $S(W_N)$ \textit{insensitive} at rate $\delta_N$
if $S(W_N) - S(W_N') = O(\delta_N)$.
Let $L_N(W_N) = \log |\det(W_N - E)|$. 
To extend the asymptotic normality of $L_N(W_N)$ to $L_N(W_N')$ it is
sufficient, via a standard smoothing argument,
to show that $\E G \circ L_N(W_N)$ is insensitive at rate $\delta_N$ for
scalar functions for which $\| G^{(j)} \|_\infty \leq b_N^j$.
For the log-determinant $\delta_N=b_N \asymp (\log N)^{-1/4}$ will work.
%\textcolor{blue}{where notation $a_N\asymp b_N$ means that $a_N\leq
%  Cb_N$ and $b_N\leq Ca_N$ for some $C$ and $N$ large.} 

In an initial \textit{approximation} step, we show that it suffices to
replace $L_N(W)$ by a function of the Stieltjes transform $s_W=N^{-1}\tr (W-z)^{-1}$
\begin{equation}  \label{eq:g-logdet}
  g(W) = N \int_{\gamma_N}^{N^{100}} \Im s_W(E+ \im \eta) \diff \eta.
\end{equation}
Here $\gamma_N = N^{-2/3-\delta}$: to show that values $0 \leq \eta \leq
\gamma_N$ can be neglected, we use an anti-concentration result 
that guarantees that with high probability, all eigenvalues are at least
$N^{-2/3 - \zeta}$--distant from $E$.
This too is proved by Lindeberg swapping, now with a second Stieltjes
functional, \cref{sec:wign-non-conc}.

     The swapping argument is now applied to show that $\E Q(W_N)$ is
insensitive for $Q$ of the form $(G\circ g)(W_N)$.
To review this in outline,
let $\gamma$ index an ordering of the independent components $\{\Re \xi_{ij}, \Im \xi_{ij}\}_{i<j}$ and $\{ \xi_{ii} \}$ of $W_N$.
Thus $\gamma$ runs over $N^2$ and $N(N+1)/2$ elements in the Hermitian and symmetric cases respectively.
By convention in each case, the first $N$ values of $\gamma$ index the
diagonal matrix entries. 
Thus $W^\gamma$ will refer to a matrix in which the elements prior to $\gamma$ come from $W_N'$ while those at $\gamma$ or later come from $W_N$.

At stage $\gamma$ in the swapping process, we can write
$W^{(0)} = W^{\gamma}$, $W^{(1)} = W^{\gamma +1}$, and
\begin{equation} \label{eq:swapdefs}
  W^{(0)} = W_0 + \frac{\xi^{(0)}}{\sqrt N} V, \qquad
  W^{(1)} = W_0 + \frac{\xi^{(1)}}{\sqrt N} V,
\end{equation}
and $W_0 = W^\gamma_{0}$ is independent of both $\xi^{(0)}$ and $\xi^{(1)}$.
In the symmetric case, $V$ is one of the elementary matrices of the form
$e_a e_a^*$ or $ e_a e_b^* + e_b e_a^*,$ 
for $1 \leq a < b \leq N$.
In the Hermitian case, we add matrices $ \im e_a e_b^* -  \im e_b e_a^*$.
The variables
$\xi^{(0)}$ and $\xi^{(1)}$ correspond to the $\gamma$th
components of $W_N$ and $W_N'$ respectively.
All matrices $W^\gamma$
% $W_{0}^\gamma$
are Wigner matrices.
%satisfying assumptions \textbf{W1-3}. \fix{not defined yet!}

To focus on individual swaps, write
\begin{equation*}
  \E Q(W) - \E Q(W')
    = \sum_\gamma \E \Delta_\gamma,
\end{equation*}
with $\Delta_\gamma = Q(W^\gamma) - Q(W^{\gamma +1})
 = Q(W^{(0)}) - Q(W^{(1)})$.
% \begin{equation*}
%   \Delta_\gamma = Q(W^\gamma) - Q(W^{\gamma +1})
%                 = Q(W^{(0)}) - Q(W^{(1)}).
% \end{equation*}

We consider $W^{(0)}$ and $W^{(1)}$ as perturbations of $W_0$.
Thus, set $W^\gamma_t = W_0^\gamma + t N^{-1/2} V_\gamma$, and 
% create notation for the resolvent and Stieltjes transforms that shows,
% when helpful, the dependence on $\gamma$. Set
% \begin{equation} \label{eq:t-defs0}
%   W^\gamma_t = W_0^\gamma + t N^{-1/2} V_\gamma, \qquad
%   R^\gamma_t = R^\gamma_t(z) = (W^\gamma_t - z)^{-1}, \qquad
%   s^\gamma_t(\eta) = N^{-1} \tr R^\gamma_t(E+\im\eta).
% \end{equation}
introduce $Q_\gamma(t) = Q(W_t^\gamma)$.
Note that this function is independent of $\xi^{(i)}$, and that
\begin{equation*}
  \Delta_\gamma = Q_\gamma(\xi^{(0)}) - Q_\gamma(\xi^{(1)}).
\end{equation*}
In a Taylor expansion of $Q_\gamma$, formal for now, this
independence implies
% of $\xi^{(i)}$ from $Q_\gamma(t)$ to write
\begin{equation*}
  \E [Q_\gamma(\xi^{(i)})]
  = \sum_j \frac{1}{j!} \E [ Q_\gamma^{(j)}(0) ] \,
                        \E\bigl([\xi^{(i)}]^j\bigr).
\end{equation*}

If moments match at order $j \leq k-1$,
that is, $\E([\xi^{(0)}]^j) = \E\left([\xi^{(1)}]^j\right)$,
then the $j$th order term in $\E \Delta_\gamma$ vanishes.
If, as one expects, $Q_\gamma^{(k)}(t)$ is of order $N^{-k/2}b_N$,
bounding the remainder term appropriately leads to the required bounds
on $\E\Delta_\gamma$. This is formalized in Proposition
\ref{prop:multi-matching}.

To show that such derivative bounds hold specifically for $Q = G \circ
g$ when $g$ 
is as in \eqref{eq:g-logdet}, we need good bounds for $\partial_t^j
g^\gamma(t)$ when 
$g^\gamma(t) = g(W^\gamma_t)$. Introduce notation for the resolvent
and Stieltjes transforms
\begin{equation} \label{eq:t-defs0}
%  W^\gamma_t = W_0^\gamma + t N^{-1/2} V_\gamma, \qquad
  R^\gamma_t = R^\gamma_t(z) = (W^\gamma_t - z)^{-1}, \qquad
  s^\gamma_t(\eta) = N^{-1} \tr R^\gamma_t(E+\im\eta).
\end{equation}
The standard resolvent perturbation argument  (equations \eqref{eq:neumann}-\eqref{coeff expl}) shows that
$\partial_t^j s^\gamma_t = c_j N^{-j/2 -1} \tr[(R^\gamma_t V)^j
R_t^\gamma]$.

This is bounded for $E$ near the edge and $\eta > N^{-2/3 - \delta}$
using the entrywise local law (see Proposition \ref{prop:P14analog}(i)).
Working at the edge allows, through use of the Ward identity, improvements in bounds because $\Im R$ is small. What results (see the proof of Proposition \ref{prop:system})
are bounds $\| \partial_t^j g^\gamma(t) \|_\infty \lesssim N^{-j/2} a_N$
with $a_N = 1$ in the log-determinant case. These bounds are
useful both for reducing the number of matching moments required to
three (for off-diagonal entries) and requiring only bounded variances
(for diagonal entries). Combining with the derivative bounds on $G$,
the chain rule shows that we obtain the desired insensitivity with
$\delta_N = a_N b_N=b_N$.

%The strategy of a) approximation by a Stieltjes edge functional,
%followed by b) Lindeberg swapping with resolvent perturbation is also
%used for other steps in our arguments: anticoncentration, and in
%JKOP2, for a fixed number of largest eigenvalues and for inverse
%moment functionals.

\subsection{Related work}
\label{sec:related-work}

% \subsection{Motivating questions}

% \label{sec:motivating-questions}

% \texttt{Make the case for AAP using both Baik-Lee transition and
% likelihood ratio as 'applied' motivations.}

%\newpage
%\subsection{(Previous) Overview of results and methods}
%\label{sec:overview-methods}

The
interest in determinants of Hermitian matrices from GUE, GOE and other
classical ensembles of Random Matrix Theory emerged in the 1960s from motivations in nuclear physics. The first published derivation %\footnote{%
%\textcolor{magenta}{In a footnote on p. 447, \cite{wigner1965distribution} points out that the joint
%distribution of the eigenvalues of GUE is \textquotedblleft probably known
%to many readers\textquotedblright .}} 
of the joint distribution of the
eigenvalues of GUE in \cite{wigner1965distribution} was spurred by the problem of
approximating the value of $\log \left\vert \det \left( W_{N}-E\right)
\right\vert ,$ where $E\in \left( -2,2\right) $.

 As pointed out by \cite{fyodorov2016fractional},
  a CLT for the GUE log statistic with singularity $E$ from a compact
  subset of $(-2,2)$ can be obtained from Theorem~1 of Krasovsky
  \cite{krasovsky2007correlations}. That theorem derives detailed
  asymptotics for the Laplace transform of the log statistic using
  Riemann-Hilbert machinery. Tao and Wu \cite{tao2012central} derive
  their CLT for $E=0$ and for Wigner matrices with atom distributions
  that match the first four moments of the normal. Bourgade and Mody
  \cite{bourgade2019gaussian} relax the conditions to require only
  matching of the first two moments.  Duy \cite{duy2017distributionsr}
  describes a very elegant proof of such a CLT for Gaussian ensembles
GUE/GOE and $E = 0$
  based on a representation of the corresponding log-determinants in
  the form of sums of independent random variables.  

For super-critical $E$ that lie outside an open set covering $[-2,2]$, the CLT for log-determinants follow from the CLT for more general linear spectral statistics with only super-critical singularities. Such a CLT is well known for classical ensembles (e.g. Johansson \cite{johansson1998fluctuations}). For extensions to Wigner matrices we refer the reader to Bai and Yao \cite{Yao2005}.

For the critical regime with $E$ local to $2$ the corresponding CLT has not been available. When this paper was close to completion, we learned about a recent work by Lambert and Paquette \cite{lambert2020a, lambert2020b} that obtains powerful asymptotic approximations to the logarithmic statistics with singularity local to the edge of the semi-circle law for Gaussian $\beta$-ensembles. Such approximations imply a CLT.

The analysis in \cite{lambert2020a, lambert2020b} starts from the recurrence for the minors of $z-\widehat{W}_N/2$, equivalent to our \eqref{determinant_recurrence_nonzero} with $z$ interpreted as $E/2$. The deterministic version of their recurrence generates monic Hermite polynomials $\pi_i(z)$ orthogonal with respect to the weight $\exp (-2Nz^2)$,
\[
\begin{pmatrix}
\pi_i(z)\\\pi_{i-1}(z)
\end{pmatrix}
=\tilde{T}_i(z)
\begin{pmatrix}
\pi_{i-1}(z)\\\pi_{i-2}(z)
\end{pmatrix}
\qquad \text{with}
\qquad
\tilde{T}_i(z)=
\begin{pmatrix}
z & -\frac{i-1}{4N}\\1 & 0
\end{pmatrix}.
\]

Lambert and Paquette point out three regimes of this recurrence: \textit{hyperbolic, parabolic}, and \textit{elliptic}, corresponding to the eigenvalues of $\tilde{T}_i(z)$ being, respectively, distinct real, coinciding or local to each other, and distinct complex. The regimes are associated with growing, Airy-type transitory, and oscillatory behavior of the Hermite polynomials, respectively. In terms of $E=2z$, the recursion remains in hyperbolic or elliptic regimes for all $i=1,...,N$ as long as $|E|>2+\varepsilon$ or $|E|<2-\varepsilon$, respectively. It enters the parabolic regime for relatively large $i$ if $|E-2|=O(N^{-2/3})$.

%\texttt{INTRO TO BE CONTINUED FROM THIS POINT. Build on Introduction of Gaussian ms and first section of combined2.tex }

\cite{lambert2020a} studies the hyperbolic regime of the recurrence for the minors of $z-\widehat{W}_N/2$.  It covers the range $E>2+\sigma_N N^{-2/3}$, where $\sigma_N \gtrsim \log^{2/3}N$ in our notations. Instead of analysing the dynamics of the Sturm ratios $M_i/M_{i-1}$ as we do, \cite{lambert2020a} base their analysis on an approximation to the product $T_N(z)...T_2(z)$, where $T_i(z)$ are the stochastic analogues of the deterministic transfer matrices $\tilde{T}_i$.

\cite{lambert2020b} extends \cite{lambert2020a} to the parabolic regime by noting that scaled versions of the minors of $z-\widehat{W}_N/2$ satisfy a finite difference equation which can be interpreted as a discretisation of the stochastic Airy equation. This yields a refined asymptotic approximation for the log-determinant.

Our CLT for the parabolic regime does not rely on the stochastic Airy equation machinery. Instead, we use asymptotics of 1-point correlation function for GUE to link hyperbolic and parabolic regimes. Although the resulting asymptotic approximations are less refined than those obtained in \cite{lambert2020b}, they do deliver the CLT for the log-determinant. In contrast to \cite{lambert2020a, lambert2020b}, we extend our asymptotic results to the general Wigner setting.

Another related and independently written paper is Augeri, Butez, and Zeitouni \cite{augeri2020}. It deals with the CLT for $\beta$-ensembles when the singularity $E$ is a fixed number in the bulk $(-2,2)$. Such a location of the singularity implies that, as $i$ goes from $1$ to $N$, the recurrence for the minors of $W_N-E$ goes through all the regimes, starting from the hyperbolic, transiting through the parabolic, and finishing in the elliptic regime. \cite{augeri2020} refers to these regimes as ``scalar'', ``transitory'', and ``oscillatory''. The analysis in the ``scalar'' regime is similar to ours. However, that of the ``transition'' regime is based on combinatoric arguments, whereas ours is using 1-point correlation and the asymptotics of the Stieltjes transform.

In contrast to \cite{augeri2020}, we do not analyze ``oscillatory'' regime because we focus on the edge singularity $E$. Although \cite{augeri2020} only consider $E$ inside $(-2,2)$, we believe that their analysis can be extended to the edge with some extra work. Unlike \cite{augeri2020}, we do extend our results to Wigner matrices.

Finally, an interesting recent paper by Bourgade, Mody, and Pain \cite{bourgade2021} obtains a CLT for the real and imaginary parts of the log determinant of $\beta$-ensemble when the singularity is in the bulk. The proof is based on a new local law result, and is completely different from the proof used in our paper. We do not know whether the proof of \cite{bourgade2021} can be extended to the case of the singularity at the edge.

\subsection{Organization and some notation} The rest of the paper is organized as follows. Section \ref{section away from edge} states a CLT for the log determinant for the special case of slowly diverging $\sigma_N$, and outlines its proof. Section \ref{section up to edge} extends this result to constant (possibly negative) $\sigma_N$, which yields Theorem~\ref{th:main}.  Section \ref{sec:proofs-gauss-ensembl} implements the strategy of the proof outlined in \cref{section away from edge} and establishes key bounds on the Stieltjes transform that are postulated in \cref{section up to edge}. Sections \ref{sec:modif-intr} and \ref{section spiked} analyze Wigner matrices and the spiked case, respectively. Relatively more technical proofs are compiled in the Supplementary Material.

Notations $a_N\ll b_N$, $a_N\lesssim b_N$, and $a_N\asymp b_N$  mean, respectively, that $a_N/b_N\rightarrow 0$, that $a_N\leq  Cb_N$ for some $C$ and $N$ large, and that $a_N\lesssim b_N$ and $b_N\lesssim a_N$.

\section{G$\beta$E: The CLT slightly away from the edge}
\label{section away from edge}

In this section we establish the following analogue of
Theorem~\ref{th:main} for general Gaussian  $\beta$-ensembles in
cases where \( (\log\log N)^{2} \ll \sigma_N \ll \log^2 N \). Hence,
the location of singularity $E=E_N$ is slightly away from the edge
in the sense that $(E-2) N^{2/3}$ slowly diverges to infinity.

%\fixb{Introduce centering and scaling constants}
%\begin{equation}
%  \label{eq:centscal}
%  \fixb{ \begin{split}
%    \mu_N & = \tfrac{1}{2} N + \sigma_N N^{1/3} - \tfrac{2}{3}
%    \sigma_N^{3/2} - \tfrac{1}{6}(\alpha-1) \log N \\
%    \tau_N & = \sqrt{\tfrac{\alpha}{3} \log N}, \qquad
%    \tilde{\tau}_N
%     = \sqrt{ \alpha \rho(\theta_N^{-2})},
%     \qquad \rho(x) = \log \tfrac{1}{2}[1+(1-x)^{-1/2}].
%  \end{split}}
%\end{equation}

\begin{theorem}
\label{theorem main}Consider matrix $\widehat{W}_N$ from a (scaled) general
Gaussian $\beta$-ensemble \eqref{tridiag form earlier} with
$\beta=2/\alpha$.
Let $D_N=\det(\widehat{W}_N-2\theta_N)$, where $2\theta_N\equiv E=2+N^{-2/3}\sigma_{N}$ with $\left(
\log \log N\right) ^{2}\ll \sigma_{N}\ll \left( \log N\right) ^{2}.$
Then,  
\begin{equation*}
  ( \log |D_N| - \mu_N)/\tilde{\tau}_N
   \overset{d}{\rightarrow} \mathcal{N}\left( 0,1\right) ,
 \end{equation*}
 where
 \[
     \tilde{\tau}_N
     = \sqrt{ \alpha \rho(\theta_N^{-2})}
     \quad \text{with} \quad \rho(x) = \log \tfrac{1}{2}[1+(1-x)^{-1/2}].
 \]
% \[
% \frac{\log \left\vert \mathcal{D}_{N}\right\vert -N/2  +\frac{\alpha-1}{6}\log N
% -\sigma_N N^{1/3}+\frac{2}{3}\sigma_N^{3/2}}{\sqrt{\alpha\log \frac{%
% \theta_{N}+\sqrt{\theta_{N}^{2}-1}}{2\sqrt{\theta_{N}^{2}-1}}}}%
% \overset{d}{\rightarrow} \mathcal{N}\left( 0,1\right) .
% \]
\end{theorem}
\begin{remark}
The modified scaling in the above CLT naturally arises from the
arguments in our proof. Note that
\begin{equation*}
 \rho(\theta_N^{-2}) = \tfrac{1}{3} \log N - \tfrac{1}{2} \log
  \sigma_N - \log 2 + O(N^{-1/3}\sigma_N^{1/2}),
\end{equation*}
% \[
% \log \frac{%
% \theta_{N}+\sqrt{\theta_{N}^{2}-1}}{2\sqrt{\theta_{N}^{2}-1}}=\frac{1}{3}\log N-\frac{1}{2}\log \sigma_N-\log 2+O(N^{-1/3}\sigma_N^{1/2}),
% \]
so that the asymptotic variance of $\log|D_N|$ is $\frac{\alpha}{3}\log N$, as in Theorem~\ref{th:main}. However, our Monte Carlo experiments (which we do not report here) indicate that the scaling in Theorem~\ref{theorem main} makes the standard normal approximation better in finite samples.
\end{remark}

The proof of Theorem \ref{theorem main} is based entirely on the recurrence equation \eqref{determinant_recurrence_nonzero} with application of some well-known deviation and concentration inequalities for sums of independent random variables. In this section, we briefly outline the main steps of the proof. Details can be found in \cref{Appendix A}.

Define normalized versions of the characteristic roots
  $\rho_j^\pm$: 
\begin{equation}
  \label{eq:rmdef}
r_{i}=1+\sqrt{1-\frac{i-1}{N\theta _{N}^{2}}}, \qquad m_{i}=1-\sqrt{1-%
\frac{i-1}{N\theta _{N}^{2}}}. 
\end{equation}
In particular, $|\rho_i^+| = \theta_N  r_i$.
Then, from \eqref{determinant_recurrence_nonzero} and the identities $r_i+m_i=2$ and $r_im_i = (i-1)/N\theta_N^2$, the normalized determinants \eqref{new scaling new} follow the recurrence%
\begin{equation}
\label{normalized recurrence}
M_{i} =\left( \alpha_i-\gamma_i-1\right)
M_{i-1}-\left(\gamma_i+ \beta_i-\delta_i\right) M_{i-2},
\end{equation}
where 
\[
\alpha_{i} =\frac{a_{i}}{\sqrt{N}\theta _{N}r_{i}},\qquad \gamma_{i} =\frac{m_{i}}{r_{i}},\qquad \beta _{i}= \frac{\sqrt{\gamma_i}c_{i-1}}{\sqrt{N}\theta_{N}r_{i-1}},\qquad
 \delta_{i}=\frac{m_{i}}{r_{i}}-\frac{m_{i}}{r_{i-1}}.
\]
With the conventions $M_0 = 1, M_{-1} =0$, and
 declaring $c_0 = 0$, so that
  $\beta_1 = 0$ along with $\gamma_1 = \delta_1 = 0$, equation \eqref{normalized recurrence} holds for $i=1, \ldots, N$.
  
Dividing both sides of \eqref{normalized recurrence} by $M_{i-1}$
% , and define $R_i=E_i/E_{i-1}+1$
yields
% \begin{align}
%     R_{i} := \frac{E_i}{E_{i-1}} + 1 
%   & =\alpha _{i}-\gamma _{i}+\frac{\gamma _{i}+\beta _{i}-\delta _{i}}{%
% 1-R_{i-1}}.  \label{R equation compact} \\
%   & = \xi_{i}+\gamma _{i}R_{i-1}+\varepsilon _{i}.
% \end{align}
\begin{equation}
  R_{i} \equiv \frac{M_i}{M_{i-1}} + 1 
  =\alpha _{i}-\gamma _{i}+\frac{\gamma _{i}+\beta _{i}-\delta _{i}}{%
1-R_{i-1}},  \label{R equation compact}
\end{equation}
which can be rewritten as a recurrence \begin{equation}
\label{linear looking}
R_{i}  =\xi_{i}+\gamma _{i}R_{i-1}+\varepsilon _{i},
\end{equation}
for $i=1,\ldots, N$, with the definitions
\begin{align}
\xi_{i} &= \alpha_{i} + \beta_{i}, \notag \\
\varepsilon _{i} &=-\delta _{i}+\left( \beta _{i}-\delta _{i}\right) \frac{%
                   R_{i-1}}{1-R_{i-1}}+\gamma _{i}\frac{R_{i-1}^{2}}{1-R_{i-1}}.
                   \label{eq:eps-def}
\end{align}
%Note that $\xi_1 = \alpha_1$ and $\varepsilon_1 = 0$.
By dropping the \textit{non-linear} term \( \varepsilon_i\) from
\eqref{linear looking}, we now define a \textit{linear} process
\(\{L_i\}_{i = 1}^{N} \)
satisfying the recursion 
\begin{equation}
\label{eq: L recursion}
L_{i} = \xi_{i}+\gamma _{i}L_{i-1}.
\end{equation}
In particular, $L_1 = \xi_1 = \alpha_1$.
Note that $\{\xi_i \}$ are independent random variables, while $\{ \gamma_i \}$ are deterministic. 
%Iterating this yields
%\begin{equation}\label{L as sum}
%  L_{i} = \; (T \xi)_i
%   = \xi _{i}+\gamma _{i}\xi _{i-1}+ \dots +
%  \gamma _{i} \dots \gamma_{2}\xi _{1}, \qquad \qquad i \geq 2,
%\end{equation}%
%and a similar iteration of the recursion for \(R_i\) yields
%$R_1 = L_1$ and
%\begin{equation}\label{R i expanded}
%  R_{i} = L_i+(T \varepsilon)_i
%  = L_{i} + \varepsilon_{i} + \gamma_{i} \varepsilon_{i-1} + \dots +
%  \gamma_{i} \dots \gamma_{3} \varepsilon_{2},
%    \qquad i\geq 2.
%  + \gamma_{i} \dots \gamma_{3} R_{2} \qquad i \geq 4.
%\end{equation}

To establish the CLT, we study  the dynamics of $ L_i$ and $R_{i}$. Our proof consists of the following three steps: 

\begin{enumerate}
\item First, derive a CLT for \( \sum_{j = 1}^{N} L_j \) with the variance of exact order \(\log N\).

\item Then, show in the regime \( \sigma_N \gg (\log\log N)^{2} \) that both  \( \max_{i} |L_i| = o_{\Pr}(N^{-1/3}) \) and \( \max_{i} |R_i| = o_{\Pr}(N^{-1/3}) \). This allows us to use Taylor's approximation for the logarithm, so that 
\[
  \log |M_N| =  \sum_{j = 1}^{N} \log |1 - R_j|
%  + \log |E_2|
  = \sum_{j = 1}^{N} (-R_j - R_j^2 / 2) + o_{\Pr}(1) \, .
\]
%We point out that all values \( R_i \), $i=1,...,N$, remain strictly smaller than \(1\) as long as the eigenvalues of the matrix \( M_N/\sqrt{N} \) remain smaller than \( 2\theta_N\), as follows from e.~g.~\cite[Lemma 2.1]{albrecht2009sturm}. In this sense the condition  \( \sigma_{N} \gg (\log \log N)^2 \) seems a bit excessive. On the other hand, we show that in this regime, the process \(R_i\) closely follows the linear process \(L_i\), which is obviously a lot easier to control.

\item Finally, prove that the sum $\sum\nolimits_{j=1}^{N} \left( -R_{j}-R_{j}^{2}/2\right) $ can be replaced with $\sum\nolimits_{j=1}^{N} -L_{j}$ at the cost of some $O_{{\Pr}}(\log\log N)$ error term and with some explicit deterministic shift.
\end{enumerate}

Achieving these objectives will show that the asymptotic behavior of $\log
\left\vert M_{N}\right\vert $ is the same as that of $-\sum%
\nolimits_{i=1}^{N}L_{i}$ up to $O_{\Pr}\left( \log\log N\right)$ and an explicit deterministic shift, and hence an appropriately centered $\log \left\vert M_{N}\right\vert $ satisfies
the same CLT as $-\sum\nolimits_{i=1}^{N}L_{i}.$
After calculating the deterministic shift between \( \log |D_N| \) and \( \log |M_N| \), we derive the CLT for the log-determinant as required. %Details of these derivations are given in \cref{Appendix A}.

\section{G(U/O)E: All the way to the edge}%
\label{section up to edge}

Theorem \ref{theorem main} covers singularities
$2 \theta_N = E = 2 + N^{-2/3} \sigma_N$ in the range
$(\log\log N)^2 \ll \sigma_N \ll \log^2 N$ for all positive $\alpha$.
We seek to extend the result to singularities at a distance of exact
order $N^{-2/3}$ away from the edge, or even (just) inside the bulk.
We consider now sequences $\sigma_N$ satisfying
\begin{equation}
  \label{eq:sigcon}
  -\gamma \leq \sigma_N \leq \bar{\sigma}_N:= (\log \log N)^3
  \qquad \qquad \text{for some} \ \gamma > 0.
\end{equation}

Our extension will rely on the
properties of GUE and GOE, and so covers only the cases $\alpha=1$ and
$\alpha=2$. 
Indeed, the main tool is uniform
approximation of the one-point function of GUE for regions
up to and containing the spectral edge, based chiefly on results of
G\"{o}etze and Tikhomirov \cite{Gotze2005}.

In Section~\ref{appendix up to edge}, we use the one-point function
approximation to obtain the following estimates on the Stieltjes
transform and its derivative near the edge. 

\begin{proposition}
\label{lemma stieltjes}Suppose that $\alpha=1$ or $\alpha=2$ and let
$\sigma_N = N^{2/3}(E-2)$ satisfy condition \eqref{eq:sigcon}. Then
\[
\sum\nolimits_{i=1}^{N}(E-\lambda_i)^{-1} - N =
O_{\Pr}\left((1+|\sigma_N|^{1/2})N^{2/3}\right) 
\]
and
\[
\sum\nolimits_{i=1}^{N}(E-\lambda_i)^{-2}=O_{\Pr}(N^{4/3}).
\]
\end{proposition}

\begin{proof}[Proof of Theorem \ref{th:main} in the range \eqref{eq:sigcon}]
  Let
  \begin{equation*}
    S_N(\sigma_N)  :=\sum_{i=1}^N \log |2+N^{-2/3}\sigma_N-\lambda_i| - \mu_N(\sigma_N)=\sum_{i=1}^N \log |E-\lambda_i| - \mu_N(\sigma_N),
  \end{equation*}
  where $\mu_N(\sigma_N)$ is given by \eqref{eq:centscal1}.
  The strategy is to use Proposition \ref{lemma stieltjes} to show that
\begin{equation}
  \label{eq:strategy}
  S_N(\bar{\sigma}_N) - S_N(\sigma_N) = O_\Pr\big((\log \log N)^6\big),
\end{equation}
so that $S_N(\sigma_N)/\sqrt{\frac{\alpha}{3} \log N}$ has the same
limiting $\mathcal{N}(0,1)$ distribution as
$S_N(\bar{\sigma}_N)/\sqrt{\frac{\alpha}{3} \log N}$, the latter being
given by Theorem \ref{theorem main}.

Abbreviate $\ell_N = (\log \log N)^6$ and note that $|\sigma_N|^{3/2}
\leq \bar{\sigma}_N^{3/2} = o(\ell_N)$.
Introducing 
\begin{align}
   \delta_N & = N^{-2/3}(\bar{\sigma}_N - \sigma_N), \notag \\
    d_i   & = \log|E-\lambda_i + \delta_N| - \log |E-\lambda_i| - \delta_N
            (E-\lambda_i)^{-1}, \notag \\
\intertext{we can decompose}
  S_N(\bar{\sigma}_N) - S_N(\sigma_N)
  & = \sum_{i=1}^N (\log|E - \lambda_i + \delta_N|
    -\log|E - \lambda_i|)
    -(\bar{\sigma}_N - \sigma_N)N^{1/3} + o(\ell_N) \notag \\
  & = \sum_{i=1}^N d_i + \delta_N \Big[\sum_{i=1}^N (E-\lambda_i)^{-1} -
    N\Big] + o(\ell_N).  \label{eq:negone} 
%    \mu_i & = 2 + N^{-2/3} \sigma_N - \lambda_i, \qquad
\end{align}

We will show that for each $\varepsilon > 0$,
with probability at least $1 - \varepsilon$,
\begin{equation}
  \label{eq:dbnd}
  \bigg|\sum_1^N d_i \bigg| \leq \delta_N^2 \sum_1^N
  (E-\lambda_i)^{-2} + o(\ell_N). 
\end{equation}

The bound \eqref{eq:strategy} then follows directly from Proposition
\ref{lemma stieltjes}, since $N^{2/3} \delta_N = O(\ell_N^{1/2})$.
Thus it remains to establish \eqref{eq:dbnd}.
For this we use some consequences of convergence to the Tracy-Widom
law formulated in the following lemma, proved in \cref{sec:proof-lemma-4}.
% also proved in Section \ref{appendix up to edge}.

\begin{lemma}\label{edge spacing lemma}
  Suppose that $W_N$ is (scaled) GUE/GOE.
Let $\gamma > 0$ be fixed, and suppose that $E = 2 + \sigma_N
N^{-2/3}$ with $\sigma_N > -\gamma$,
and that $\bar{E} = 2 + N^{-2/3} \bar{\sigma}_N$.
Then for each $\epsilon>0$ small,
there exists $k = k(\epsilon,\gamma)$ such that for large $N$,
  % Then, for any  $\gamma > 0$ , there exists \( k = k(\epsilon, \gamma) \) such that whenever \( \sigma_N > -\gamma \) and for large enough \( N \),
\begin{equation}
  \label{eq:first-pair}
  \Pr( \lambda_1 > \bar{E} - N^{-2/3} ) < \epsilon, \qquad
    \Pr( \lambda_k > E ) < \epsilon \, .
\end{equation}
% \[
%     \Pr( \lambda_k > E ) < \epsilon \, .
% \]
Moreover, there are constants \(c_1 = c_1(\epsilon, \gamma) \) small and \( C_1 = C_1(\epsilon, \gamma)\) large, such that for large enough \(N\),
\begin{equation}
  \label{eq:second-pair}
    \Pr\Big( \min_{i \leq N} |\lambda_i-E| < c_1N^{-2/3} \Big) < \epsilon, 
    \quad
    \Pr\Big( \max_{i \leq k} |\lambda_i-E| > (C_1 +
      |\sigma_{N}|)N^{-2/3} \Big) < \epsilon\, .   
\end{equation}
\end{lemma}

Turning to \eqref{eq:dbnd}, our first goal is to establish
probabilistic bounds on $\vert d_i \vert$ for $i=1,...,N$.
Let $\mu_i = E-\lambda_i$.
Given $\varepsilon$ and $\gamma$, Lemma \ref{edge spacing lemma} yields
$k, c_1, C_1$ such that the event
$\mathcal{E} = \mathcal{E}(k, c_1, C_1)$ given by
\[
\mathcal{E} = 
%\mathcal{E}(k, c_1, C_1)=
\Big\{ \lambda_1 \leq \bar{E} - N^{-2/3}, \ 
\mu _{k}>0,\;\;\;
  \min_{i=1,...,N}\vert N^{2/3}\mu _{i}\vert \geq
  c_{1},\;\;\;\max_{i=1,...,k}\vert N^{2/3}\mu 
_{i}\vert \leq C_{1}+|\sigma_N|\Big\} 
\]%
has probabality at least $1 - \varepsilon$. On event $\mathcal{E}$, 
for $i\geq k,$ the bound $|\log(1+x)-x| \leq x^2/2$ for $x \geq 0$
implies that
\[
  |d_{i}| \leq  \tfrac{1}{2} \delta_N^2 \mu_i^{-2}.
  % \frac{N^{-4/3}\left( \sigma^{*}_N
  %           -\tilde{\sigma}_N\right) ^{2}}{%2\mu _{i}^{2}}.
\]%
% This follows from the
% Taylor remainder bound $|\log(1+x)-x| \leq x^2/2$ for $x \geq 0$.
Further, note that 
\begin{eqnarray*}
\left\vert d_{i}\right\vert  &=&\big| \log |
              \bar{\sigma}_N-\sigma_N+N^{2/3}\mu_{i}|
                                 -\log | N^{2/3}\mu _{i}|
          - \delta_N \mu_i^{-1} \big|. %\\
%&\leq &\log \big( 2\left( \log \log N\right) ^{3}\big) +\log (C_{1}
%+|\sigma_N|)+ c_1^{-1}\left( \log \log N\right) ^{3}.
\end{eqnarray*}%
On the other hand, still on $\mathcal{E}$, for any $i \leq k$ and all
  sufficiently large $N$, the first logarithm
is non-negative and
no larger
than $\log(3\bar{\sigma}_N+C_1)$;
the second one is no larger in absolute value than
  $|\log c_1|+\log(\bar{\sigma}_N+C_1)$; and the last term on the
  right hand side of the above display is no larger in absolute value
  than $2\bar{\sigma}_N/c_1$. Each of these bounds is
  $O(\ell_N^{1/2})$. 

Hence, overall on $\mathcal{E}$,%
\begin{equation}
\left\vert \sum\nolimits_{i=1}^{N}d_{i}\right\vert \leq
% \sum\nolimits_{i=1}^{N}\left\vert d_{i}\right\vert \leq
\delta_N^2 \sum\nolimits_{i=1}^{N} \mu_i^{-2} + C_2 \ell_N^{1/2}
\label{d_bound}
\end{equation}%
for some constant $C_2 = C_2(\varepsilon,\gamma)$.
Since $\Pr \left( \mathcal{E}\right) \geq 1-\varepsilon ,$ the latter
inequality holds for sufficiently large $N,$ with probability at least $%
1-\varepsilon $.

This yields \eqref{eq:dbnd} and completes the proof.
\end{proof}

\section{Proofs for Gaussian ensembles}
\label{sec:proofs-gauss-ensembl}

  \subsection{Proofs from Section \ref{section away from edge}}
  \label{Appendix A}

  This section implements the three steps of the analysis, described in
Section \ref{section away from edge}, that lead to Theorem \ref{theorem main}.
%After some preliminaries, we start from the most
%straightforward step (step 3), which 
%calls for deriving a CLT for the sum of the linear process $L_i$.

\subsubsection{Preliminaries}
%The notation $a_N = O(b_N)$ or $a_N \lesssim b_N$ means that for some
%$C >0$ and large $N > N(C)$, we have $a_N \leq C b_N$.
%In such bounds the value of $C$ may change from one appearance to
%another.
For \(p \geq 1\), denote by \( \| X \|_{p} = (\E |X|^{p})^{1/p} \) the
\(p\)-norm of a random variable \( X \).

\smallskip
\textit{Sub-gaussianity.} We say that a centred random variable $X$ belongs to the
\emph{sub-gamma} family
${SG}(v,u)$ for $v, u > 0$, if
\begin{equation}\label{definitions SG}
    \log \E e^{tX} \leq \frac{t^2 v}{2(1 - tu)},
    \qquad
    \forall t: \; |t| < \frac{1}{u} \, .
  \end{equation}
If $X \in SG(v,u)$ then $X \in SG(v',u')$ for each $v' \geq v$,$u' \geq u$.
 For arbitrary $c \in \R$, we have $cX\in SG\left(
c^{2}v_{X},|c| u_{X}\right) .$  
If $X\in SG\left( v_{X},u_{X}\right) $ and $Y\in SG\left( v_{Y},u_{Y}\right) 
$ are independent, then $X+Y\in SG\left( v_{X}+v_{Y},\max \left\{
u_{X},u_{Y}\right\} \right).$ If $X\sim \mathcal{N}\left(
0,1\right) ,$ then $X\in SG\left( 1,0\right) $, and if $X\sim \chi ^{2}\left( d\right) -d,$
then $X\in SG\left( 2d,2\right) .$
We refer to chapter~{2.4} of \cite{boucheron2013concentration} for this
last result.

\smallskip
\textit{The recurrence parameters.} 
We often view the deterministic sequences $r_i, m_i,
\gamma_i, \delta_i,  g_i$ etc.  as discretizations of functions
evaluated at $x_i = (i-1)/(N \theta_N^2)$, with step size
$\Delta_N = 1/(N \theta_N^2)$.
%Note that $x_1 = 0$ and that $x_N = 1-2w_N N^{-2/3} +
%O(N^{-1})$. 
For example, $r_i = r(x_i)$ with $r(x) = 1 + \sqrt{1-x}$
for $x \in [0,1)$, where this function is concave decreasing.

\smallskip
\textit{The operator $T$.} \
Iterating \eqref{eq: L recursion} yields
\begin{align}\label{L as sum}
  L_{i} &= \; (T \xi)_i
   = \xi _{i}+\gamma _{i}\xi _{i-1}+ \dots +
  \gamma _{i} \dots \gamma_{2}\xi _{1}, \qquad \qquad i \geq 2,\\
  L_{1}&=(T\xi)_1=\xi_{1}.\notag
\end{align}%
The generalized exponential moving average used in \eqref{L as sum} is a linear map $T = T_N: \R^N \to \R^N$
formally defined by
\begin{equation*}
  Ta_i = (Ta)_i =
  \begin{cases}
    a_i + \gamma_i a_{i-1} + \cdots + \gamma_i \cdots \gamma_2 a_1  & i \geq 2 \\
    a_1 & i = 1.
  \end{cases}
\end{equation*}  
The corresponding matrix $T = (T_{ij})$ is lower triangular, with
entries
\begin{equation}
  \label{eq:gami-to-j}
  T_{ij} = \gamma_{i:j+1}, \qquad
  \gamma_{i:j} =
  \begin{cases}
    \gamma_i \gamma_{i-1} \cdots \gamma_j  &  i \geq j \\
    1          & i = j-1  \\
    0          & i < j-1.
  \end{cases}
\end{equation}
Resummation yields
\begin{equation}
  \label{eq:resummation}
  \sum_{i=1}^N Ta_i = \sum_{j=1}^N  g_{j+1} a_j,
\end{equation}
where $g_{j} = (T^* \mathbf{1})_{j-1}$ is given by
\begin{equation}
  \label{eq:gjdef}
  g_j = 1 + \gamma_j + \ldots + \gamma_N \cdots \gamma_j
      = 1 + \sum_{i=j}^N \gamma_{i:j}
\end{equation}
for $2 \leq j \leq N$ and $g_{N+1} = 1$.
Since $\gamma_i$ is increasing, we have
\begin{equation}
  \label{eq:T-bound}
  |Ta_i| \leq \frac{1}{1-\gamma_i} \max_{j \leq i} |a_j|,
\end{equation}
which simplifies to $Ta_i \leq (1-\gamma_i)^{-1} a_i$ for an
increasing sequence $a_i \geq 0$.
% Note for later that 
If $T_i$
satisfies a recurrence
\begin{equation}
  \label{eq:recurr}
    T_i = \gamma_i T_{i-1} + a_i, \qquad \qquad T_1 = a_1
\end{equation}
then $T_i = (Ta)_i$ and bound \eqref{eq:T-bound} applies.

\smallskip
\textit{Estimates by integrals.} \
With
$\Delta_N = 1/(N \theta_N^2)$, if $f(x)$ is increasing, then
\begin{equation}
  \label{eq:riemann}
  \sum_{i=a}^b f(x_i) \Delta_N \leq \int_{x_a}^{x_b+\Delta_N} f(x)
  \diff x.
\end{equation}
In the corresponding lower bound the integral has limits $x_a-\Delta_N$ and
$x_b$.
In particular,
\begin{align}
  \frac{1}{N} \sum_{i=1}^N \frac{1}{(r_i-1)^\beta}
    & < \theta_N^2 \int_0^{x_N+\Delta_N} \frac{\diff
      x}{(r(x)-1)^\beta}
      = \theta_N^2 \int_0^{\theta_N^{-2}} \frac{\diff
      x}{(1-x)^{\beta/2}}  \notag
  \\
    & \leq
      \begin{cases}
        C_\beta \theta_N^2   &  \text{if } \beta \in (0,2) \\
        C_\beta w_N^{1-\beta/2} N^{(\beta-2)/3}  & \text{if } \beta > 2,
      \end{cases}
  \label{eq:ri-sums}
\end{align}
for sufficiently large $N$ with $C_\beta = 2|\beta-2|^{-1}$ and $w_N=\sigma_N/2$. 
Alternatively, 
the error bound for the trapezoid rule, in integral form, states
\begin{equation}
  \label{eq:trapezoidal-bound}
  \begin{split}
  \Bigg| \sum_{i=a}^b f(x_i) \Delta_N - \int_{x_a}^{x_b} f(x) \diff x \Bigg|
  & \leq \frac{\Delta_N^2}{8} \int_{x_a}^{x_b} |f''(x)| \diff x
  + \frac{\Delta_N}{2} [|f(x_a)|+|f(x_b)|]  \\
  & = \varepsilon_{N1}(f) + \varepsilon_{N2}(f).    
  \end{split}
\end{equation}

\subsubsection{Step one: CLT for  $\sum_{i = 1}^{N}L_{i}$ }\label{CLT for L section}
Equation \eqref{eq:resummation} applied to $L_i = T \xi_i$ yields
\begin{equation*}
    \sum_{i=1}^{N} L_{i}=\sum_{i=1}^{N}g_{i+1}\xi _{i},
\end{equation*}
where $g_{i}$ is as in \eqref{eq:gjdef},
and we recall that 
$\xi_{i}$ are independent zero-mean random variables.
Lyapunov's CLT implies that 
\[
\frac{\sum_{i=1}^{N}L_{i}}{\sqrt{\sum_{i=1}^{N}g_{i+1}^{2}\E\xi
_{i}^{2}}}\overset{d}{\rightarrow} N\left( 0,1\right) 
\]%
as long as 
\begin{equation}
\sum_{i=1}^{N}g_{i+1}^{4}\E\xi _{i}^{4}/\left(
\sum_{i=1}^{N}g_{i+1}^{2}\E\xi _{i}^{2}\right) ^{2}\rightarrow 0.
\label{Lyapunov}
\end{equation}
Let us establish (\ref{Lyapunov}), and find $\sum_{i=1}^{N}g_{i+1}^{2}%
\E\xi _{i}^{2}$.

The proofs of the following two lemmas are given in Subsections \ref{proof of lemma moments} and \ref{proof of lemma si}.

\begin{lemma}
\label{lemma moments of ksi}For any integer $q\geq 1$ and all
sufficiently large $N$, there exist 
constants $c_2, C_{q}>0,$ such that for all $1\leq i\leq N,$ 
\[
\E\xi _{i}^{2q}\leq C_{q}\alpha^{q}N^{-q}\text{ and }\E\xi
_{i}^{2}> c_2 \alpha N^{-1}
\]
\end{lemma}

\begin{lemma}
\label{lemma si}Let $w_N=\sigma_{N}/2$ so that $\theta _{N}=1+N^{-2/3}w_{N}$. Suppose that $\left( \log
  \log N\right) ^{2}\ll w_{N}\ll \left( \log N\right) ^{2}.$
Then, for all $1\leq
i\leq N-N^{1/3},$ any $k>0,$ and all sufficiently large $N,$%
\begin{equation}
g_{i}>\frac{r_{i}}{2\left( r_{i}-1\right) }\left( 1-\log ^{-k}N\right) .
\label{lower bound on si}
\end{equation}%
Further, for all $1\leq i\leq N$ and all sufficiently large $N$,%
\begin{equation}
g_{i}<\frac{r_{i}}{2\left( r_{i}-1\right) }\left( 1+w_{N}^{-3/2}\right) .
\label{upper bound on si}
\end{equation}
There is also a trivial bound, sharper than \eqref{upper bound on si} for $ N - N^{1/3} \leq i \leq N+1$:
\begin{equation}
  \label{eq:g-trivial}
  g_i \leq N-i+2.
\end{equation}
\end{lemma}

Let $n_0 = N - [N^{1/3}]-1$. These lemmas yield, along with Riemann
sum bounds like \eqref{eq:riemann},~\eqref{eq:ri-sums} 
\begin{equation*}
  \sum_{i=1}^{N}g_{i+1}^{2}\E\xi _{i}^{2}
  \gtrsim \frac{1}{N\theta_N^2} \sum_{i=2}^{n_0} \frac{1}{(r_i-1)^2}
  > \int_{x_1}^{x_{n_0}} \frac{\diff x}{1-x}
  = -\log(1-x_{n_0})\gtrsim \log N.
\end{equation*}
For the last inequality, note that $1-x_{n_0} = 1-\theta_N^{-2} + O(N^{-2/3})$. On the other hand,
\begin{equation} \label{eq:thetam2}
  1 - \theta_N^{-2} = 2 w_N N^{-2/3} + O(w_N^2 N^{-4/3}).
\end{equation}
%so that, from \eqref{eq:thetam2},
%$1-x_{n_0} \sim 2w_N N^{-2/3}$. 

\begin{remark}
The logarithmic growth of $\sum_{i=1}^{N}g_{i+1}^{2}\E\xi _{i}^{2}$
is a consequence of our choosing $\theta _{N}$ local to one. Had it been
separated from one, the asymptotic variance of $\sum_{i=1}^{N}L_{i}$ would
be constant. This agrees well with the fact that linear spectral statistics
without singularities close to the edge of the semi-circle law do not need
scaling for the convergence to normality.
\end{remark}

Similarly, from \eqref{eq:ri-sums} and \eqref{upper bound on si}
\begin{equation*}
  \sum_{i=1}^{N}g_{i+1}^{4}\E\xi _{i}^{4}
  \lesssim \frac{1}{N^2}\sum_{i=2}^{N}\frac{1}{%
    \left( r_{i}-1\right) ^{4}}
  \lesssim \frac{1}{N^{1/3}w_N}.
\end{equation*}
Hence,%
\[
\sum_{i=1}^{N}g_{i+1}^{4}\E\xi _{i}^{4}/\left(
\sum_{i=1}^{N}g_{i+1}^{2}\E\xi _{i}^{2}\right) ^{2} \lesssim \frac{1}{%
N^{1/3}w_{N}\log ^{2}N}\rightarrow 0,
\]%
which establishes the Lyapunov condition (\ref{Lyapunov}). Let us now
approximate $\sum_{i=1}^{N}g_{i+1}^{2}\E\xi _{i}^{2}.$

Since, as we have just shown, $\sum_{i=1}^{N}g_{i+1}^{2}\E\xi
_{i}^{2}\gtrsim \log N,$ we will tolerate approximation errors of magnitude $%
o\left( \log N\right) .$ The following lemma is established in Subsection \ref{proof of lemma variance}.

\begin{lemma}
\label{lemma variance} Under assumptions of Lemma \ref{lemma si} for
all $1\leq i\leq N,$% 
\[
\E\xi _{i}^{2}=\frac{2\alpha}{N\theta _{N}^{2}r_{i}^{3}}\left(
  1+\varepsilon _{i}\right), \qquad
  %0 < \epsilon_i < \frac{1}{N(r_i-1)} \lesssim N^{-2/3}.
  |\varepsilon_i|<\frac{1}{N(r_i-1)}
  \lesssim N^{-2/3}.
\]%
\end{lemma}

Combining this lemma  with Lemma \ref{lemma si} yields
\begin{equation*}
  g_{i+1}^2 \E \xi_i^2
    =
    \begin{cases}
      \dfrac{\alpha}{2N\theta_N^2} \dfrac{r_{i+1}^2}{(r_{i+1}-1)^2
        r_i^3} \big[1+O(w_N^{-3/2}) \big] &  1 \leq i \leq N - N^{1/3}
      \\[10pt]
      O(N^{-1}(N-i+1)^2)   &      N - N^{1/3} \leq i \leq N      
    \end{cases}
\end{equation*}
Over the second range, the sum
$\sum g_{i+1}^2 \E \xi_i^2 = O(1)$, which will 
 be negligible. Over the first range, introduce
$f(x) = 1/[(r(x)-1)^2 r(x)]    =1/[(1-x)(1+\sqrt{1-x})]$. 
Monotonicity of $r(x)$ yields 
\begin{equation}
  \label{eq:fbds}
  f(x_i) \leq \frac{r_{i+1}^2}{(r_{i+1}-1)^2} \frac{1}{r_i^3} \leq
  f(x_{i+1}).
\end{equation}
Apply the trapezoidal rule bounds \eqref{eq:trapezoidal-bound}
with $|f''(x)| 
\lesssim (1-x)^{-3}$, and $\epsilon_{N1}(f) = O(N^{-2+4/3})$, while
$\epsilon_{N2}(f) = O(N^{-1+2/3})$.
With $x_a =0$ and 
$x_b = (N - [N^{1/3}])/(N\theta_N^2) = \theta_N^{-2} + O(N^{-2/3})$,
\begin{align*}
  \frac{1}{N \theta_N^2} \sum_{i=a}^b f(x_i)
  & = \int_{1-x_b}^1 \frac{\diff x}{(1+\sqrt{x})x} +O(N^{-1/3})
 = 2 \rho(x_b) + O(N^{-1/3}),
  % \\
  % & = 2 \log \Big[ \frac{1}{2} \Big(1+\frac{1}{\sqrt{1-x_b}}\Big) \Big] +O(N^{-1/3}).
\end{align*}
where $\rho(x) = \log [ \frac{1}{2} (1+1/\sqrt{1-x})]$.
Jumps of $1$ in $a, b$ to cover the two sides of
\eqref{eq:fbds} do  not alter the approximation.
In the range of interest, $0 \leq \rho'(x) \leq (1-x)^{-1} \leq
(1-\theta_N^{-2})^{-1} = O(w_N^{-1} N^{2/3})$, so that
$\rho(x_b) = \rho(\theta_N^{-2}) + O(w_N^{-1})$. 
In summary, 
\begin{align*}
  \sum_{i=1}^N g_{i+1}^2 \E \xi_i^2
    & = [\alpha \rho(\theta_N^{-2}) + O(w_N^{-1})][1+O(w_N^{-3/2})] +
      O(1) \\
    & = \alpha \rho(\theta_N^{-2}) + O(w_N^{-3/2} \log N).
    % & = \alpha \log \frac{\theta_N +
    %   \sqrt{\theta_N^2-1}}{2\sqrt{\theta_N^2-1}} + O\Big( \frac{\log
    %   N}{w_N^{3/2}}\Big).
\end{align*}

Recalling that $\tilde{\tau}_N^2 = \alpha \rho(\theta_N^{-2})$,
we have established the following theorem. 

\begin{theorem}
\label{theorem CLT for L}
If $\theta _{N}=1+N^{-2/3}w_{N}$ with $\left( \log \log N\right)
^{2}\ll w_{N}\ll \left( \log N\right) ^{2}.$ Then,
\begin{equation*}
  \tilde{\tau}_N^{-1} \sum_{i=1}^{N}L_{i}   \stackrel{d}{\to} \mathcal{N}(0,1).
\end{equation*}
% \[
% \frac{\sum_{i=1}^{N}L_{i}}{\sqrt{\alpha \log \frac{\theta _{N}+\sqrt{\theta
% _{N}^{2}-1}}{2\sqrt{\theta _{N}^{2}-1}}}}\overset{d}{\rightarrow} \mathcal{N}\left( 0,1\right) .
% \]
\end{theorem}

%We now turn to the second step of the analysis proposed in Section \ref{section away from edge}.

\subsubsection{Step 2a: Uniform bound on $L_{i}$}
\begin{lemma}
\label{Lemma Egor}If $\theta _{N}=1+N^{-2/3}w_{N}$ with $\left(\log
\log N\right)^{2}\ll w_{N}\ll \left(\log N\right)^{2}.$ Then,%
\[
\max_{1\leq i\leq N}\left\vert L_{i}\right\vert =o_{\Pr}\left( N^{-1/3}\right).
\]
\end{lemma}
To this end, we show that $\xi_i$ and $L_i$ are sub-gamma
  variables and apply exponential tail inequalities. The proof of the following lemma is given in \cref{proof of lemma L0}.

\begin{lemma}
\label{lemma L0}For any $1\leq i\leq N,$ $\xi_i\in SG\left(v_i,u_i\right)$ and $L_{i}\in SG\left(
v_{Li},u_{Li}\right) $ with%
\begin{equation}
v_i=\frac{2\alpha}{N\theta_N^2r_i^3},\quad v_{Li}=\frac{\alpha}{2N\theta _{N}^{2}\left( r_{i}-1\right) }\quad \text{ and }\quad u_i=u_{Li}=%
\frac{\alpha}{N\theta _{N}^{2}r_{i}^{2}}.  \label{medium i L0}
\end{equation}
\end{lemma}

\begin{proof}[Proof of Lemma \ref{Lemma Egor}]
Suppose that $L \in SG(v,u)$ and that $1 > 2vt > ut$, as will happen
in our uses below. Then
\begin{equation}
  \label{eq:blm}
  \Pr(|L| > 2\sqrt{2vt})
    \leq \Pr(|L| > \sqrt{2vt} +ut)
    \leq 2e^{-t},
\end{equation}
where the second inequality, valid for $t>0$, is essentially the
display prior to Theorem 2.3 in \cite{boucheron2013concentration}.
When applied to $L = L_i$, the bound
\begin{equation*}
  B_i(t) = 2\sqrt{2v_{Li}t} \lesssim \sqrt{t/N} (1-x_i)^{-1/4}
\end{equation*}
increases with $i$, and for $i \leq n_0 = N - \left[N^{1/3} \log^{2+\eta} N\right]$, it
is of order $\sqrt{t} N^{-1/3} \log^{-1/2-\eta/4}N$.
We may take a union bound over such $i$ in \eqref{eq:blm} by replacing
$t \leftarrow t + \log N$ because $B_{n_0}(t+\log N) = o(N^{-1/3})$.
More precisely, with $\eta = 1$ and $t = 5 \log \log N + \log N$,
there exists an absolute constant $C_1$ such that
\begin{equation}
  \label{eq:bdA}
  \Pr \left( \max_{i \leq n_0} |L_i| > \frac{C_1}{N^{1/3} \log^{1/4} N} \right) \leq \frac{2}{\log^5 N}.
\end{equation}

For $i$ close to $N$, however, this approach fails since
$B_N(t) \asymp \sqrt t N^{-1/3} w_N^{-1/4}$, which is no longer
$o(N^{-1/3})$ unless $t = O(\log \log N)$.

Instead, we use a simple form of chaining plus a version of the
Kolmogorov maximal inequality. We pick indices $n_0 < n_1 < \ldots <n_K =
N$ and use a bound of the form
\begin{equation}
  \label{eq:chain}
  \Pr \Big( \max_{i > n_0} |L_i| > \epsilon_1 + \epsilon_2 \Big)
  \leq \Pr \Big( \max_{0 \leq k < K} |L_{n_k}| > \epsilon_1 \Big)
       + \sum_{k=0}^{K-1} \Pr \Big( \max_{n_k < j \leq n_{k+1}}
       |L_j-L_{n_k}| > \epsilon_2 \Big).
\end{equation}

Iterating the relation $L_i = \xi_i + \gamma_i L_{i-1}$ and recalling definition \eqref{eq:gami-to-j} of $\gamma_{j:i}$, we have for
$i < j$
\begin{equation*}
  L_j = \gamma_{j:i+1} L_i + \sum_{k=i+1}^j \gamma_{j:k+1} \xi_k,
\end{equation*}
Since $\gamma_{j:i+1} < 1$, we then have
\begin{equation*}
  L_j - L_i < L_j/\gamma_{j:i+1} - L_i =: L_{i,j},
\end{equation*}
where
\begin{equation}
  \label{eq:Lj}
  L_{i,j} = \sum_{k=i+1}^j \xi_k/\gamma_{k:i+1}
\end{equation}
are partial sums of independent random variables.
Set $i = n_k$.  We then have
\begin{equation}
  \label{eq:maximal}
  \Pr \Big( \max_{n_k < j \leq n_{k+1}} |L_j - L_{n_k}| > 4 \epsilon \Big)
  \leq \Pr \Big( \max_{n_k < j \leq n_{k+1}} |L_{n_k,j}| > 4
  \epsilon \Big)
  \leq 4 \max_{n_k < j \leq n_{k+1}} \Pr ( |L_{n_k,j}| > \epsilon ),
\end{equation}
where the second inequality uses a maximal inequality, Theorem~1 from
\cite{etemadi1985some}. 

Now choose $K = [\log^5 N]$ and for $k = 1, \ldots, K$, let $n_k$ be the
closest integer to $n_0 + kN^{1/3} \log^{-2}N$. For these intervals 
  the products $\gamma_{j:i+1}$ are not too small: in
  Subsection \ref{proof of Lemma gammas Egor} we prove

\begin{lemma} \label{Lemma gammas Egor}
  Under assumptions of Lemma \ref{Lemma Egor}, if $N - N^{1/3} \log^3 N \leq i < j \leq N$ and
  $j-i \leq N^{1/3} \log^{-2}N+1$ then for large $N$ 
  \begin{equation*}
    \gamma_{j:i+1} \geq 1/2.
  \end{equation*}
\end{lemma}

For $ j \in (n_k,n_{k+1}]$ we therefore have $\gamma_{j:n_k+1} \geq
1/2$, and so from \eqref{eq:Lj} and \cref{lemma L0}
\begin{equation*}
  L_{n_k,j} \in SG \Big( \frac{8\alpha+1}{(N\theta_N)^{2/3} \log^2 N},
  \frac{2\alpha}{N \theta_N^2} \Big).
\end{equation*}
For some absolute constant $C_2$ and large $N$, the tail bound
\eqref{eq:blm} then implies, with $t = 10 \log \log N$,
\begin{equation*}
  \max_{n_k < j \leq n_{k+1}} \Pr \bigg( |L_{n_k,j}| >
  \frac{C_2 \sqrt{ \log \log N}}{N^{1/3} \log N} \bigg)
    \leq \frac{2}{\log^{10}N}.
\end{equation*}
With the same $t$ and recalling $B_{n_k}(t) \leq B_N(t) \asymp
\sqrt{t} N^{-1/3} w_N^{-1/2}$, we can find $C_3$ so that 
\begin{equation*}
  \Pr \bigg( |L_{n_k}| >
  \frac{C_3 \sqrt{ \log \log N}}{N^{1/3} w_N^{1/2}} \bigg)
    \leq \frac{2}{\log^{10}N}.
\end{equation*}
For $i > n_0$, use the last two bounds and the maximal inequality
\eqref{eq:maximal}, and for $i \leq n_0$ recall \eqref{eq:bdA}.
We conclude that on an event of probability at least $1 - 12
\log^{-5}N$ we have
\begin{equation}
\label{eq:result of lemma}
  \max_{i \leq N} |L_i| \leq C_4 N^{-1/3} \varepsilon_N
\end{equation}
where  $\varepsilon_N = \max ( \log^{-1/4} N, \sqrt{\log \log N}/w_N^{1/2} )
   \leq 1/\sqrt{\log \log N}$
under the assumptions on $w_N$. This completes the proof of Lemma
\ref{Lemma Egor}.
\end{proof}

\subsubsection{Step 2b: Uniform bound on $R_{i}$}

We show that $R_i$ is close enough to $L_i$ so that an analogous uniform
bound holds.
\begin{lemma}
\label{Lemma Egor R}Under the assumptions of Lemma \ref{Lemma Egor}, 
%with probability $1-o(1),$
\[
\max_{1\leq i\leq N}\left\vert R_{i}\right\vert =o_{\Pr}\left( N^{-1/3}\right)
. 
\]%
\end{lemma}

The starting point for analysis of the nonlinear process $R_i$ is
the perturbation representation $R_i = L_i + (T\varepsilon)_i$. Note that decomposition \eqref{eq:eps-def} expresses $\varepsilon_i$ in the
form $\varepsilon_i = \varepsilon(R_{i-1}; \beta_i,\gamma_i,\delta_i)$.
Consider a Winsorized version of $R_i$
\begin{equation*}
  \bar R_i = \phi_{N^{-1/3}/2}(R_i), \qquad
  \phi_u(x) =
  \begin{cases}
    -u & x<-u \\
    x & |x| \leq u \\
    u & x>u,
  \end{cases}
\end{equation*}
and create a modified series from $\bar{R}_{i-1}$:
\begin{equation}
  \label{eq:e-bar}
  \begin{split}
    \bar{\varepsilon}_i & = \varepsilon(\bar{R}_{i-1};
    \beta_i,\gamma_i,\delta_i) \\
    \tilde{R}_i  & = L_i + T \bar{\varepsilon}_i.
  \end{split}
\end{equation}
We will
show
\begin{equation}
  \label{eq:Rtildesmall}
  \max_{1 \leq i \leq N} | \tilde{R}_i | = o_\Pr(N^{-1/3}).
\end{equation}

A key observation is that on the event
$\tilde{\mathcal{R}}_N = \{ \max_{1 \leq i \leq N} | \tilde{R}_i |
\leq N^{-1/3}/2 \}$, we have $\tilde{R}_i = R_i$ for $i = 1, \ldots,
N$.
Indeed, $\bar \varepsilon_1 = 0$ since $\varepsilon_1 = 0$, and
$\tilde{R}_1 = L_1 = R_1$ so $|R_1| = |\tilde{R}_1| \leq N^{-1/3}/2$
on $\tilde{\mathcal{R}}_N$ and hence $\bar R_1 = R_1$.
Then \eqref{eq:e-bar} implies $\bar \varepsilon_2 = \varepsilon_2$, so
$\tilde{R}_2 = R_2$ and so again $|R_2| = |\tilde{R}_2| \leq
N^{-1/3}/2$, and so $\bar R_2 = R_2$.
Hence $\bar \varepsilon_3 = \varepsilon_3$ and so on, so that
eventually $\tilde{R}_i = R_i = \bar R_i$ for $i = 1, \ldots, N$. 

But this observation and \eqref{eq:Rtildesmall} imply Lemma \ref{Lemma Egor R}, since for any $\epsilon>0$
\begin{equation*}
  \Pr(\|R\|_\infty \leq \epsilon N^{-1/3})
  \geq \Pr(\|R\|_\infty \leq \epsilon N^{-1/3}, \tilde{\mathcal{R}}_N)
  = \Pr(\|\tilde{R}\|_\infty \leq \epsilon N^{-1/3}) \to 1.
\end{equation*}

We now outline the proof of \eqref{eq:Rtildesmall}, leaving more technical details to \cref{proof of eq 45}. By \cref{Lemma Egor}, it is sufficient to prove that $\max_{1\leq i\leq N}|T\bar{\varepsilon}_i|=o_{\Pr}(N^{-1/3})$. The terms  $\bar{\varepsilon}_i$ may be decomposed by
rewriting a Winsorized version of \eqref{eq:eps-def} as
follows, after defining $\bar{R}_{i-1}^{(1)} = \bar{R}_{i-1}/(1-\bar{R}_{i-1})$:
  \begin{equation}
    \label{eq:epsdecomp}
    \bar{\varepsilon}_i = - \delta_i + \beta_i \bar{R}_{i-1}^{(1)} +
    \bar{R}_{i-1}^{(1)}(\gamma_i \bar{R}_{i-1}^2 - \delta_i) + \gamma_i \bar{R}_{i-1}^2,
    \qquad i \geq 1.
  \end{equation}
Let $\bar{\varepsilon}_i^{\rm m}, \bar{\varepsilon}_i^{\rm s}$ and
$\bar{\varepsilon}_i^{\rm q}$ respectively denote the last three terms.
Using linearity of $T$, we have
\begin{equation}
\label{eq: epsdecomp short}
  T \bar{\varepsilon}_i = - T\delta_i + T\bar{\varepsilon}_i^{\rm m}
                    + T\bar{\varepsilon}_i^{\rm s}  + T\bar{\varepsilon}_i^{\rm
                      q}. 
\end{equation}

It is relatively straightforward to establish sufficiently tight bounds on the terms $T\delta_i$, $T\bar{\varepsilon}_i^{\rm s}$, $T\bar{\varepsilon}_i^{\rm q}$ using \eqref{eq:T-bound} (see \cref{proof of eq 45}). For $T\bar{\varepsilon}_i^{\rm m}$, note that $\bar{\varepsilon}_i^{\rm m} = \beta_i \bar{R}_{i-1}^{(1)}$ is
a martingale difference, 
since $\beta_i$ has mean $0$ and is independent of $\bar{R}_{i-1}^{(1)}$.
%,which is $\mathcal{F}_{i-1}$-measurable.
The term $T\bar{\varepsilon}_i^{\rm m}$ can be viewed as a special case of the quantity
\begin{equation*}
  T(\beta \mathsf{R})_i = \sum_{j=1}^i \gamma_{i:j+1} \mathsf{R}_j \beta_j,
\end{equation*}
with $\mathsf{R}_j$  measurable in the sigma-field generated by
$\alpha_1, \beta_1, \ldots, \alpha_{i-1}, \beta_{i-1}$ and
$\| \mathsf{R}_j \|_p \leq
\rho_N$.
%\footnote{In \eqref{eq:epsdecomp}, $\mathsf{R}_j = R_{j-1}/(1-R_{j-1})$, but the moment bound is
%not guaranteed.}
This is a sum $\sum_1^i X_j$ of martingale differences with $p$-th
moments, and the Marcinkiewicz-Zygmund-type
inequality of \cite[Theorem 2.1]{rio2009moment} says that
$\| \sum_1^i X_j \|_p^2 \leq (p-1) \sum_1^i \|X_j\|_p^2$. 

In \cref{proof of eq 45}, we use this inequality together with the Markov inequality and the union bound to show that there exists $C>0$ such that for all sufficiently large
$N$, with probability at least $1 - 1/N$, 
\begin{equation*}
  \max_i |T\bar \varepsilon_i^{\rm m} |
  \leq CN^{-2/3}\log^{3/2}N.%=O(N^{-2/3} \log^{3/2}N).
\end{equation*} 
Since the right hand side is obviously $o(N^{-1/3})$, this finishes the proof of \eqref{eq:Rtildesmall}.

\subsubsection{Step 3: Linear approximation for $\log |M_{N}| $}
Recall that 
\[
\log \left\vert M_{N}\right\vert =\sum\nolimits_{i=1}^{N}\log \left\vert
  1-R_{i}\right\vert.
\]%
Since \( \max_{i} |R_i| = o_{\Pr}(N^{-1/3}) \), we have a uniform Taylor's approximation
\[
    \log |1 - R_i| = -R_i - {R_i}^{2}/2 + o_\Pr(N^{-1}) \, .
\]
Summing up,
\begin{equation}\label{log expansion}
    \log |M_N| = \sum_{i = 1}^{N} (-R_i - {R_i^2}/{2}) + o_{\Pr}(1) \, .
\end{equation}
In the rest of this subsection, our goal is to show that we can replace each term \( -R_i - R_i^2 / 2\) with the linear process \( L_i\), {with inclusion of a deterministic shift. To be precise, we will show that}
\begin{equation}\label{replacing R by L}
    \sum_{i = 1}^{N} (-R_{i} - R_{i}^{2} / 2) + \sum_{i = 1}^{N} L_{i} = {\frac{1 - \alpha}{6} \log N} + O_{\Pr}(\log\log N) \, .
\end{equation}

Similarly to \eqref{eq: epsdecomp short}, we have
\begin{equation}
\label{eq:epsdecomp new}
T \varepsilon_i = - T\delta_i + T\varepsilon_i^{\rm m}
                    + T\varepsilon_i^{\rm s}  + T\varepsilon_i^{\rm
                      q},
\end{equation}
where with the notation $R^{(1)}_{j}=R_j/(1-R_j)$, we have
\[
\varepsilon_i^{\rm m}=\beta_iR^{(1)}_{i-1},\quad \varepsilon_i^{\rm s}=R^{(1)}_{i-1}(\gamma_iR_{i-1}^2-\delta_i),\quad \varepsilon_i^{\rm
                      q}=\gamma_iR_{i-1}^2.
\]
The decomposition \eqref{eq:epsdecomp new} leads to
\begin{equation*}
  \sum R_i = \sum L_i + T \varepsilon_i^{\rm m} + T \varepsilon_i^{\rm s}
                      + T \varepsilon_i^{\rm q} - T \delta_i.
\end{equation*}
We will see that $\sum T \varepsilon_i^{\rm m} + T \varepsilon_i^{\rm s} =
O_\Pr(1)$. Note however that $\varepsilon_i^{\rm q} = \gamma_i R_{i-1}^2$
are positive, and will contribute to the deterministic shift.
We therefore further decompose $T \varepsilon_i^{\rm q}$ using
\begin{equation*}
  \varepsilon_i^{\rm q}
  = \varepsilon_i^{\rm qd}
  + (\varepsilon_i^{\rm qL} - \varepsilon_i^{\rm qE})
  + \varepsilon_i^{\rm qE},
\end{equation*}
where the \textbf{q}uadratic ``\textbf{d}ifference'', ``\textbf{L}inear
approximation'' and ``\textbf{E}xpectation'' terms are respectively
given by
\begin{equation*}
  \varepsilon_i^{\rm qd} = \gamma_i(R_{i-1}^2- L_{i-1}^2) \qquad
  \varepsilon_i^{\rm qL} = \gamma_i L_{i-1}^2 \qquad
  \varepsilon_i^{\rm qE} = \gamma_i \E L_{i-1}^2.
\end{equation*}

An analysis similar to one we used to prove Lemma~\ref{Lemma Egor R} leads to the following lemma. The proof is rather technical and is given in Subsection \ref{proof of lemma star}.

\begin{lemma}
\label{Lemma star}
Under the assumptions of Lemma \ref{Lemma Egor},%
\begin{align}
  & \sum_i T \varepsilon_i^{\rm m} + T \varepsilon_i^{\rm s}
    + T \varepsilon_i^{\rm qd} = O_\Pr(1), \label{eq:msqd sum}\\
  & \sum_i R_i^2 = O_\Pr(1), \label{eq: R2 sum}\\
    &\sum_i T \varepsilon_i^{\rm qL} - T \varepsilon_i^{\rm qE}
        = o_\Pr(1)  \label{eq:TeqL}\\
  &\sum_i T \varepsilon_i^{\rm qE}
         = \alpha \sum_i T \delta_i +  o_\Pr(1).  \label{eq:TeqE}
\end{align}
\end{lemma}

Combining the results of \cref{Lemma star}, we arrive at
\begin{equation*}
  \sum_i R_i + R_i^2/2
    = \sum_i L_i + (\alpha-1) \sum_i T \delta_i + O_\Pr(1).
\end{equation*}
Remarkably, for $\alpha = 1$, this will be the end of the proof. When
$\alpha \neq 1$, the remaining sum \( (\alpha - 1) \sum_{i = 1}^{N}
T \delta_i \) results in an additional shift. The proof of the
following lemma is postponed to Section~\ref{proof lemma Tdeltai}. 

\begin{lemma}\label{lemma Tdeltai}
It holds, for large enough $N$,
\[
    \sum_{i = 1}^{N} T\delta_i = \frac{1}{6} \log N + O(\log\log N) \, .
\]
\end{lemma}

Now the CLT for $\log|M_N|$ follows from Theorem ~\ref{theorem CLT for
  L}.

\begin{corollary}
\label{theorem subordinate}Under the assumptions of Theorem \ref{theorem CLT for L},%
\[
\frac{\log \left\vert M_{N}\right\vert +\frac{\alpha-1}{6} \log N} {\tilde{\tau}_N}\overset{d}{\rightarrow} \mathcal{N}\left(
0,1\right) .
\]
\end{corollary}

\subsubsection{CLT for  $\log |D_{N}| $}
\label{sec:CLT forD}

From \eqref{new scaling new} and \eqref{eq:rmdef},
  $D_N = M_N \theta_N^N \prod_1^N r_i$. Hence
%From the definitions of $\mathcal{D}_N$ and $E_N$, we have
\begin{equation}
    \label{Dn}
\log \left\vert D_{N}\right\vert =\log \left\vert M_{N}\right\vert
+ N\log \theta _{N}
+\sum\nolimits_{i=1}^{N} \log \left( 1+\sqrt{1-x_{i}}\right). 
\end{equation}%
In the trapezoidal approximation to
$\Delta_N \sum_1^N \log(1+\sqrt{1-x_i})$ we have
$\varepsilon_{N2}(f) = O(\Delta_N)$ and,
since $|f''(x)| \asymp (1-x)^{-3/2}$ for $0 < x < 1$, also
$\varepsilon_{N1}(f) = O(N^{1/3} \Delta_N^2)$.
So from \eqref{eq:trapezoidal-bound}
\begin{align*}
  \sum_1^N \log(1+\sqrt{1-x_i})
  &  = N \theta_N^2 \int_{1-x_N}^1 \log(1+ \sqrt u) \diff u + O(1) \\
  & = \tfrac{1}{2} N \theta_N^2 - \tfrac{2}{3}(2w_N)^{3/2} + O(1).
%  & = \frac{N \theta_N^2}{2} - \frac{2}{3}(2w_N)^{3/2} + O(1).
\end{align*}
At the second line we used
$\int_0^1 \log(1+\sqrt u) \diff u = \frac{1}{2}$ and, with
$a_N = 1-x_N = 2w_N N^{-2/3} + O(N^{-1})$, %a_N^2)$, 
 also
$\int_0^{a_N} \log(1+\sqrt u) \diff u = \frac{2}{3}a_N^{3/2} + O(a_N^2)$.

Using this in \eqref{Dn} together with \eqref{eq:thetam2}
 and
\[
N \log \theta_{N} =
w_{N}N^{1/3}+O\big( w_{N}^{2}N^{-1/3}\big)
\]
yields, recalling that $2w_N=\sigma_N$,
\[
\log \left\vert D_{N}\right\vert =\log \left\vert
  M_{N}\right\vert+ \tfrac{1}{2}N + \sigma_N N^{1/3}-\tfrac{2}{3}
  \sigma_N^{3/2} + O(1).
\]
Using Corollary \ref{theorem subordinate}, we obtain Theorem \ref{theorem main}.

\subsection{Proofs from Section~\ref{section up to edge}}
\label{appendix up to edge}

\subsubsection{Preliminaries}
\textit{Correlation functions.} \ Let $P_N(x_1,...,x_N)$ be a joint density of \textit{unordered}
eigenvalues $l_1,...,l_N$ of G(U/O)E (scaled so that $\max l_i$ is close
to 2 for large $N$). Following \cite{tracy1998correlation}, the
$k$-point correlation function is defined as 
\[
R_k(x_1,...,x_k)=\frac{N!}{(N-k)!}\int...\int P_N(x_1,...,x_N)\mathrm{d}x_{k+1}...\mathrm{d}x_N.
\]
Note that this is not a probability density: it has total integral $N!/(N-k)!$.

For any integrable function $F(x_1,...,x_k)$, we have
\begin{equation}
  \label{Iain1}
  \E F(l_1,...,l_k)
  % =\int...\int F(x_1,...,x_k)P_N(x_1,...,x_N)\mathrm{d}x_1...\mathrm{d}x_N \notag \\
=\frac{(N-k)!}{N!}\int...\int F(x_1,...,x_k)R_k(x_1,...,x_k)\mathrm{d}x_1...\mathrm{d}x_k.
\end{equation}
% \begin{eqnarray}
% \E F(l_1,...,l_k)=\int...\int F(x_1,...,x_k)P_N(x_1,...,x_N)\mathrm{d}x_1...\mathrm{d}x_N \notag \\
% =\frac{(N-k)!}{N!}\int...\int F(x_1,...,x_k)R_k(x_1,...,x_k)\mathrm{d}x_1...\mathrm{d}x_k. \label{Iain1}
% \end{eqnarray}

Write $\rho_N(\lambda) = \rho_{N,\alpha}(\lambda) = N^{-1}
R_1(\lambda)$ for the normalized one-point correlation function,
interpreted as the ``mean density'' of the eigenvalues. The expected
value of a linear spectral statistic can be written as
\begin{equation}
  \label{eq:linear-stat}
  \E \Big[ N^{-1} \sum_{i=1}^N f(\lambda_i) \Big]
  =   \E \Big[ N^{-1} \sum_{i=1}^N f(l_i) \Big]
  = \int f(\lambda)
    \rho_{N,\alpha}(\lambda) \diff \lambda. 
\end{equation}

A key tool in approximating such expectations will be a uniform bound,
due to G\"otze and Tikhomirov, for 
the deviation of the one-point function in GUE %(and GOE)
from the semi-circle density
$p_{\rm SC}(x) = (2\pi)^{-1} \sqrt{4-x^2} \mathbf{1}_{|x| \leq 2}$.
Indeed, \cite[Theorem 1.2]{Gotze2005}
show the existence of
positive absolute constants $a, A$
% $\gamma, C > 0$
such that for all $|x| \leq 2 - aN^{-2/3} $,
\begin{equation}\label{goetze_tihkomirov1}
    |\rho_{N}(x) - p_{SC}(x)| \leq \frac{A}{N(4 - x^2)} \, .
\end{equation}

% \fixm{\texttt{
% NB: The following equation was excised, but seems to be needed for an argument in \cref{prf:Efeta}, where it is referenced.
% In the original text, this referred to Section 6.5 of our transition paper.
% }}

% \fixb{At the edge, the one-point function decays at least exponentially:
%  for large enough $N$,
% for all $ s > -\gamma $,
% \begin{equation}\label{rhoNedge}
%     \rho_{N}(2 + sN^{-2/3}) \leq C(\gamma) N^{-1/3} e^{-2s}\, ,
% \end{equation}
% we refer to \cref{sec:1pt-edge-decay} of for details.}

Determinantal correlation functions imply the following elementary variance bound, which we prove in \cref{proof:lem:variance-bound}.

\begin{lemma}
  \label{lem:variance-bound}   For eigenvalues from GUE,
  \begin{equation}
    \label{variance bound}
    \mathrm{Var}\Big[N^{-1}\sum_{i=1}^{N}f(l_i)\Big] \leq
  N^{-1}\int f^{2}(x)\rho_N(x)\mathrm{d}x. 
\end{equation}
\end{lemma}

\begin{remark}
For the usual linear statistic, with $f$ not depending on $N$ and
analytic in a neighborhood of $[-2,2]$, this is a terrible bound
since then
$\mathrm{Var}\left[\sum\nolimits_{1}^{N}f(l_i)\right]=O(1)$. But in
our critical case settings, it seems to give the right order, and will
become useful below. 
\end{remark}

\textit{Edge bounds.} \ For both GUE and GOE, we have
\begin{equation}
  \label{eq:onept-bds}
  \rho_{N,\alpha}(2 + s N^{-2/3})
    \leq
    \begin{cases}
      C N^{-1/3} e^{-2s}  &  \qquad \qquad \qquad s > -1 \\
      C N^{-1/3} |s|^{1/2} &   \qquad  -N^{2/3-\varepsilon} < s \leq -1, 
    \end{cases}
\end{equation}
for any $0<\varepsilon<2/3$.
The bounds for GUE follow directly from the G\"{o}tze-Tikhomirov
bounds and Tracy-Widom asymptotics of the Hermite functions at the
edge. 
%(2.3-2.4) of [JKOP21].
For GOE the one-point function $\rho_{N,2}$ differs from
$\rho_{N,1}$ by a term involving integrals of scaled Hermite
functions, and this can again be analyzed by bounds on Hermite
polynomials. These bounds are implicit in \cite{johnstone2012fast} and \cite{Gotze2005}, but for the
reader's convenience some discussion appears in \cref{sec:edge-bounds-one}. %\textcolor{red}{[question: are the bounds implicit in Gotze-Tikhomirov? Do they consider region outside the bulk?]}

\subsubsection{Gaussian Non-concentration}
\label{sec:gauss-non-conc}
Let us introduce new notation
\begin{equation}
\label{def: sigmaN check}
\check{\sigma}_N=(\log N)^{O(\log\log N)}.
\end{equation}

\begin{lemma}
  \label{lem:anticoncGOE}
  Suppose that $W_N$ is a matrix drawn from either GOE or GUE,
  divided by $\sqrt{N}$ to have support
  of the limiting spectral distribution on $[-2,2]$.
  Suppose also that
$|E - 2| \leq N^{-2/3} \check{\sigma}_N$. 
  Then for each $c_0 > 0$ and each $d \in (0, c_0)$, we have
  for $N > N(d)$,
    \begin{equation}  \label{eq:anticonc}
      \mathbf{P} ( \min_j |\lambda_j - E| \leq N^{-2/3 - c_0}) \leq N^{-d}.
    \end{equation}
  % (ii) If $\gamma>0$ is fixed and $E > - \gamma$, then for each $c_1 >
  % 0$ small, there exists $C = C(\gamma)$ such that for $N$ large 
  %   \begin{equation}  \label{eq:anticonc1}
  %     \mathbf{P} ( \min_j |\lambda_j - E| \leq c_1 N^{-2/3}) \leq Cc_1.
  %   \end{equation}
\end{lemma}

\begin{proof}
%Let $ \rho_{N, j}(x) $ be the one-point correlation function for GOE,
%$j = 1$, and GUE, $j = 2$.
% both ensambles are scaled to have support
% of the limiting distribution on $[-2,2]$.
%Case (i) is  
The lemma is an immediate consequence of the bounds
\begin{equation}\label{rho12_bound}
    \rho_{N, \alpha}(E) \lesssim N^{-1/3} \check{\sigma}_N^{1/2},
    \qquad \alpha = 1, 2 \, ,
\end{equation}
holding uniformly for all $ E: |E - 2| \leq N^{-2/3}
\check{\sigma}_N  $.
Bounds \eqref{rho12_bound} follow directly from \eqref{eq:onept-bds}.
%are established in Section \ref{sec:edge-bounds-one}.

Indeed, for a matrix $M$ and set $I$, let the number of eigenvalues of $M$ in
$I$ be $N_M(I)$.
%\[
%N_M(I)=\operatorname{tr}\left[\mathbf{1}_I(M)\right].
%\]
Set $I = [E-N^{-2/3 - c_0},E+N^{-2/3 - c_0}]$. We have
\begin{align*}
  \mathbf{P} ( \min_j |\lambda_j - E| \leq N^{-2/3 - c_0})
  & = \mathbf{P} (N_{W_N}(I) \geq 1)  \\
  & \leq \E N_{W_N}(I)
  = N \int_I \rho_{N, \alpha}(E) \diff E
  \lesssim 2 \check{\sigma}_N^{1/2} N^{-c_0}
  \leq N^{-d}
\end{align*}
for each $d < c_0$ and $N$ large.
% Case (ii) follows similarly, replacing $N^{-c_0}$ by $c_1$, and
% also $\check{\sigma}_N^{1/2}$ by $C(\gamma)$ in the
% last line. 
\end{proof}

\subsubsection{Proof of Proposition   \ref{lemma stieltjes}}
First, we prove the proposition for GUE case $\alpha = 1$.
To establish Proposition~\ref{lemma stieltjes} via \eqref{variance bound}, we will calculate the expectation of the truncated statistics 
\begin{equation*}
    L_{lN} = \frac{1}{N} \sum_{j = 1}^{N} f_c^l(\lambda_j),
    \qquad l = 1, 2, 
\end{equation*}
where
\begin{equation*}
  f_c(\lambda) = \frac{1}{E-\lambda} \bo \{|E-\lambda| > c N^{-2/3} \}.
\end{equation*}
% $$f_{\eta}(\lambda) = \phi(2 + N^{-2/3} \tilde{\sigma}_N - \lambda),
% \qquad
% \phi(t) = \frac{1}{x} \Ind_{|x| \geq \eta} .
% $$
The truncation makes the function  integrable with respect to the
density $ \rho_{N}(x) $.

The truncation is typically harmless: if there are no eigenvalues near
$E$, more precisely if $N_W(E-cN^{-2/3},E+cN^{-2/3}) = 0$, then
\begin{equation*}
   \sum_{j = 1}^{N} f_c^l(\lambda_j) =  \sum_{j = 1}^{N} (E-\lambda_j)^{-l}.
\end{equation*}
Lemma \ref{edge spacing lemma} assures that for $\varepsilon > 0$
small, there exists $c = c(\varepsilon,\gamma)$ small so that equality
holds with probability at least $1-\varepsilon$. 
Therefore, 
%Criterion \fix{C.2} for bounds in probability then shows that it
it suffices to show that Proposition \ref{lemma stieltjes} holds with
$N^{-1} \sum (E-\lambda_j)^{-l}$ replaced by $L_{lN}(c)$ for each $c > 0$
fixed. 

The main work, contained in the next lemma, is to control the expected
values of $L_{lN}$.

\begin{lemma}\label{lemma Efeta}
%  Let $ c N^{-2/3} < \eta < C N^{-2/3}$ for some fixed $c, C > 0$.
  Suppose that $ \sigma_N $ satisfies \eqref{eq:sigcon}.
  Then for each $c > 0$ we have
\[
    \E L_{lN} = \E \frac{1}{N} \sum_{j = 1}^{N} f_c^l(\lambda_j) = \left\{
    \begin{aligned}
    & 1 + O\left((1 + |\sigma_{N}|^{1/2})N^{-1/3}\right),
    \qquad& &l = 1 \\ 
    & O(N^{1/3}), \qquad && l = 2 \, . 
    \end{aligned}
    \right.
\]
\end{lemma}

Lemma \ref{lemma Efeta} and the variance bound \eqref{variance bound}
quickly yield Proposition \ref{lemma stieltjes}.
Indeed, for $L_{2N} > 0$, Proposition \ref{lemma stieltjes} holds since
$ L_{2N} = O_\Pr( \E L_{2N})$.
For $L_{1N}$, bound \eqref{variance bound} implies
\begin{equation*}
  \Var( L_{1N}) \leq N^{-1} \E L_{2N} = O(N^{-2/3}),
\end{equation*}
and we conclude Proposition \ref{lemma stieltjes} from
$ L_{1N} - \E L_{1N} = O_\Pr(\sqrt{ \Var( L_{1N})} )$. 

\begin{proof}[Proof of Lemma \ref{lemma Efeta} (for GUE)] 
  \label{prf:Efeta}
  First let us bound the error of replacing $ \rho_{N} $ with $ \psc $ in the integral
\[
    \E \frac{1}{N} \sum_{i = 1}^{N} f_c^l(\lambda_j) = \int
    f_c^l(\lambda) \rho_{N}(\lambda) d \lambda \, . 
\]
Abbreviate $\epsilon_{N} = N^{-2/3} \sigma_N$,  $\delta_N =
a N^{-2/3}$,  and decompose $\R$ into
$I_N = [-2+\delta_N,2-\delta_N]$ along with
$J_N = (2-\delta_N,\infty)$ and $J_N^- = (-\infty, -2+\delta_N)$,
and write $ g = f_c^l $. Then,
\begin{equation*}
  \int g \rho_N - g \psc
  = \int_{I_N} g(\rho_N - \psc)
     + \int_{J_N \cup J_N^-} g \rho_N 
     - \int_{J_N \cup J_N^-} g \psc. 
\end{equation*}
First, we have by \eqref{eq:onept-bds},
\begin{align*}
  \int_{J_{N}} |g| \rho_N
  & \leq (\sup_{J_N} |g|) \int_{J_N} \rho_{N}
    \lesssim N^{2l/3 - 2/3}  \int_{- a}^{\infty} \rho_N(2 +
    sN^{-2/3}) \mathrm{d}s \\
  & \lesssim N^{2l/3 - 1} \int_{-a}^{\infty} e^{-2s} \mathrm{d}s 
    \lesssim N^{2l/3 - 1} \, .
\end{align*}

Similar bounds hold for the integrals over $ J_{N}^{-} $ and for those
with respect to $\psc$.
For the middle interval, we use the G\"{o}tze and Tikhomirov bound
\eqref{goetze_tihkomirov1}, 
\[
  \int_{I_N} |g (\rho_N - \psc)|
  \lesssim \frac{1}{N} \int_{-2 + \delta_N}^{2 - \delta_N} \frac{|f_{c}(\lambda)^{l}|}{4 - \lambda^2} \mathrm{d}\lambda.
    %\lesssim \frac{1}{N}
    %\int_{-2 + \delta_N}^{0} \frac{f_{\eta}(\lambda)^{l}}{2 + \lambda} \mathrm{d}\lambda +
    %\frac{1}{N} \int_{0}^{2 - \delta_N} \frac{f_{\eta}(\lambda)^{l}}{2 - \lambda} \mathrm{d}\lambda \, .
\]
Observe that on $ [-2 + \delta_N, 0] $,  we have $ 0 \leq
f_c^l(\lambda) \leq 1 $.
Therefore, 
\[
    \frac{1}{N} \int_{-2 + \delta_N}^{0}
    \frac{|f_c^l(\lambda)|}{4- \lambda^2} \mathrm{d}\lambda
    \leq \frac{1}{N}
    \int_{-2 + \delta_N}^{0} \frac{\mathrm{d}\lambda}{2 +
      \lambda} = O(N^{-1} \log \delta_N) = o(N^{-1/3}) \, . 
\]
To deal with the remaining part of the integral,
make the change of variable $\lambda = 2 - uN^{-2/3}$ and note that
$f_c(\lambda) = (u+\sigma_N)^{-1} N^{2/3}$ except for $|u+\sigma_N|
\leq c$, where it vanishes.
Thus
\begin{align*}
  \int_0^{2-\delta_N} \frac{|f_c^l(\lambda)|}{4-\lambda^2} \diff \lambda
  &\leq \int_0^{2-\delta_N} \frac{|f_c^l(\lambda)|}{2-\lambda} \diff \lambda\\
  &=N^{2l/3} \int_a^{2N^{2/3}}
    \frac{\bo\{|u+\sigma_N|\geq c\}}{|u+\sigma_N|^l} \frac{\diff u}{u}
  \leq C N^{2l/3},  
\end{align*}
where, for example, we may take $C = C_{\gamma,c}
= c^{-l} \int_a^{\gamma+1} u^{-1} \diff u + \int_{\gamma+1}^\infty (u-\gamma)^{-l}
u^{-1} \diff u$.  

So we have established the following approximation
\[
  \int f_c^l(\lambda) \rho_{N}(\lambda) d \lambda
  = \int f_c^l(\lambda) \psc(\lambda) d \lambda + O\left( N^{2l/3 - 1} \right),
\]
and it remains to analyse the integral with respect to the semi-circle density. 

Let $m(z) = \int (\lambda-z)^{-1} \psc(\lambda) \diff \lambda
= (-z+\sqrt{z^2-4})/2$ denote the Stieltjes transform of
$\psc(\lambda)$.
When $\sigma_N \geq c$, we simply have, since $E = 2 + \epsilon_N
= 2 + \sigma_N N^{-2/3}$,
\begin{equation*}
  \int f_c^l(\lambda) \psc (\lambda) \diff \lambda
  = \int_{-2}^2 \frac{1}{(E-\lambda)^l}\psc (\lambda) \diff \lambda
  =
  \begin{cases}
    -m(2+\epsilon_N) = 1+O(\epsilon_N^{1/2})  & l=1 \\
    \ \ m'(2+\epsilon_N) = O(\epsilon_N^{-1/2}) & l = 2.
  \end{cases}
\end{equation*}

When $\sigma_N < c$, note that $f_c(\lambda) = 0$ for $\lambda \in
[\lambda_-,\lambda_+]$ with
$\lambda_\pm = 2 + N^{-2/3}(\sigma_N \pm c)$, and 
in particular $\lambda_- < 2$.
The square-root decay of $\psc$ near $2$ implies that
\begin{equation}
  \label{eq:fcl-right}
  \int_{\lambda_-}^2 f_c^l(\lambda) \psc (\lambda) \diff \lambda
  \leq c^{-l} N^{2l/3} \int_{\lambda_-}^2 \psc (\lambda) \diff \lambda
  \lesssim N^{2l/3-1}.
\end{equation}
Consider now $l=1$. The change of variable
$\lambda = \lambda_- - N^{-2/3}x=E-N^{-2/3}(x+c)$ yields, along with
$\psc \lesssim \sqrt{2-\lambda}$,
\begin{align*}
  \int_{-2}^{\lambda_-} \Big| f_c(\lambda) - \frac{1}{2-\lambda} \Big|
  \psc(\lambda) \diff \lambda
  & \lesssim |\sigma_N| N^{-2/3} \int_{-2}^{\lambda_-} \frac{\diff
    \lambda}{(E-\lambda)\sqrt{2-\lambda}} \\
  & \leq |\sigma_N| N^{-1/3} \int_0^\infty \frac{\diff x}{(x+c)
    \sqrt{x}}
    \lesssim N^{-1/3}.
\end{align*}
The same bound for $\psc (\lambda)$ also gives
\begin{equation*}
  \int_{-2}^{\lambda_-} \frac{\psc(\lambda)}{2-\lambda}  \diff \lambda
   = -m(2) - \int_{\lambda_-}^2 \frac{\psc(\lambda)}{2-\lambda}
   \diff \lambda
   = 1 + O(N^{-1/3}).
\end{equation*}
Combining the last three displays yields, for $-\gamma \leq \sigma_N<c$,
\begin{equation*}
  \int_{-2}^2 f_c(\lambda) \psc(\lambda) \diff \lambda
    = 1 + O_{c,\gamma}(N^{-1/3}).
\end{equation*}

We turn to $l=2$, still with $\sigma_N < c$. Again setting $\lambda =
\lambda_- - N^{-2/3}x$, we have
\begin{align*}
  \int_{-2}^{\lambda_-} f_c^2(\lambda) \psc(\lambda) \diff \lambda
  & \lesssim \int_{-2}^{\lambda_-} (2-\lambda+N^{-2/3} \sigma_N)^{-2}
    \sqrt{2-\lambda} \diff \lambda \\
  & \lesssim  N^{1/3} \int_0^\infty (c+x)^{-2}(x+c-\sigma_N)^{1/2}
    \diff x
    \lesssim N^{1/3}.
\end{align*}
Together with \eqref{eq:fcl-right} this shows that
$\int f_c^2(\lambda) \psc(\lambda) \diff \lambda = O(N^{1/3})$ and
completes the proof. 
\end{proof}

This finishes our proof of Proposition \ref{lemma stieltjes} for the GUE case.

For the GOE case, the proposition follows from the following theorem. A proof of this theorem can be found in Section~\ref{section gue to goe}.

\begin{theorem}
\label{theorem Damian}
Let $M_N^{\mathbb{C}}$ and $M_N^{\mathbb{R}}$ be $N\times N$ (unscaled) GUE and GOE matrices, respectively. Suppose that $f_N$ is a series of functions such that
\[
f_N\left(M_N^{\mathbb{C}}\right)=a_N+O_{\Pr}(b_N),
\]
for some sequences $a_N$ and $b_N$. Then,
\[
f_N\left(M_N^{\mathbb{R}}\right)
  = a_N+O_{\Pr}(b_N+\mathrm{TV}(f_N)),
\]
where $\mathrm{TV}(f_N)$ is the total variation of $f_N$.
\end{theorem}
Indeed, the theorem and the fact that scaling of the argument does not
change the total variation of a function yield the equivalents of
Lemmas \ref{lem:variance-bound} and \ref{lemma Efeta}
% and Corollary \ref{corollary variance}
for GOE. These equivalents, combined with the anticoncentration bound
of Lemma \ref{edge spacing lemma} 
%the following lemma,
imply Proposition \ref{lemma stieltjes} for GOE.

\section{Extension to Wigner matrices}
%\section{W: Modified Introduction}
\label{sec:modif-intr}

\subsection{Lindeberg swapping formalism for asymptotically flat  \(Q\)}
\label{sec:lindeberg-swapping}

~\newline
\textit{Definitions.} A Wigner matrix is an Hermitian $N \times N$ matrix
$W_N=(\xi_{ij}/\sqrt{N})$ satisfying

\quad (i) the upper-triangular components
$\{\operatorname{Re}\xi_{ij},\operatorname{Im}\xi_{ij}\}_{i<j}$ and 
$\{\xi_{ii}\}$ are independent
\newline  random variables with mean zero,

\quad (ii) $\mathbf{E}|\xi_{ij}|^2=1$ for $i\neq j$ and
$\mathbf{E}\xi_{ii}^2\leq B$ for some absolute constant $B$;

\quad (iii) a moment bound uniform in $N$: for all \(p \in \Z_{>0}\), there
is a constant \(C_p\) such that 
\begin{equation} \label{eq:W3}
  \E \abs{\Re\xi_{ij}}^p,
  \E \abs{\Im\xi_{ij}}^p 
  \leq C_p.
\end{equation}
%\end{definition}

This definition is standard, e.g.
\cite[][Def 2.2]{benaych2018}, except that we also require
independence of $\Re\xi_{ij}$ and $\Im \xi_{ij}$ to simplify our
swapping arguments.
Condition (ii) allows for zero variances on the diagonal, as in the
SSK model of \cite{Kosterlitz1976}.

The moments of two Wigner matrices $W_N, W_N'$
\textit{match to order $m$} if for  integer $0 < a \leq m$
\begin{equation*}
  \E (\Re \xi_{ij})^a = \E (\Re \xi_{ij}')^a, \qquad
  \E (\Im \xi_{ij})^a = \E (\Im \xi_{ij}')^a
\end{equation*}
for all $1 \leq i < j \leq N$. Note that this constrains only the
\textit{off-diagonal} entries. The diagonal entries already match to
order one by assumption, which is all that
we need.

An event sequence \(A_N\) holds \textit{with high
  probability} if there exists a \(d > 0\) such that 
  \[
    \Pr(A_N\comp) \lesssim N^{-d}.
  \]
An event \(B_N\) holds \textit{with overwhelming probability} (w.o.p.)
if, for all \(A > 0\), 
  \[
    \Pr(B_N\comp) \lesssim N^{-A}.
  \]
If $X_N \lesssim c_N$ w.o.p. and there are constants $C_0, C_2$ such
that eventually $c_N \geq N^{-C_0}$ and $\E X_N \leq N^{C_2}$, then
$\E X_N \lesssim c_N$. [For proof, see e.g. \cite[][Lemma
7.1]{benaych2018}.]
Here  and later ``$X_N\lesssim c_N$ w.o.p.'' means
that there exists $C$ such that event $X_N\leq C c_N$ holds
w.o.p. Similarly for statements like $X_N = O(c_N)$ w.o.p. 
% in what follows, statements like
% ``$a_N=b_N+O(N^\delta)$ w.o.p.'' mean that there exist $C$ such that
% event $|a_N-b_N|\leq CN^\delta$ holds w.o.p.  

% Throughout, we use the following notation \textcolor{red}{[Possibly,move definitions of some standard notation to the end of the introduction session]}:
% \[
%   \check{\sigma}_N
%   = (\log N)^{O(\log\log N)}.
% \]

 \medskip
 In \Cref{prop:multi-matching} and its consequence
 Proposition
 %\ref{prop:multi-matching} and
 \ref{prop:joint-convergence}, we 
make the swapping argument explicit for abstract
\(Q\) satisfying generic asymptotic `flatness' derivative bounds.
In the next subsection we assemble tools -- resolvent perturbation and
local law -- with the goal of establishing, in Proposition \ref{prop:system}, 
the necessary flatness bounds for some specific choices
of \(Q\) needed for our later applications.

Fix $c_0>0$ and set $ \| F \|_{c_0} = \sup \{ |F(t)|, |t| \leq   N^{c_0} \}$.
Let  \(\delta_N \rightarrow 0\) in such a way that \(\delta_N \gtrsim
N^{-c_1}\) for some \(c_1 > 0\).
Let $Q$ be a function on $N \times N$ Hermitian/symmetric matrices
taking values in $[0,1]$.
Let Wigner matrices $W_N, W_N'$ be given and define
$Q_\gamma(t) = Q(W_t^\gamma)$ as in Section \ref{sec:outline-approach}.
We say that $Q$ satisfies \textbf{condition F} or $F(\delta_N)$ if for
all $\gamma$ and $1 \leq k \leq 4$ we have w.o.p. that 
  \begin{align}
    \label{eq:D-cond-1}
    \norm{Q_{\gamma}^{(k)}}_{c_0}
    \lesssim N^{-\frac{k}{2}} \delta_{N}.  \tag{F}
  \end{align}

% \begin{proposition}
%   \label{prop:single-matching}
%   Let \(W_N, W_N'\) be Wigner matrices satisfying \textcolor{blue}{\textbf{W3}} whose off-diagonal entries have real and imaginary parts whose moments match up to third order.
%   Let \(c_0, c_1 > 0\)  be fixed and suppose that $Q$ satisfies
%   condition F. Then
%   \begin{align}
%     \label{eq:Q-swap-1}
%     \E Q(W_N) - \E Q(W_N')
%     \lesssim  \delta_{N}.
%   \end{align}
% \end{proposition}
\begin{proposition}
  \label{prop:multi-matching}
  Let \(W_N, W_N'\) be Wigner matrices
  % satisfying \textcolor{blue}{\textbf{W3}} whose off-diagonal entries
  % have real and imaginary parts
  whose moments match to third order.
  Let \(c_0, c_1> 0\) be fixed and for
each \(j = 1, \dotsc, m\), let \(Q_j \maps \C^{N \times N} \rightarrow
[0, 1]\) satisfy condition $F(\delta_{j,N})$.  If \(Q = \prod_{j=1}^m
Q_j\), then, 
  \begin{align}
    \label{eq:Q-swap}
    \E Q(W_N) - \E Q(W_N')
    \lesssim \max_{j=1, \ldots m} \delta_{j,N}.
%    \lesssim \sum_{j=1}^m \delta_{j,N}.
  \end{align}
\end{proposition}

\begin{proof}
  Consider first the case $m=1$. 
  We set $\Delta_{\gamma i}=Q(W^{(i)})-Q(W_0)$, and decompose
  \begin{equation*}
    \E Q(W_N) - \E Q(W_N')
       =  \sum_\gamma \E (\Delta_{\gamma 0} - \Delta_{\gamma 1}).
  \end{equation*}
Let $E_N = E_N(W_0^\gamma)$ denote the overwhelming probability event
\eqref{eq:D-cond-1} and then introduce `good' events
$G_{Ni} = E_N \cap \{ |\xi^{(i)}| \leq N^{c_0} \}$. 
Let $A$ be a fixed constant such that $N^{2-A} \lesssim \delta_N$.
Using boundedness of $Q$ and the moment bound \eqref{eq:W3}, with $p$
chosen so that $pc_0 > A$, we have
\begin{equation}  \label{eq:initial}
  \E (\Delta_{\gamma 0} - \Delta_{\gamma 1})
    = \E ( \Delta_{\gamma 0} \bo(G_{N0})) - \E ( \Delta_{\gamma 1}
    \bo(G_{N1}))  + O(N^{-A}).
\end{equation}

As before, set $Q_\gamma(t) = Q(W^\gamma_t)$, so that
$\Delta_{\gamma i} = Q_\gamma(\xi^{(i)}) - Q_\gamma(0)$.
By Taylor expansion,
\begin{equation*}
  \Delta_{\gamma i}
  = \sum_{j=1}^{k-1} \frac{1}{j!} Q_\gamma^{(j)}(0) (\xi^{(i)})^j
                 +\frac{1}{k!} Q_\gamma^{(k)}(\xi^*) (\xi^{(i)})^k,
\end{equation*}
for some $\xi^*$ with $|\xi^*| \leq |\xi^{(i)}|$.
Both $Q_\gamma(t)$ and event $E_N$ are independent of
$\xi^{(i)}$, so
\begin{align*}
  \E [ Q_\gamma^{(j)}(0) (\xi^{(i)})^j \bo(G_{Ni}) ]
   & = \E [Q_\gamma^{(j)}(0) \bo(E_N)] \, \E [ (\xi^{(i)})^j
     \bo( \abs{\xi^{(i)}} \leq N^{c_0}) ]  \\
   & = \E [Q_\gamma^{(j)}(0) \bo(E_N)] \, \E [ (\xi^{(i)})^j] +
     O(N^{-j/2} \delta_N \cdot N^{-A}).
\end{align*}
where we used the fact that
\(\E[\abs{\xi^{(i)}}^j\bo(\abs{\xi^{(i)}} > N^{c_0})]
\leq C_p N^{-c_0(p-j)} = O(N^{-A})\) for suitable $p$, as follows from
the Markov inequality and \eqref{eq:W3}.
For the remainder, on event $G_{Ni}$ we also have
$\abs{Q_\gamma^{(k)}(\xi^*)}
 \leq \| Q_\gamma^{(k)} \|_{c_0} \lesssim N^{-k/2} \delta_N$, and
 hence
 \begin{equation*}
   | \E [ Q_\gamma^{(k)}(\xi^*) (\xi^{(i)})^k \bo(G_{Ni}) ] |
     \lesssim N^{-k/2} \delta_N.
 \end{equation*}

Summarizing, we have 
\begin{equation*}
  \E[\Delta_{\gamma i}\bo(G_{Ni})]
   = \sum_{j=1}^{k-1} \frac{1}{j!} \E [Q_\gamma^{(j)}(0) \bo(E_N)] \,
   \E [ (\xi^{(i)})^j] + O(N^{-k/2} \delta_N +N^{-A}). 
\end{equation*}

Choose $k = k(\gamma)$ so that $\E (\xi^{(1)})^j = \E (\xi^{(0)})^j  $
for $1 \leq j \leq k-1$. Then the sums cancel and \eqref{eq:initial} yields
\begin{equation*}
  \E (\Delta_{\gamma 0} - \Delta_{\gamma 1})
    = O(N^{-k/2} \delta_N + N^{-A}).
\end{equation*}

For the $O(N^{2})$ off-diagonal terms, moment matching to third order allows
$k(\gamma) = 4$, while for the $N$ diagonal terms, we take
$k(\gamma) = 2$, since then only $\E \xi^{(i)} = 0$.
Summing over all $\gamma$, we obtain
\begin{equation} \label{eq:m-eq-1}
   \E Q(W_N) - \E Q(W_N')
    = O(\delta_N + N^{2 - A}) = O(\delta_N)
\end{equation}
from the choice of $A$.

For $m > 1$, apply the product rule, use \eqref{eq:D-cond-1}
and $\| Q_{j, \gamma} \|_{c_0} \leq 1$:
%  This follows simply from the product rule, wherein
  \begin{align*}
    % Q_\gamma^{(k)}(t)
    % &= \sum_{\ell_1 + \dotsb + \ell_m = k} \binom{k}{\ell_1, \dotsc, \ell_m} \prod_{j=1}^m Q_{j,\gamma}^{(\ell_j)}(t), \\
    \norm{Q_\gamma^{(k)}}_{c_0}
    &\lesssim \sum_{\ell_1 + \dotsb + \ell_m = k} \binom{k}{\ell_1,
      \dotsc, \ell_m} \prod_{\substack{1 \leq j \leq m \\ \ell_j \geq
    1}} N^{-\frac{\ell_j}{2}} \delta_{j,N}
    % \\
    %   &
        \lesssim N^{-\frac{k}{2}} \max_{j=1,...,m}\delta_{j,N},
  \end{align*}
  % where the last line follows from the fact that the dominating terms in the sum come from the terms with a single \(\delta_{j,N}\) factor.
  % But this is exactly the condition required for 
  % swapping with three matching moments, and so we conclude \cref{eq:Q-swap}.\qed
Thus $Q$ satisfies $F(\max_j \delta_{j,N})$ and the result follows
from \eqref{eq:m-eq-1}.
\end{proof}

%To avoid repetition later,
We use Proposition
\ref{prop:multi-matching} to formulate a criterion that allows joint
convergence in distribution of vector functions of $W_N$ to be
transferred to the 
corresponding functions of $W_N'$. A proof of the following proposition is in \cref{sec:proof-proposition-21}.

\begin{proposition}
  \label{prop:joint-convergence}
  Let \(W_N, W_N'\) be Wigner matrices
  % satisfying \textcolor{blue}{\textbf{W3}} whose
  % off-diagonal entries have real and imaginary parts
  whose moments match up to third order.
  Let $\bxi_N = \bxi_N(W_N)$ and $\bxi_N' = \bxi_N(W_N')$ both be
  $\mR^m$ valued random vectors.
  Suppose that $\bm{\xi}_N \stackrel{\rm d}{\to} \bm{\xi}$, and that
  each component $\xi_j$ of the limit has a continuous distribution
  function.
  
  Let $\eta_N \to 0$ be given, and suppose that for each $1 \leq j
  \leq m$ and $s \in \mR$ 
  %there exist functions $Q_j^\pm(\cdot,s)$   satisfying condition $F(\delta_{j,N})$ such that for $W = W_N, W_N'$, w.o.p.
  %\begin{alignat}{3}
  %\label{eq:Qjp}
  %  \bo \{ \xi_{Nj}(W) \leq s\}
  %& \leq Q_j^+(W,s)  &
  %& \leq \bo \{ \xi_{Nj}(W) \leq s + \eta_N\}  \\
 %   \bo \{ \xi_{Nj}(W) \leq s - \eta_N\}
 % & \leq Q_j^-(W,s)  &
 % & \leq \bo \{ \xi_{Nj}(W) \leq s \}
 %    \label{eq:Qjm}
 %  \end{alignat}
 there exists a function $Q_j(\cdot,s)$   satisfying condition $F(\delta_{j,N})$ such that for $W = W_N, W_N'$, w.o.p.
 \begin{equation}
     \label{eq:Qj}
      \bo \{ \xi_{Nj}(W) \leq s-\eta_N\}
   \leq Q_j(W,s)  
  \leq \bo \{ \xi_{Nj}(W) \leq s + \eta_N\}.
 \end{equation}
   
   Then we also have (joint) convergence
   $\bm{\xi}_N' \stackrel{\rm d}{\to} \bm{\xi}$.
\end{proposition}
%\begin{proof}
%  Let $\delta_N = \max_j \delta_{j,N}$.
%  It suffices to show that for each $\mathbf{s} = (s_j)$ that
%  \begin{equation}
%    \label{eq:sandwich}
%    \Pr( \bxi_N \leq \mathbf{s}-\eta_N) - O(\delta_N)
%    \leq \Pr( \bxi_N' \leq \mathbf{s})
%    \leq \Pr( \bxi_N \leq \mathbf{s}+\eta_N) + O(\delta_N).
%  \end{equation}
%Indeed, we then have
%\begin{equation*}
%  |\Pr( \bxi_N' \leq \mathbf{s}) - \Pr( \bxi_N \leq \mathbf{s}) |
%  \leq 
%  %\sum_j \mP ( |\xi_{nj}-s_j| \leq \eta_N )
%  \sum_j \Pr ( |\bxi_{Nj}-s_j| \leq \eta_N )
%  + O(\delta_N)
%  \to 0,
%\end{equation*}
%since each limiting component $\xi_j$ has a continuous distribution
%function.
%
%We verify the upper bound in \eqref{eq:sandwich}.
%For each $A > 0$ large, we have from \eqref{eq:Qjp} for $W_N'$, then
%Proposition  \ref{prop:multi-matching}
%and then \eqref{eq:Qjp} again, now for $W_N$,
%that
%\begin{align*}
%  \Pr( \bxi_N' \leq \mathbf{s})
%  & \leq \E \prod_j Q_j^+(W_N', s_j) + O(N^{-A}) \\
%  & \leq \E \prod_j Q_j^+(W_N, s_j) + O(\delta_N) 
%  \leq \Pr( \bxi_N \leq \mathbf{s}+\eta_N) +O(\delta_N).
%\end{align*}
%The lower bound in \eqref{eq:sandwich} follows similarly, using
%\eqref{eq:Qjm}.
%\end{proof}

\subsection{Flatness for Stieltjes functionals}
\label{sec:results-about-funct}

\subsubsection{Resolvent perturbation: deterministic bounds}
\label{sec:resolv-pert-determ}

We recall and modify some bounds of \cite{tao2012central} on stability
of Hermitian 
matrices with respect to perturbation in one or two entries, using
Ward's identity to improve the bounds at the edge.
% \textcolor{blue}{ Further, exploiting the Ward identity, e.g. \cite[][eq. (3.6)]{benaych2018}, we obtain better bounds than those in \cite{tao2012central} for $E$ near the edge. }
% \fix{Perhaps should emphasize what's different in our treatment.}

Let $M_0$ be a Hermitian $N \times N$ matrix, $z= E+ \im \eta \in
\C_+$ and $V$ an elementary matrix as defined after \eqref{eq:swapdefs}.
Set $M_t = M_0 + t N^{-1/2}V$ and $R_t = R_t(z) = (M_t-z)^{-1}$,
and $s_t(z) = N^{-1} \tr R_t(z)$.
Recall from \cite{tao2012central} the definitions of the matrix norms
$\|A\|_{(q,p)}$, and in particular
\begin{equation*}
  \|A\|_{(\infty,1)}
  =  \max_{1 \leq i,j \leq N} | A_{ij} |, \qquad
   \| A \|_{(\infty,2)}
  = \max_i \Big( \sum_j |A_{ij}|^2 \Big)^{1/2}.
\end{equation*}
Note also that if $V$ is an elementary matrix, then
\begin{align}
  \label{eq:tr-bd}
  |\tr(AV)|
  & = |\tr(VA)| \leq 2 \|A\|_\str \\
  \|AVB\|_\str
  & \leq 2 \|A\|_\str \|V\|_\str \|B\|_\str.
    \label{eq:mult-bd}
\end{align}

Let $\kappa_N(z,t) = |t| N^{-1/2} \| R_t \|_\str$.
Lemma 12 of \cite{tao2012central} says that if $\kappa_N(z,t) \to 0$ as $N\to\infty$, then for
large $N$
\begin{equation}
  \label{eq:neumann}
  R_{t+u} = R_t + \sum_{j=1}^\infty \Big(\frac{-u}{\sqrt N} \Big)^j
  (R_tV)^j R_t,
\end{equation}
with the right side being absolutely convergent. In addition, for $1
\leq p \leq \infty$,
\begin{equation}
  \label{eq:t-to-0}
  \|R_t\|_{(\infty,p)}
    \leq \|R_0\|_{(\infty,p)} \exp \{ 2|t|N^{-1/2} \|R_0\|_\str \}.
\end{equation}
Here the factor 2 arises from the use of
$\| V\|_{(1,\infty)} \leq 2$ in the proof of Lemma 12
  of \cite{tao2012central}.
% \eqref{eq:mult-bd} in the
% Tao-Vu argument.
The same bound holds with the roles of $R_0$ and $R_t$ reversed.

Expansion \eqref{eq:neumann} allows evaluation of $t$-derivatives of
$s_t(z)$. Indeed
\begin{align}
\partial_t^j s_t(z)
  & = j! N^{-j/2} c_j(z,t)  \label{eq:del-st}
    \\
  c_j(z,t)
   & = (-1)^j N^{-1} \tr((R_tV)^jR_t).  \label{coeff expl}
\end{align}
% We have
% \begin{align}
%   s_{t+u}(z)
%    & = s_t(z) + \sum_{j=1}^\infty \Big(\frac{u}{\sqrt N} \Big)^j
%      c_j(z,t)  \notag \\
%  \intertext{with}  \label{coeff expl}
%   c_j(z,t)
%    & = (-1)^j N^{-1} \tr((R_tV)^jR_t) \\
%  \intertext{and}   \partial_t^j s_t(z)
%     & = j! N^{-j/2} c_j(z,t)  \label{eq:del-st}
% \end{align}

The following variant of \cite[][Proposition 13]{tao2012central}
yields uniform bounds on 
$c_j$ in terms of $\| \Im R \|_\infty = \max_{1 \leq i \leq N} | \Im
R_{ii} |$, which allows tighter bounds near the edge.
% \textcolor{blue}{We modify their proof to constrain $c_j$ in terms of the imaginary part of the unperturbed resolvent's diagonal, which delivers tighter bounds for $E$ near the edge.}
\begin{proposition} \label{prop:TV12variant}
  Let $c_0$ and $A$ be small and positive.
%  $\check{\sigma}_N = (\log N)^{O(\log \log N)}$,
  and define
  \begin{align}
    \mathbf{S}_{e}(A)
    & = \{ z=E+ \im \eta \in \C:~ |E-2| \leq N^{-2/3+ A},
    % \check{\sigma}_N,
      \       \eta > N^{-2/3 -A} \}
      \label{eq:edge-dom}
    \\
    \kappa_N
    & = \sup_{|t|\leq N^{c_0}, z \in \mathbf{S}_e(A)}
              |t| \|R_0\|_\str / \sqrt N.  \notag
  \end{align}
  Then for $z \in \mathbf{S}_e(A)$ and $|t| \leq N^{c_0}$,
  \begin{equation}
        \abs{c_j(z,t)}
     \leq (N\eta)^{-1} 2^je^{2(j+1)\kappa_N} \norm{R_0}_{(\infty,1)}^{j-1} \norm{\Im R_0}_\infty.   \label{eq:cj-better}
  \end{equation}
%  where $\| \Im R \|_\infty = \max_{1 \leq i \leq N} | \Im R_{ii} |$.
\end{proposition}

\begin{proof}
% Before beginning the proof, some further remarks on matrix norms. 
  From the cyclic property of traces, then \eqref{eq:tr-bd} and
  \eqref{eq:mult-bd}, we have
    \begin{equation*}
    \abs{\tr ((R_tV)^{j}R_t)}
    = \abs{\tr (V(R_tV)^{j-1} R_t^2)}
    \leq 2 \norm{(R_tV)^{j-1} R_t^2}_\str
    \leq 2^{j} \norm{R_t}_{(\infty,1)}^{j-1} \norm{R_t^2}_{(\infty,1)}.
  \end{equation*}
The Ward identity, e.g. 
\cite[][eq. (3.6)]{benaych2018} says that
  \begin{equation*}
\sum_j|R_{ij}|^2=\eta^{-1}\Im R_{ii}
  \end{equation*}
is valid for any resolvent matrix $R=(W-E-\im \eta)^{-1}$ with $\eta\neq
0$ and Hermitian (or symmetric) $W$. For $\eta > 0$, we have
\begin{equation}
  \label{eq:Wcor}
  \| R \|_{(\infty,2)}^2
  = \max_i \sum_j |R_{ij}|^2
  = \eta^{-1} \| \Im R \|_\infty.
\end{equation}

If $B$ is a normal matrix, (i.e.\ $B^*B = BB^*$), then 
\begin{equation}
  \label{eq:csnormal}
  \|AB\|_{(\infty,1)} \leq \|A \|_{(\infty,2)} \|B \|_{(\infty,2)}.
\end{equation}
This uses the  Cauchy-Schwarz bound
$|(AB)_{ij}|^2 \leq  \sum_k |A_{ik}|^2 \sum_k |B_{kj}|^2$, since
$B$ normal implies
$\sum_k |B_{kj}|^2  = \sum_k |B_{jk}|^2 \leq \|B \|_{(\infty,2)}^2$ .
% \begin{equation*}
%   |(AB)_{ij}| \leq \Big( \sum_k |A_{ik}|^2 \sum_k |B_{kj}|^2 \Big)^{1/2},
% \end{equation*}
% and since $B$ is normal,
% \begin{equation*}
%   \sum_k |B_{kj}|^2
%   = (B^*B)_{jj} = (BB^*)_{jj}
%   = \sum_k |B_{jk}|^2
%   \leq \|B \|_{(\infty,2)}^2.
% \end{equation*}

  The resolvent of a Hermitian matrix is normal, so from
  \eqref{eq:csnormal},  \eqref{eq:t-to-0}, and then
  \eqref{eq:Wcor} we have
    \begin{equation*}
    \|R_t^2\|_{(\infty,1)}
    \leq \|R_t \|_{(\infty,2)}^2
    \leq e^{4 \kappa_N} \|R_0\|_{(\infty,2)}^2
    =  \eta^{-1} e^{4 \kappa_N} \| \Im R_0 \|_\infty .
  \end{equation*}
Combine the last display with the first of the proof and then
  refer to \eqref{eq:t-to-0} to bound
  $\|R_t \|_{(\infty,1)} \leq e^{2\kappa_N} \|R_0 \|_{(\infty,1)}$ to
  arrive at \eqref{eq:cj-better}.
\end{proof}

\subsubsection{Local law}
\label{sec:local-law-}

We will need the local law for Wigner matrices and some of its
important consequences, in particular at the spectral edge.

\begin{proposition}
  \label{prop:P14analog}
  Let $W_N$ be a Wigner matrix. \\
%  satisfying conditions \textbf{W1} and \textbf{W3}. \\
    (i) (local law)  Let $R(z) = (W_N-zI)^{-1}$
    denote the resolvent matrix and $s_{sc}(z)$ the Stieltjes
    transform of the semicircle law. Fix $\tau > 0$ small.
For each $\epsilon > 0$, we have
    w.o.p.
    \begin{align*}
    % &\frac{1}{N}\operatorname{tr}R(z)=s_{sc}(z)+O\left(\frac{N^{\varepsilon}}{N\eta}\right),\\
      &R_{ij} = s_{sc}(z) \delta_{ij} + O(N^\epsilon \Psi(z)),
    \end{align*}
uniformly for $z \in \mathbf{S}(\tau) = \{E + \im \eta:
|E| < \tau^{-1}, N^{-1+\tau} \leq \eta \leq \tau^{-1} \}$ and $i,j = 1,  \ldots , N$, where
\begin{equation*}
  \Psi(z) = \sqrt{ \frac{\Im s_{sc}(z)}{N\eta}} + \frac{1}{N\eta}.
\end{equation*}
    (ii) (semi-circle law on small scales) \ For each $\epsilon > 0$,
    we have w.o.p. that
    \begin{equation*}
      \mathcal{N}_{W_N}(I) = N \int_I \rho_{\rm sc}(\diff x)  + O(N^\epsilon),
    \end{equation*}
    uniformly for all intervals $I \subset \mR$, where $\mathcal{N}_{W_N}(I)$ denotes the number of eigenvalues of $W_N$ in $I$ and $\rho_{sc}$ denotes the semi-circle law. \\
%     (iii)  (eigenvalue rigidity) Let $\gamma_j$ be the typical location of the $j$-th eigenvalue $W_N$, which we denote as $\lambda_j$. That is
% \[
% N\int_{-\infty}^{\gamma_j}\rho_{sc}(\diff x)=N+1-j.
% \]
% Then, for any $\varepsilon>0$, w.o.p.
% \[
% |\lambda_j-\gamma_j|<N^{-2/3+\varepsilon}(j\wedge (N+1-j))^{-1/3}
% \]
% uniformly over $j=1,...,N$. \\
%     (iv) (complete eigenvector delocalization) \ For each $\epsilon > 0$,
%     we have w.o.p. that
%     \begin{equation*}
%       \sup_{1 \leq i \leq N} \| u_i(W_N) \|_\infty \leq N^{-1/2 + \epsilon},
%     \end{equation*} where $u_i(W_N)$ is a normalized eigenvector corresponding to $\lambda_i$ \textcolor{blue}{ and $\norm{\cdot}_{\infty}$ is a short hand notation for $\norm{\cdot}_{(\infty,1)}$}.\\
    % \textcolor{red}{Here again the concept of w.o.p. is applied to an
    %   asymptotic statement rather than an event, indexed by $N$, which
    %   is slightly confusing. \cite{benaych2018} avoid this by using
    %   $\prec$.}

    (iii) (at the edge.) Let $A > 0$ be small and fixed, and let
%    $\check{\sigma}_N = (\log N)^{O(\log \log N)}$, and let
    $\mathbf{S}_e(A)$ be the edge domain \eqref{eq:edge-dom}.
    For each $\epsilon > 0$ and uniformly for $z = E + \im \eta \in
    \mathbf{S}_e(A)$, we have w.o.p.
    \begin{equation}
      \label{eq:local-cor}
      \|R\|_{(\infty,1)} \lesssim  1 \wedge \eta^{-1}, \qquad
      \| \Im R \|_\infty
     \lesssim 
        (\eta^{1/2} + N^{-1/3 + \epsilon + A}) \wedge \eta^{-1}.  
      \end{equation}
        Let $W_0 = W - \xi N^{-1/2}V$ with $V$ an elementary matrix and
   $\xi$ satisfying moment bounds \eqref{eq:W3}. Set $R_0 = (W_0 - z
   I)^{-1}$. Then the bounds 
   \eqref{eq:local-cor} apply to $R_0$ also.
\
\end{proposition}
   \begin{remark} \ 
     For clarity, we emphasize that these are simultaneous high
     probability bounds for all $z$ in the indicated ranges. For
     example, 
     % $\mathbf{S}_e = \{ z \in \C :~ |E-2| \leq N^{-2/3}
     % \check{\sigma}_N, N^{-1/2} \leq \eta \leq 1 \}$,
     then w.o.p.
     \begin{equation*}
       \sup_{z \in \mathbf{S}_e(A)} \, (1 \vee \Im z)
       \| R(z) \|_{(\infty,1)}     \lesssim 1.
     \end{equation*}
   Such statements follow from the $N^2$-Lipschitz continuity of
   $R_{ij}(z), s_{\rm sc}(z)$ and of the right side bounds over the
   indicated ranges, c.f. e.g. \cite[][Remark 2.7]{benaych2018}.
   \end{remark}

\begin{proof}
  For (i) and (ii), see e.g. \cite[][Theorems 2.6, 2.8]{benaych2018}.
     We turn to (iii). 
     Basic bounds on $s_{sc}(z)$ , e.g. \cite[][Lemma 6.2]{erdos2017dynamical},
     establish for $\eta > 0$, $|E| \leq 10$ and $\kappa = | |E| - 2|$
     that 
  \begin{align*}
    |s_{sc}(z)| & \leq 1, \qquad
    \Im s_{sc}(z)  \lesssim \sqrt{\kappa + \eta}.
  \end{align*}
  For $N^{-2/3 - A} \leq \eta \leq 1$, we have
%  $(N\eta)^{-1} \leq N^{-1/3 +A}$. If also $\eta \leq 1$ we have
  $\Psi(z) \lesssim (N \eta)^{-1/2} \leq N^{-1/6 + A/2}$ and
so from the local law $\| R \|_{(\infty,1)} \lesssim 1$. 
For $\eta \geq 1$, just use the elementary bound
$|R_{jk}| \leq \eta^{-1}$ arising from the spectral decomposition
\begin{equation} \label{eq:spec}
  R_{jk}(E+\im \eta)
    = \sum_{l=1}^N \frac{u_l(j) u_l^*(k)}{\lambda_l-E-\im \eta},
\end{equation}
where $u_l(j)$ denotes the $j$-th component of the eigenvector $u_l$ corresponding to $\lambda_l(W_N)$.

For $(\Im R)_{jj}$ we exploit the improved bounds on $\Im s_{sc}$
at the edge.
Since $\kappa  \leq N^{-2/3+A}$,
%$\check{\sigma}_N$,
%For $N^{-2/3-A} \leq \eta  \textcolor{blue}{\leq 1}$,  we use
we have 
$\Im s_{sc} \lesssim N^{-1/3 + A/2}  + \eta^{1/2}$ and
$N \eta \geq N^{1/3- A}$, and  conclude
\begin{equation*}
  \Psi(z)
   \lesssim \frac{N^{-1/6+A/4}} {\sqrt{N\eta}} +
   \frac{1}{\sqrt{N \eta^{1/2}}} + \frac{1}{N\eta}
%   \leq N^{-1/3}( N^{A/2} \check{\sigma}^{1/4} + N^{A/4} + N^A)
   \lesssim N^{-1/3 + A}.   
 \end{equation*}
Hence, the second part of \eqref{eq:local-cor} follows from the local law.

Turning to $R_0$, we put $\Delta = W - W_0 = N^{-1/2} \xi V$ and use the
resolvent identity
% in the form
$R_0 = R + R\Delta R + R_0(\Delta R)^2$.
% \begin{equation*}
%   R_0 = R + R\Delta R + R_0(\Delta R)^2.
% \end{equation*}
Write $\| \cdot \|_*$ for $\| \cdot \|_{(\infty,1)}$.
Even for $R_0$, the bound $\|R_0 \|_* \leq \eta^{-1}$ follows from
\eqref{eq:spec} as before. 
So to conclude the rest of 
\eqref{eq:local-cor} for $R_0$, it suffices to show that w.o.p.
$\|R_0 - R \|_* \lesssim N^{-1/3 + \epsilon + A}$ for $ N^{-2/3 -
  A} \leq \eta \leq 1$.

We have the trivial bound $\|R_0 \|_* \leq \eta^{-1} \leq N^{2/3 +
  A}$.
Since \eqref{eq:W3} implies that $|\xi| \leq N^{\epsilon/2}$ w.o.p., 
% $\xi \prec 1$,
we have that $\| \Delta \| \lesssim N^{-1/2 +
  \epsilon/2}$, and along with $\| R \|_* \lesssim 1$, and
bound \eqref{eq:mult-bd} for elementary matrices, we
find that w.o.p. both
\begin{equation*}
  \| R \Delta R \|_* \lesssim N^{-1/2 + \epsilon/2}, \qquad
  \| R_0(\Delta R)^2 \|_* \lesssim N^{2/3 + A -1+ \epsilon}
  \lesssim N^{-1/3 + \epsilon + A}.  \qquad \qedhere
\end{equation*}
   \end{proof}

\subsubsection{Stieltjes functionals}
\label{sec:enum-funct}

We return to establishing flatness condition \eqref{eq:D-cond-1} for
certain functionals $Q = G \circ g$.
Let $W$ be an Hermitian matrix and $s_W(z)$ its empirical Stieltjes
transform. In the following proposition, we consider examples of \textit{Stieltjes functionals}
$g(W) = \Lambda(s_W)$ for some continuous linear functional
\(\Lambda\) acting on functions holomorphic on \(\C^+\).

The first two of these examples will be used in the next subsection to extend the non-concentration property for the eigenvalues of G(U/O)E matrices (Lemma~\ref{lem:anticoncGOE}) to Wigner matrices and, using this, to extend the log determinant CLT to Wigner matrices. The last two examples are key to the analysis of the SSK model in the companion paper \cite{johnstone2021}. There, we need to extend results on the $k$-th largest eigenvalue and the trace of the inverse powers of $z-W_N$ from G(U/O)E to Wigner matrices.

\begin{proposition}
  \label{prop:system}
  Let \(W\) be a Wigner matrix.
  % satisfying \textbf{W1} and \textbf{W3}.
  Let $\epsilon > 0$ and $0 < c_0 < \frac{1}{2}$ and
let \(E \in \R\) be such that $\abs{E - 2} \lesssim
N^{-2/3 +A}.$
%\check{\sigma}_N N^{-\frac{2}{3}}\) and let \(\epsilon > 0\), \(0 < c_0 < 1/2\).

  For each of the following statistics, define functions \(g \maps
  \C^{N\times N} \rightarrow \R\), \(G\maps \R \rightarrow \R\) and a
  sequence \(\delta_N\) according to the following specifications,
  in each case for $1 \leq j \leq 4$:
  \begin{enumerate}
  \item Log-determinant: with $\gamma_N = N^{-2/3-\epsilon}$,
    \begin{equation*}
      g(W) = N \int_{\gamma_N}^{N^{100}} \Im s_W(E+ \im \eta) \diff \eta,
      \qquad \|G^{(j)}\|_\infty \leq (\log N)^{-j/4},
      \qquad \delta_N = (\log N)^{-1/4}.
    \end{equation*}
  \item Eigenvalue counting: with $\eta = N^{-2/3-9\epsilon}$ and
    $E_i, i = 1,2$ such that $|E_i - 2| \leq N^{-2/3 + 10 \epsilon}$,
        \begin{equation*}
      g(W) = \frac{N}{\pi} \int_{E_1}^{E_2} \Im s_W(x + \im \eta) \diff x,
      \qquad \|G^{(j)}\|_\infty \leq (\log N)^{Cj},
      \qquad \delta_N = N^{-\frac{1}{3} + O(\epsilon)}.
    \end{equation*}
  % \item Non-concentration: with $\eta = N^{-2/3-9\epsilon}, \ell =
  %   \tfrac{1}{2} N^{-2/3-\epsilon}$, and $E^\pm = E \pm N^{-2/3 - c_0}$ for \(\epsilon = 2c_0\),
  %   \begin{equation*}
  %     g(W) = N \int_{E^- + \ell}^{E^+ - \ell} \Im s_W(y+ \im \eta) \diff y,
  %     \qquad \|G^{(j)}\|_\infty \leq C^{j},
  %     \qquad \delta_N = N^{-\frac{1}{3} + O(\epsilon)}.
  %   \end{equation*}
  % \item \(k\)th-largest eigenvalue: with \(s \in \R\), $\eta = N^{-2/3-9\epsilon}, E_1 =
  %   2+N^{-2/3}(s+N^{-\epsilon}/2)$ and $E_2 = E_+ = 2 +
  %   2N^{-2/3+\epsilon}$,
  %   \begin{equation*}
  %     g(W) = N \int_{E_1}^{E_2} \Im s_W(y + \im \eta) \diff y,
  %     \qquad \|G^{(j)}\|_\infty \leq (\log N)^{Cj},
  %     \qquad \delta_N = N^{-\frac{1}{3} + O(\epsilon)}.
  %   \end{equation*}
  \item Inverse moments:
    with \(\eta = N^{-2/3-\epsilon}\) and $l\in\mathbb{Z}_+$,
    \begin{equation*}
      g(W)
      = N^{-\frac{2}{3}l + 1} \Re s_W^{(l-1)}(E + \im \eta),
      \qquad \|G^{(j)}\|_\infty \leq (\log N)^{Cj}, \\
      \qquad \delta_N = N^{-\frac{1}{3} + O(\epsilon)}.
    \end{equation*}
  \end{enumerate}

  In each of the cases listed above, the corresponding function \(Q = G \circ g\) satisfies the condition of \cref{eq:D-cond-1}.
  That is, for \(1 \leq k \leq 4\), it follows w.o.p.\ that
  \begin{align*}
    \norm{Q_\gamma^{(k)}}_{c_0}
    \lesssim N^{-\frac{k}{2}} \delta_N.
  \end{align*}

  \proof
  Define $g^\gamma(t) = g(W_t^\gamma)$ so that
  $Q_\gamma(t) = G(g_\gamma(t))$.
  In order to bound $Q_\gamma^{(k)}(t)$ we start with bounds for
  $\partial_t^j g^\gamma$.
  Recalling \eqref{eq:t-defs0}, we have $g^\gamma(t) =
  \Lambda(s_t^\gamma)$. Standard results on differentiation of
  integrals and then \eqref{eq:del-st} imply that
  % In all cases above, we have
  % \begin{align*}
  %   g^\gamma(t)
  %   &= \Lambda(s_t^\gamma)
  % \end{align*}
  % for some continuous linear functional \(\Lambda\) acting on functions holomorphic on the upper half-plane \(\C^+\).
  % Hence, it follows from \cref{eq:del-st} that
  \begin{align*}
    \partial_t^j g^\gamma(t)
    = \Lambda(\partial_t^j s_t^\gamma) 
    = j! N^{-j/2} \Lambda(c^\gamma_j(\adot, t)), 
    % \norm{\partial_t^j g^\gamma}_{c_0}
    % &\lesssim N^{-\frac{j}{2}} \sup_{\abs{t} \leq N^{c_0}}
    %   \abs{\Lambda(c^\gamma_j(\adot, t))} ,
  \end{align*}
  where from \eqref{coeff expl} and \eqref{eq:t-defs0})
  % \textcolor{blue}{where  $c_j^\gamma(z,t)$ is as defined in \eqref{coeff expl}, but with $V$ and $R_t$ replaced by $V_\gamma$ and $R_t^\gamma$, respectively (for the definition of the latter see \eqref{eq:t-defs0}). That is, 
  \[
 c_j^\gamma(z,t) =(-1)^jN^{-1}\tr((R_t^\gamma V_\gamma)^jR_t^\gamma).
  \]

  Hence, to bound \(\norm{\partial_t^j g^\gamma}\), it suffices to use bounds on the coefficients \(c_j^\gamma\). We will omit the superscript $\gamma$ to simplify notations.
  From Propositions \ref{prop:TV12variant} and \ref{prop:P14analog},
%  and Corollary  \ref{cor:local-cor},
  for fixed \(A > 0\), 
  \begin{equation}
    \label{eq:c-bounds}
    N \abs{c_j(E + \im \eta,t)}
    \lesssim
    \begin{cases}
      \eta^{-1/2}+ N^{-1/3 + \epsilon + A}\eta^{-1}  & N^{-2/3 - A} \leq \eta \leq 1 \\
      %\eta^{-1/2}               & N^{-1/2} \leq \eta \leq 1 \\
      \eta^{-j-1}              & \eta \geq 1,
    \end{cases}
  \end{equation}
uniformly in $|t|\leq N^{c_0}$. Note that there is no dependence on $j$ for $\eta \leq 1$.

  In the log-determinant case, we have that
  \begin{align*}
    \Lambda(f)
    &= \int_{\gamma_N}^{N^{100}} N \Im f(y + \im \eta) \diff \eta.
  \end{align*}
  Evaluated at \(c_j(\adot, t)\), we use \eqref{eq:c-bounds}
  with $A = 2 \epsilon$ to obtain
%  bound the absolute value above by a constant times
  \begin{align*}
    \int_{\gamma_N}^{N^{100}} N\abs{c_j(E + \im \eta, t)} \diff \eta
    &\lesssim \int_{\gamma_N}^{1} (\eta^{-1/2} + N^{-1/3 + \epsilon +
      A}  \eta^{-1}) \diff \eta 
    % + \int_{\textcolor{red}{\gamma_N}}^1 \eta^{-1/2} \diff \eta
    + \int_1^{N^{100}} \eta^{-j-1} \diff \eta \\
    &\lesssim 1.
  \end{align*}

  For the remaining integrals, we need only the following consequence
  of 
 %  \(N^{-\frac{2}{3} - A}  \leq \eta \leq 1\) case
  of \cref{eq:c-bounds}. 
%  In particular, we use the consequence that for these values of \(\eta\)
  \begin{equation}
    \label{eq:c-weak-bound}
    N\abs{c_j(E + \im \eta, t)}
    \lesssim N^{\frac{1}{3} + \epsilon + 2A}.
  \end{equation}
%  We use \cref{eq:c-weak-bound} with
  Set \(A = 10 \epsilon\) for the eigenvalue counting case.
%  both the \(k\)th-largest eigenvalue and non-concentration cases.
  This yields
  \begin{align*}
    \int_{E_1}^{E_2}\frac{N}{\pi} \abs{c_j(y + \im\eta, t)} \diff y
    &\lesssim N^{\frac{1}{3} + \epsilon + 2A} N^{-\frac{2}{3} + \epsilon} 
    = N^{-\frac{1}{3} + O(\epsilon)}.
  \end{align*}
  % for the $k$th largest eigenvalue.
  % Similarly we reach the same final bound for non-concentration since
  % $E_+ - E_- = 2 N^{-2/3 - \epsilon/2}$.
  % \begin{align*}
  %   \int_{E_- + \ell}^{E_+ - \ell}N\abs{c_j(y + \im\eta, t)} \diff y
  %   &\lesssim N^{\frac{1}{3} + \epsilon + 2A} N^{-\frac{2}{3} - \epsilon} \\
  %   &\lesssim N^{-\frac{1}{3} + O(\epsilon)}.
  % \end{align*}

  For inverse moments, we have
$\Lambda(c_j(\cdot,t)) = N^{-2l/3+1} \Re c_j^{(l-1)}(E+ \im \eta)$. 
  % by Cauchy's integral formula that
  % \begin{align*}
  %   \Lambda(f)
  %   &= N^{-\frac{2}{3}l + 1} \Re f^{(l-1)}(E + \im \eta) \\
  %   &= \Re \Bigl[\frac{(l-1)!}{2\pi \im} \oint_\Gamma \frac{N f(w)}{N^{\frac{2}{3}l}(w - E - \im \eta)^l} \diff w\Bigr],
  % \end{align*}
Let \(\Gamma\) be a contour of radius \(N^{-\frac{2}{3} - 2\epsilon}\) around \(E + \im \eta\).
  In this way, each \(c_j\) is analytic on the interior of \(\Gamma\),
  and so we use Cauchy's integral formula and
  \eqref{eq:c-weak-bound} with \(A = 2\epsilon\) to see  that
  \begin{align*}
    N^{-\frac{2}{3}l + 1}\abs{c_j^{(l-1)}(E + \im \eta)}
    &\leq \frac{(l-1)!}{2\pi} \oint_{\Gamma} \frac{N\abs{c_j(w)}}{N^{\frac{2}{3} l}|w - E - \im \eta|^l} \abs{\diff w} 
    \lesssim N^{-\frac{1}{3} + O(\epsilon)}.
  \end{align*}

 In sum suppose that, for sequences \(a_N\) and \(b_N\) such
 that \(a_Nb_N  \rightarrow 0\), we have w.o.p. 
  \begin{align*}
    \norm{\partial_t^j g^\gamma(t)}_{c_0}
    &\lesssim N^{-\frac{j}{2}} a_N, \qquad 
    \norm{G^{(j)}}_\infty
    \lesssim b_N^j
  \end{align*}
In the proof so far, we have seen that the above conditions hold with the following values of \(a_N\) and \(b_N\) for some constant \(C\):
  \begin{enumerate}
  \item Log-determinant: \qquad \quad \(a_N = 1, \qquad \qquad \ \ \ 
    b_N = (\log N)^{-1/4}\).
  \item Eigenvalue counting: \quad \ \ \(a_N = N^{-\frac{1}{3} +
      O(\epsilon)}, \quad b_N = (\log N)^C\).
  % \item Non-concentration: \qquad \ \(a_N = N^{-\frac{1}{3} +
  %     O(\epsilon)},  \quad b_N = C\).
  % \item \(k\)th-largest eigenvalue: \quad \(a_N = N^{-\frac{1}{3} +
  %     O(\epsilon)}, \quad b_N = (\log N)^C\).
  \item Inverse moments: \qquad \quad \(a_N = N^{-\frac{1}{3} +
      O(\epsilon)},  \quad  b_N = (\log N)^C\).
  \end{enumerate}

  We apply Fa\`{a} di Bruno's formula to compute bounds for
  \(\partial_t^k (G \circ g^\gamma)(t)\). 
  Let \(\mathcal{M}_k = \{m \in \Z_{\geq 0}^k : \sum_{j=1}^k j
  m_j = k\}\), so that $m_+ = m_1 + \dotsb + m_k \geq 1$ for each \(m
  \in   \mathcal{M}_k\).
  Then for certain combinatorial constants \(C_{km}\) we have that,
  uniformly in  \(\abs{t} \leq N^{c_0}\),  
  \begin{align*}
    \abs{\partial_t^k (G \circ g^\gamma)(t)}
    &\leq \sum_{m \in \mathcal{M}_k} C_{km} \abs{G^{(m_+)}(g^\gamma(t))} \cdot \prod_{j=1}^k \abs{g^{(j)}(t)^{m_j}} \\
    &\lesssim \sum_{m\in \mathcal{M}_k} C_{km} b_N^{m_+} \prod_{j=1}^k N^{-\frac{j m_j}{2}} a_N^{m_j} 
    = N^{-\frac{k}{2}}  \sum_{m\in \mathcal{M}_k} C_{km} (a_N b_N)^{m_+} 
    \lesssim N^{-\frac{k}{2}} a_N b_N,
  \end{align*}
%  where the last inequality  uses \(m_+ \geq 1\)
  Hence, the conclusion in each case follows with
  \(\delta_N = a_Nb_N\). \qed 
\end{proposition}

\subsection{Concluding the extension to Wigner matrices}
\label{sec:dist-conclusions}

\subsubsection{Wigner Non-concentration}
\label{sec:wign-non-conc}

\begin{proposition}  \label{prop:wign-non-conc-2}
  Let $W_N'$ be a Wigner matrix
  % satisfying \textbf{W1-3}
  whose off-diagonal
  % entries have real and imaginary parts whose
  moments match GOE or GUE to third order.
  Call its eigenvalues \(\lambda_1', \dotsc, \lambda_N'\).
  Let \(E \in \R\) be such that \(\abs{E - 2} \lesssim
  N^{-\frac{2}{3}} \check{\sigma}_N\), with
  $\check{\sigma}_N = (\log N)^{O(\log \log N)}$.
Then there exists a \(c_1\) such that, for each \(c_0 \in  (0, c_1)\), there exists \(d > 0\) such that, for \(N\) large,
  \begin{align} \label{eq:nonconcW}
    \Pr(\min_{j = 1, \dotsc, N} \abs{\lambda_j' - E} \leq N^{-\frac{2}{3} - c_0})
    &\leq N^{-d}.
  \end{align}
\end{proposition}
\begin{proof}
  Define the eigenvalue counting function $\mathcal{N}_W(E_1,E_2)  =
  \# \{j: E_1 \leq \lambda_j(W) \leq  E_2 \}$.
  The event in \eqref{eq:nonconcW} has the form
  $\mathcal{N}_W(E_1,E_2)  \geq 1$.
  The first step is to approximate this using the Stieltjes transform.

  % Let \(W_N\) be a GOE or GUE matrix and let \(W_N'\) be a Wigner matrix satisfying \textbf{W1-3} whose off-diagonal entries have real and imaginary parts whose moments match those of \(W_N\) to third order.
  % Denote by \(\lambda_1, \dotsc, \lambda_N\) and \(\lambda_1', \dotsc, \lambda_N'\) the eigenvalues of \(W_N\) and \(W_N'\) respectively.
  % Let \(c_0 > 0\).
  
  Let \(\epsilon = 2c_0\) and define \(\ell =
  \frac{1}{2}N^{-\frac{2}{3} - \epsilon}\), and \(\eta = N^{-\frac{2}{3} -
    9\epsilon}\). 
  Let \(E_1, E_2 \in \R\) be such that \(\abs{E_1 - 2}, \abs{E_2 - 2} \lesssim N^{-\frac{2}{3}} \check{\sigma}_N\) and \(E_2 - E_1 \geq 2\ell\).

  A suitable approximation is given by Corollary
  17.3 of \cite{erdos2017dynamical} (based on the local law and
  eigenvalue rigidity)\footnote{The definition of Wigner matrices in \cite{erdos2017dynamical} is slightly different from ours, so its proof needs minor adjustments to take the difference into account. Specifically, one needs to use suitable versions of the local law and eigenvalue rigidity (theorems 2.6 and 2.9 from \cite{benaych2018}). We refer the interested reader to section A4 of \cite{johnstone2021} for details.}, which we apply twice
  % For the Wigner matrix \(W = W_N\) or \(W_N'\),
  % define $\mathcal{N}_W(E_1,E_2)  = \# \{j: E_1 \leq \lambda_j(W) \leq
  % E_2 \}$. We apply  twice,
  with \(E = E_1\) and \(E_2\) respectively.
  Subtracting the latter bounds from the former, this yields w.o.p.\ that
  \begin{equation*}
%    \label{eq:2-side-bounds}
    \frac{N}{\pi} \int_{E_1 + \ell}^{E_2 - \ell} \Im s_W(y + \im \eta)
    \diff y - 2 N^{-\epsilon}
    \leq \mathcal{N}_W(E_1, E_2)
    \leq \frac{N}{\pi} \int_{E_1 - \ell}^{E_2 +
      \ell} \Im s_W(y + \im \eta) \diff y + 2 N^{-\epsilon}. 
  \end{equation*}

  Let \(E^\pm = E \pm 2 N^{-\frac{2}{3} - c_0}\), and
  define the function
  \begin{align*}
    g(W)
    &= \frac{N}{\pi} \int_{E^{-} + \ell}^{E^+ - \ell} \Im s_W(y + \im
      \eta) \diff y, 
  \end{align*}
  Applying these bounds with \((E_1, E_2) = (E^-, E^+)\) and \((E^- + 2\ell, E^+ - 2\ell)\), we conclude that, w.o.p.,
  \begin{gather}
    \label{eq:count-bounds}
    \mathcal{N}_W(E^{-} + 2\ell, E^{+} - 2\ell) - 2 N^{-\epsilon}
 %   \leq N\int_{E^{-} + \ell}^{E^{+} - \ell} \Im s_W(y + \im \eta) \diff y
    \leq g(W)
    \leq \mathcal{N}_W(E^{-}, E^{+}) + 2 N^{-\epsilon}.
  \end{gather}

  Let \(G\) be a smooth increasing function such that
  \begin{align*}
    G(x)
    = \begin{cases}
      1 & \text{if } x \geq 2/3, \\
      0 & \text{if } x \leq 1/3.
    \end{cases}
  \end{align*}

  Taking \(Q = G \circ g\) and applying \(G\) to each side of
  \cref{eq:count-bounds}, we then have that, w.o.p., 
  \begin{align*}
    \Ind\{\mathcal{N}_W(E^{-} + 2\ell, E^{+} - 2\ell) \geq 1\}    
    \leq Q(W)
    \leq \Ind\{\mathcal{N}_W(E^{-}, E^{+}) \geq 1 \}.
  \end{align*}

  Now we can use Propositions \ref{prop:multi-matching} and
  \ref{prop:system}(2) to compare $Q(W_N')$ with $Q(W_N)$, for $W_N$
  drawn from G(O/U)E with egienvalues $\lambda_j$. 
  % Using the above, together with \cref{prop:system}(3) we have,
  For  any \(A > 0\), we have
  \begin{align*}
    \Pr(\min_{j} \abs{\lambda_j' - E}
     \leq 2N^{-\frac{2}{3} - c_0} -
    2\ell)  
    & = \Pr \{ \mathcal{N}_{W'_N} (E^{-} + 2\ell, E^{+} - 2\ell) \geq
      1\} \\
    &\leq \E Q(W_N') + O(N^{-A}) \\
    &\leq \E Q(W_N) + O(N^{-\frac{1}{3} + O(\epsilon)}) \\
    &\leq \Pr(\min_{j} \abs{\lambda_j - E} \leq 2 N^{-\frac{2}{3} -
      c_0}) + O(N^{-\frac{1}{3} + O(\epsilon)}) \\
%      \numberthis{eq:gaussian-non-conc-bound} \\
    &  \leq \tfrac{1}{2}N^{-d} + O(N^{-1/3 + O(\epsilon)})
       \leq N^{-d}
  \end{align*}
At the last line we applied the non-concentration bound for G(O/U)E,
Section \ref{sec:gauss-non-conc}.

For $N$ large, we have $2N^{-\frac{2}{3} - c_0} - 2\ell \geq
N^{-2/3-c_0}$ and so the final bound \eqref{eq:nonconcW} follows from
these inequalities.
\end{proof}

The next lemma shows that non-concentration implies control of
inverse power sums at around their typical magnitude. The proof is by
standard dyadic decomposition (see \cref{sec:proof-lemma-28}).

\begin{lemma}
  \label{lem:neg-mom-bd}  
Let $\{ \lambda_j \}$ be the eigenvalues of a Wigner matrix $W_N$ 
whose off-diagonal moments match GOE or GUE to third order.
%%  satisfying  conditions \textbf{W1-3}.
  Suppose that $|E - 2| \leq N^{-2/3}  \check{\sigma}_N$.
  % Let $\delta = N^{-2/3 - \zeta}$ for some $\zeta \in (0,\frac{1}{3})$, so
  % that the event $A_N = \{ \min_j |E - \lambda_j| > \delta \}$ occurs with high
  % probability.
  Then there exist constants $\{C_r\}$ such that for each $\epsilon >
  0$ small, with high probability we have
  % on the event $A_N$ we have for each $0<\epsilon<1/3-\zeta$ w.o.p.
  \begin{equation*}
    S_r(E) :=
    \sum_{j=1}^N \frac{1}{|\lambda_j-E|^r}
    \leq
    \begin{cases}
      C_1 N \quad & \text{if } \ r = 1 \\
      C_r N^{2r/3 + (r +1) \epsilon} & \text{if } \ r \geq 2.
    \end{cases}
  \end{equation*}
  The bounds also hold for $S_r(E')$ uniformly in $|E'-E| \leq
\delta/2$ with $\delta=N^{-2/3-\epsilon}$, by increasing $C_r$ to $2^r C_r$.
\end{lemma}

\subsubsection{Log-determinant}
\label{sec:log-determ-larg}

We  derive the central limit theorem for the log-determinant 
\begin{align*}
  L_N(W_N')
  &= \log \abs{\det(W_N' - E)}
\end{align*}
for a Wigner matrix \(W_N'\) and $E = E_N = 2 + \sigma_N N^{-2/3}$.
%with $-\gamma \leq \sigma \ll \log^2 N$.
% with $\ts_N$
% a monotone sequence for which $\ts_N \geq -C$ for some positive constant $C$ and
% $\ts_N = o(\log^2 N)$. 
Recall the scaling constants $\mu_N, \tau_N$ from \eqref{eq:centscal1}.
% \begin{equation*}
%   \mu_N = \tfrac{1}{2} N + \ts_N N^{1/3} - \tfrac{2}{3} \ts_N^{3/2} -
%   \tfrac{1}{6}(\alpha-1)\log N, \qquad
%   \tau_N = (\tfrac{\alpha}{3} \log N)^{1/2}.
% \end{equation*}
Let $W_N$ be drawn from (scaled) GOE or GUE.
From Theorem \ref{th:main}, which we have already established for the Gaussian ensembles, we have
%\cite[][Theorem 1]{johnstone2020logarithmic}, we have
\begin{equation}
  \label{eq:W-clt}
   \check{L}_N(W_N) = \tau_N^{-1}(L_N(W_N)-\mu_N) \stackrel{\rm d}{\to} \mathcal{N}(0,1).
\end{equation}

\begin{proposition}[Log determinant CLT] \label{prop:zero-diag}
  Let $W_N'$ be a Wigner matrix 
whose off-diagonal moments match GOE or GUE to third order.
% satisfying \textbf{W1-4}.
Let $E = E_N = 2 + \sigma_N N^{-2/3}$ with
  with $-\gamma \leq \sigma_N \ll \log^2 N$.
%   $\ts_N$
% a monotone sequence for which $\ts_N \geq -C$ for some positive constant $C$ and
% $\ts_N = o(\log^2 N)$.

  Then
  \begin{equation}
  \label{eq:zW-clt}
  \tau_N^{-1} (\log\abs{\det(W_N' - E)} - \mu_N)  \stackrel{\mathrm{d}}{\rightarrow} \mathcal{N}(0,1).
\end{equation}

\proof

To rewrite the log-determinant in terms of an integral of the
Stieltjes transform, note that
$(\diff/\diff \eta) \log |\lambda - E - \im \eta|
= \Im [(\lambda-E- \im \eta)^{-1}]$, which yields
\cite[][eq.\ (46)]{tao2012central}
\begin{equation*}
  L_N(W)
  = \log \abs{\det(W-E - \im N^{100})} - N\int_0^{N^{100}} \Im s_W(E + \im \eta) \diff \eta.
\end{equation*}

% Arguing as in \cite[p.\ 93]{tao2012central}, we have for any Wigner
% matrix \(W\) satisfying \textbf{W1-3}
The uniform moment bounds \eqref{eq:W3} imply that
\begin{align*}
  \log \abs{\det(W-E - \im N^{100})}
  &= N \log(N^{100}) + O_{\Pr}(N^{-50}).
\end{align*}

Moreover, for each \(\epsilon > 0\) small, if we take \(\gamma_N =
N^{-\frac{2}{3} - 2\epsilon}\), then non-concentration implies that
the contribution to the integral from $\eta \leq \gamma_N$ 
is negligible.
Indeed, with $\lambda_j = \lambda_j(W_N)$, 
\begin{equation*}
  N \Im s_{W}(E+ \im \eta)
%    = \sum_{j=1}^N \Im \frac{1}{\lambda_j-E- \im \eta}
    = \eta \sum_{j=1}^N \frac{1}{(\lambda_j-E)^2+\eta^2}
    \leq \eta \sum_{j=1}^N \frac{1}{(\lambda_j-E)^2}.
\end{equation*}
By \cref{lem:neg-mom-bd}
% with $\zeta = \epsilon$
we thus have
\begin{align*}
  \abs[\Big]{N\int_{0}^{\gamma_N} \Im s_W(E + \im \eta) \diff \eta}
%   &\leq \int_0^{\gamma_N} \eta \sum_{j=1}^N
%   \frac{1}{(\lambda_j-E)^2} \diff \eta
   & \leq \tfrac{1}{2} \gamma_N^2 S_2(E)    
  =  O_{\Pr}(\gamma_N^2 N^{\frac{4}{3} + 3\epsilon}) 
   = o_{\Pr}(1).
\end{align*}

To summarize, if we define the Stieltjes functional
\begin{equation*}
    g(W)
  = N \int_{\gamma_N}^{N^{100}} \Im s_W(E + \im \eta) \diff \eta, 
\end{equation*}
set \(\bar{\mu}_N = \mu_N + N \log(N^{100})\) and define
\begin{equation*}
  \xi_N(W) = \tau_N^{-1} (g(W_N)-\bar{\mu}_N),
\end{equation*}
then we have shown that $\check{L}_N(W_N)  = \xi_N(W_N)  + o_\Pr(1)$. 

We carry out the Lindeberg swapping with $g(W)$.
Let $H: \mR \to [0,1]$ be a smooth decreasing function such that 
\begin{equation*}
  H(x) =
  \begin{cases}
    1 & \text{if } \ x \leq -\eta_N \\
    0 & \text{if } \ x \geq \eta_N.
  \end{cases}
\end{equation*}
%and let $H^-(x) = H^+(x-\eta_N)$. 
For $s \in \mR$ define $G_s(x) = H(\tau_N^{-1}(x-\bar
\mu_N)-s)$.
One verifies that
\begin{equation}
  \label{eq:logdetsw}
  \bo \{\xi_N(W) \leq s-\eta_N \} 
  \leq G_s(g(W))
  \leq \bo \{\xi_N(W) \leq s + \eta_N\} .
\end{equation}
Setting $Q(W,s) = G_s(g(W))$, we obtain bound \eqref{eq:Qj}.
%Bound \eqref{eq:Qjm} follows similarly using instead $H^-$ with
%$Q^-(W,s) = G_s^-(g(W))$.

Observe that $\| G_s^{(j)} \|_\infty \lesssim (\tau_N \eta_N)^{-j}
\lesssim (\log N)^{-j/4}$ if we choose $\eta_N = \tau_N^{-1/2}$.
Then Proposition \ref{prop:system} (1) implies that
$Q(\cdot,s)$ satisfy condition F with $\delta_N = (\log N)^{-1/4}$. 
From Proposition \ref{prop:joint-convergence} we conclude that
$\xi_N(W_N')$ and hence $\check{L}(W_N')$ have the same limiting
distribution as $\xi_N(W_N)$ and $\check{L}_N(W_N)$. 
Thus the validity of Theorem \ref{th:main} for Gaussian ensembles implies \cref{eq:zW-clt}. \qedsymbol

\end{proposition}

\section{The spiked case}
\label{section spiked}

In this section, we consider deformed Wigner matrices
\[
W_{h,N}=W_N+h\mathbf{vv^\ast},
\]
where $\mathbf{v}$ is an arbitrary deterministic vector from
$\mathbb{R}^N$ ($\alpha=2$) or $\mathbb{C}^N$ ($\alpha=1$) with unit
norm, $\norm{\mathbf{v}}=1$ and $h\neq 1$ is a fixed non-critical spike.
%and $W_N$ satisfies \textbf{W1-4}.
Since $\mathbf{v}$ is arbitrary, a version of the isotropic local law
of Knowles and Yin \cite{knowles2013isotropic} plays a key role.

%\textit{Isotropic local law.} 
\begin{proposition}
\label{thm: isotropic}
Let $W_N$ be a Wigner matrix
whose off-diagonal moments match GOE or GUE to third order.
% satisfying conditions \textbf{W1}, \textbf{W2}, and \=textbf{W3}. In
% addition, suppose that the off-diagonal entries of $W_N$ have zero
% third moment in the sense that
% $\E\xi_{ij}^3=\E\xi_{ij}^2\overline{\xi_{ij}}=0$ for any $i\neq j$.
We have:
\newline
(i) (isotropic local law) Fix $\tau>0$. Then for each $\epsilon > 0$, we have
    w.o.p.
\begin{equation}
    \label{isotropic law}
    \mathbf{v}^\ast R(z)\mathbf{w}=s_{sc}(z)\mathbf{v}^\ast\mathbf{w}+O (N^{\varepsilon}\Psi(z))
\end{equation}
uniformly for $z \in \mathbf{S}(\tau)$ and
% normalized (to have unit Euclidean norm)
deterministic vectors $\mathbf{v,w}$ of unit Euclidean
length in $\mathbb{C}^N$.
% Outside the spectrum, we have w.o.p. 
% \begin{equation}
%   \label{eq:outside-spec}
%   \mathbf{v}^\ast R(z)\mathbf{w}=s_{sc}(z)\mathbf{v}^\ast\mathbf{w}
%   +O(N^{\varepsilon-1/2}(\kappa + \eta)^{-1/4})
% \end{equation}
% uniformly for $z \in \mathbf{S}^\circ(\tau)\equiv \{ z \in \mathbb{C}: 2 \leq |E| \leq
%   \tau^{-1}, \   \eta > 0 \}$ and
%   normalized $\bv, \mathbf{w}$, where $\kappa=\big||E|-2\big|$.
\newline
(ii) (isotropic delocalization) Let $\mathbf{u}^{(j)}$ be the $j$-th principal normalized eigenvector of $W_N$. Then, for each $\varepsilon>0$, we have w.o.p
\begin{equation}
    \label{isotropic delocalization}
    \max_{j}|\mathbf{v}^\ast \mathbf{u}^{(j)}|^2=O(N^{\varepsilon-1})
\end{equation}
uniformly for normalized deterministic vectors $\mathbf{v}\in\mathbb{C}^N$.
\end{proposition}

% A proof of this proposition is given in the Appendix (\cref{sec:
%   appendix}). It
The proof is a direct modification, summarized in \cref{sec:B3}, of those of
Case A of Theorems 2.2 and 2.5 of \cite{knowles2013isotropic},
% where similar results are 
established for Wigner matrices with sub-exponential entries whose
third moments match those of GUE/GOE.
In contrast, we 
allow for arbitrary bounded variance profile along the main diagonal and
assume bounded moments \eqref{eq:W3} instead of the sub-exponentiality.
% \textcolor{blue}{For
%   the off-diagonal elements we require both \textbf{W1} and
%   \textbf{W2} to hold as in \cite{knowles2013isotropic}.} 

\begin{remark} While \cite{bloemendal2014isotropic} proves the
isotropic law for generalized Wigner matrices (without moment matching
to GUE), the assumptions on variances exclude our setting 
% such matrices are assumed to have a doubly stochastic matrix
% of variances. This does not match our setting, even after a scaling, 
in which the diagonal variances follow arbitrary profiles.
\end{remark}

% Note that although the isotropic law has been established for generalized Wigner matrices in \cite{bloemendal2014isotropic}, such matrices must have a double stochastic matrix of variances, which does not fit (even after a scaling) in our framework, where the diagonal variances follow arbitrary profiles.

\begin{proposition}[Spiked log determinant CLT] \label{prop:zero-diag spiked}
%  Suppose that $W_N$ is drawn from (scaled) GOE or GUE, while
  Let $W_{h,N}$ be a non-critically spiked Wigner matrix as defined above.
  % whose off-diagonal entries have real and imaginary parts whose moments match GOE or GUE to third order.
  Let $E = E_N = 2 + \sigma_N N^{-2/3}$
    with $-\gamma \leq \sigma_N \ll \log^2 N$.
  % with $\ts_N$ a monotone sequence for which $\ts_N \geq -C$ for some
  % positive constant $C$ and $\ts_N = o(\log^2 N)$.
  Then
 % Suppose that
 %  $E = 2 + \ts_N$ with $\ts_N$ a monotone sequence for which
 %  $\ts_N \geq 0$ and $\ts_N = o(\log^2 N)$.
 % With   $\delta_N = \sqrt{\sigma_N}$ we have
  \begin{equation}
  \label{eq:zW-clt-spiked}
  \tau_N^{-1} (\log\abs{\det(W_{h,N} - E)} - \mu_N)  \drightarrow \mathcal{N}(0,1).
\end{equation}

\proof
Note that, for any $E$ such that $W_N-E$ is invertible, we have
\[
\log|\det (W_{h,N}-E)|=\log|\det(W_N-E)|+\log|1+h\mathbf{v^\ast}(W_N-E)^{-1}\mathbf{v}|.
\]
To transfer to $W_{h,N}$ the CLT established for $W_N$ in
Proposition \ref{prop:zero-diag}, it is sufficient to prove that 
% By \cref{prop:zero-diag}, the first log on the right hand side satisfies the CLT after the stated centering and scaling. Hence, to show that the log on the left hand side satisfies the same CLT, it is sufficient to prove that 
\begin{equation}
 \label{spike extra1}   
\tau_N^{-1}\log|1+h\mathbf{v^\ast}(W_N-E)^{-1}\mathbf{v}|\overset{P}{\rightarrow}0.
\end{equation}
Fix $\epsilon > 0$ small, and let $z = E + \im  \eta$ with
$\eta = N^{-1/3-3\epsilon}$. Since $R(z) = (W_N-z)^{-1}$, we have
%Let us show that \cref{thm: isotropic} yields \cref{spike extra1}.
%Let $G(z)=(W_N-z)^{-1}$ with 
%For $z=E+\im \eta$ with $E=2+\tilde{\sigma}_N N^{-2/3}$, 
%We have
\[
\mathbf{v}^\ast R(E)\mathbf{v}=s_{sc}(z)+\left(\mathbf{v}^\ast R(z)\mathbf{v}-s_{sc}(z)\right)+\mathbf{v}^\ast (R(E)-G(z))\mathbf{v}.
\]
Now use spectral decomposition \eqref{eq:spec}, the non-concentraion
bound of Lemma \ref{lem:neg-mom-bd} and finally isotropic
delocalization  (\ref{isotropic delocalization}) to conclude that
with high probability,
% Let $\epsilon>0$ be an arbitrarily small constant. Then, for the last term, we have the following bound
\begin{align*}
  |\mathbf{v}^\ast (R(E)-R(z))\mathbf{v}|
  &\leq \sum_j |\mathbf{v}^\ast \mathbf{u}^{(j)}|^2\left\vert
    \frac{1}{E-\lambda_j}-\frac{1}{E+\im \eta-\lambda_j}\right\vert\\ 
  &\leq \eta \max_j |\mathbf{v}^\ast \mathbf{u}^{(j)}|^2
    \sum_j\frac{1}{(E-\lambda_j)^2} \\
  &  \leq \eta \max_j |\mathbf{v}^\ast \mathbf{u}^{(j)}|^2
    N^{4/3+\epsilon}
  = o_{\Pr}(1).
\end{align*}
% with high probability, where the last inequality follows from
% \cref{lem:neg-mom-bd}.
%Using now \cref{isotropic delocalization}, we conclude that
% \[
% |\mathbf{v}^\ast (R(E)-R(z))\mathbf{v}|=o_{\Pr}(1)
% \]
% as long as, say, $\eta= N^{-1/3-\delta}$ for some fixed small positive $\delta$. Further, for such $\eta$, \cref{isotropic law} yields
Further, for such $\eta$, we have $\Psi(z) \lesssim (N \eta)^{-1/2}$
and so \eqref{isotropic delocalization} yields
\[
\mathbf{v}^\ast R(z)\mathbf{v}-s_{sc}(z)=o_{\Pr}(1).
\]
Finally, $s_{sc}(z)=-1+ O(|E-2|^{1/2} + \eta^{1/2}) = -1 + o(1)$,
(e.g. \cite[Lemma 6.2]{erdos2017dynamical}), and so
% since for such $E$ and $\eta$, $s_{sc}(z)=-1+o(1)$, we conclude that
\[
\mathbf{v}^\ast R(E)\mathbf{v}= -1+o_{\Pr}(1).
\]
Therefore \eqref{spike extra1} holds, since for any fixed $h\neq 1$,
\[
\log |1+h\mathbf{v}^\ast
(W_N-E)^{-1}\mathbf{v}|=\log|1-h|+o_{\Pr}(1)=O_{\Pr}(1). \qedhere
\]
% and \cref{spike extra1} holds.
\end{proposition}

\begin{remark} For the Gaussian ensembles and $\sigma_N$ slowly diverging to infinity, \cref{prop:zero-diag spiked} can be extended to the critical case $h=1$ using the tri-diagonal representations of the spiked GUE/GOE. In such an extension, an extra shift $-\frac{1}{3}\log N$ will ensure the convergence of the normalized log determinant to the standard normal distribution (see \cref{sec:proofs-section-6}).
\end{remark}

%%%%%%%%%%%%%%%%%%%%%%%%%%%%%%%%%%%%%%%%%%%%%%
%% Support information, if any,             %%
%% should be provided in the                %%
%% Acknowledgements section.                %%
%%%%%%%%%%%%%%%%%%%%%%%%%%%%%%%%%%%%%%%%%%%%%%
\begin{acks}[Acknowledgments]
  The authors would like to thank
  Ofer Zeitouni for drawing our attention to
    references \cite{augeri2020} and \cite{bourgade2021}.
\end{acks}
%%%%%%%%%%%%%%%%%%%%%%%%%%%%%%%%%%%%%%%%%%%%%%
%% Funding information, if any,             %%
%% should be provided in the                %%
%% funding section.                         %%
%%%%%%%%%%%%%%%%%%%%%%%%%%%%%%%%%%%%%%%%%%%%%%
\begin{funding}
 The first and fourth authors were supported in part by NSF grant
  DMS 1811614.

%The second author was supported in part by ...
\end{funding}

%%%%%%%%%%%%%%%%%%%%%%%%%%%%%%%%%%%%%%%%%%%%%%
%% Supplementary Material, including data   %%
%% sets and code, should be provided in     %%
%% {supplement} environment with title      %%
%% and short description. It cannot be      %%
%% available exclusively as external link.  %%
%% All Supplementary Material must be       %%
%% available to the reader on Project       %%
%% Euclid with the published article.       %%
%%%%%%%%%%%%%%%%%%%%%%%%%%%%%%%%%%%%%%%%%%%%%%
%\begin{supplement}
%\stitle{???}
%\sdescription{???.}
%\end{supplement}

%%  MR numbers will be added by VTeX.                      %%
%%                                                         %%
%%  Use \cite{...} to cite references in text.             %%
%%                                                         %%
%%%%%%%%%%%%%%%%%%%%%%%%%%%%%%%%%%%%%%%%%%%%%%%%%%%%%%%%%%%%%

%% if your bibliography is in bibtex format, uncomment commands:
\bibliographystyle{imsart-number} % Style BST file (imsart-number.bst or imsart-nameyear.bst)
\bibliography{library_complete}       % Bibliography file (usually '*.bib')

\begin{thebibliography}{51}
% BibTex style file: imsart-number.bst, 2017-11-03
% Default style options (sort=1,type=number).
% Used options (sort=1,type=number).

\bibitem{albrecht2009sturm}
\begin{barticle}[author]
\bauthor{\bsnm{Albrecht},~\bfnm{James~T.}\binits{J.~T.}},
  \bauthor{\bsnm{Chan},~\bfnm{Cy~P.}\binits{C.~P.}} \AND
  \bauthor{\bsnm{Edelman},~\bfnm{Alan}\binits{A.}}
(\byear{2009}).
\btitle{{Sturm sequences and random eigenvalue distributions}}.
\bjournal{Foundations of Computational Mathematics}
\bvolume{9}
\bpages{461--483}.
\end{barticle}
\endbibitem

\bibitem{anderson1997monotoneity}
\begin{barticle}[author]
\bauthor{\bsnm{Anderson},~\bfnm{G.}\binits{G.}} \AND
  \bauthor{\bsnm{Qiu},~\bfnm{S.~L.}\binits{S.~L.}}
(\byear{1997}).
\btitle{{A monotoneity property of the gamma function}}.
\bjournal{Proceedings of the American Mathematical Society}
\bvolume{125}
\bpages{3355--3362}.
\end{barticle}
\endbibitem

\bibitem{anderson2010introduction}
\begin{bbook}[author]
\bauthor{\bsnm{Anderson},~\bfnm{Greg~W.}\binits{G.~W.}},
  \bauthor{\bsnm{Guionnet},~\bfnm{Alice}\binits{A.}} \AND
  \bauthor{\bsnm{Zeitouni},~\bfnm{Ofer}\binits{O.}}
(\byear{2010}).
\btitle{{An Introduction to Random Matrices}}.
\bpublisher{Cambridge university press}.
\end{bbook}
\endbibitem

\bibitem{augeri2020}
\begin{barticle}[author]
\bauthor{\bsnm{Augeri},~\bfnm{Fanny}\binits{F.}},
  \bauthor{\bsnm{Butez},~\bfnm{Raphael}\binits{R.}} \AND
  \bauthor{\bsnm{Zeitouni},~\bfnm{Ofer}\binits{O.}}
(\byear{2020}).
\btitle{{A CLT for the characteristic polynomial of random Jacobi matrices, and
  the G$\beta$E models for $\beta>0$}}.
\bjournal{arXiv preprint arXiv:2011.06870v2}.
\end{barticle}
\endbibitem

\bibitem{basi09}
\begin{bbook}[author]
\bauthor{\bsnm{Bai},~\bfnm{Z.~D.}\binits{Z.~D.}} \AND
  \bauthor{\bsnm{Silverstein},~\bfnm{Jack}\binits{J.}}
(\byear{2009}).
\btitle{{Spectral Analysis of Large Dimensional Random Matrices}},
\bedition{Second} ed.
\bpublisher{Springer}, \baddress{New York}.
\end{bbook}
\endbibitem

\bibitem{BaiBenPech}
\begin{barticle}[author]
\bauthor{\bsnm{Baik},~\bfnm{Jinho}\binits{J.}},
  \bauthor{\bsnm{Ben~Arous},~\bfnm{G{\'{e}}rard}\binits{G.}} \AND
  \bauthor{\bsnm{P{\'{e}}ch{\'{e}}},~\bfnm{Sandrine}\binits{S.}}
(\byear{2005}).
\btitle{{{Phase transition of the largest eigenvalue for nonnull complex sample
  covariance matrices}}}.
\bjournal{Annals of Probability}
\bvolume{33}
\bpages{1643--1697}.
\bdoi{10.1214/009117905000000233}
\end{barticle}
\endbibitem

\bibitem{baik2016fluctuations}
\begin{barticle}[author]
\bauthor{\bsnm{Baik},~\bfnm{Jinho}\binits{J.}} \AND
  \bauthor{\bsnm{Lee},~\bfnm{Ji~Oon}\binits{J.~O.}}
(\byear{2016}).
\btitle{{Fluctuations of the free energy of the spherical
  {S}herrington--{K}irkpatrick model}}.
\bjournal{Journal of Statistical Physics}
\bvolume{165}
\bpages{185--224}.
\end{barticle}
\endbibitem

\bibitem{Baik2017}
\begin{barticle}[author]
\bauthor{\bsnm{Baik},~\bfnm{Jinho}\binits{J.}} \AND
  \bauthor{\bsnm{Lee},~\bfnm{Ji~Oon}\binits{J.~O.}}
(\byear{2017}).
\btitle{{{Fluctuations of the Free Energy of the Spherical
  Sherrington–Kirkpatrick Model with Ferromagnetic Interaction}}}.
\bjournal{Annales Henri Poincare}
\bvolume{18}
\bpages{1867--1917}.
\bdoi{10.1007/s00023-017-0562-5}
\end{barticle}
\endbibitem

\bibitem{BaikLeeWu}
\begin{barticle}[author]
\bauthor{\bsnm{Baik},~\bfnm{Jinho}\binits{J.}},
  \bauthor{\bsnm{Lee},~\bfnm{Oon~Ji}\binits{O.~J.}} \AND
  \bauthor{\bsnm{Wu},~\bfnm{Hao}\binits{H.}}
(\byear{2018}).
\btitle{{{Ferromagnetic to Paramagnetic Transition in Spherical Spin Glass}}}.
\bjournal{Journal of Statistical Physics}
\bvolume{173}
\bpages{1484--1522}.
\bdoi{10.1007/s10955-018-2150-6}
\end{barticle}
\endbibitem

\bibitem{benaych2018}
\begin{barticle}[author]
\bauthor{\bsnm{Benaych-George},~\bfnm{Florent}\binits{F.}} \AND
  \bauthor{\bsnm{Knowles},~\bfnm{Antti}\binits{A.}}
(\byear{2018}).
\btitle{{Lectures on the local semicircle law for Wigner matrices}}.
\bjournal{arXiv:1601.04055v4}.
\end{barticle}
\endbibitem

\bibitem{bloemendal2014isotropic}
\begin{barticle}[author]
\bauthor{\bsnm{Bloemendal},~\bfnm{Alex}\binits{A.}},
  \bauthor{\bsnm{L\'{a}szl\'{o}},~\bfnm{Erd\H{o}s}\binits{E.}},
  \bauthor{\bsnm{Knowles},~\bfnm{Antti}\binits{A.}},
  \bauthor{\bsnm{Yau},~\bfnm{Horng-Tzer}\binits{H.-T.}} \AND
  \bauthor{\bsnm{Yin},~\bfnm{Jun}\binits{J.}}
(\byear{2014}).
\btitle{{Isotropic local laws for sample covariance and generalized Wigner
  matrices}}.
\bjournal{Electronic Journal of Probability}
\bvolume{19}
\bpages{1--53}.
\end{barticle}
\endbibitem

\bibitem{boucheron2013concentration}
\begin{bbook}[author]
\bauthor{\bsnm{Boucheron},~\bfnm{St{\'e}phane}\binits{S.}},
  \bauthor{\bsnm{Lugosi},~\bfnm{G{\'a}bor}\binits{G.}} \AND
  \bauthor{\bsnm{Massart},~\bfnm{Pascal}\binits{P.}}
(\byear{2013}).
\btitle{{Concentration inequalities: A nonasymptotic theory of independence}}.
\bpublisher{Oxford university press}.
\end{bbook}
\endbibitem

\bibitem{bourgade2019gaussian}
\begin{barticle}[author]
\bauthor{\bsnm{Bourgade},~\bfnm{P.}\binits{P.}} \AND
  \bauthor{\bsnm{Mody},~\bfnm{K.}\binits{K.}}
(\byear{2019}).
\btitle{{Gaussian fluctuations of the determinant of {W}igner matrices}}.
\bjournal{Electronic Journal of Probability}
\bvolume{24}
\bpages{1--28}.
\end{barticle}
\endbibitem

\bibitem{bourgade2021}
\begin{barticle}[author]
\bauthor{\bsnm{Bourgade},~\bfnm{Paul}\binits{P.}},
  \bauthor{\bsnm{Mody},~\bfnm{Krishnan}\binits{K.}} \AND
  \bauthor{\bsnm{Pain},~\bfnm{Michel}\binits{M.}}
(\byear{2021}).
\btitle{{Optimal local law and central limit theorem for $\beta$-ensembles}}.
\bjournal{arXiv preprint arXiv:2103.06841v2}.
\end{barticle}
\endbibitem

\bibitem{duy2017distributionsr}
\begin{bincollection}[author]
\bauthor{\bsnm{Duy},~\bfnm{T.~K.}\binits{T.~K.}}
(\byear{2017}).
\btitle{{Distributions of the determinants of {G}aussian beta ensembles}}.
In \bbooktitle{2023 Spectral and Scattering Theory and Related Topics}
\bpages{77--85}.
\bpublisher{RIMS Kokyuroku}.
\end{bincollection}
\endbibitem

\bibitem{erde60}
\begin{barticle}[author]
\bauthor{\bsnm{Erd\'{e}lyi},~\bfnm{A.}\binits{A.}}
(\byear{1960}).
\btitle{{Asymptotic solutions of differential equations with transition points
  or singularities}}.
\bjournal{J. Mathematical Phys.}
\bvolume{1}
\bpages{16--26}.
\bdoi{10.1063/1.1703631}
\bmrnumber{111915}
\end{barticle}
\endbibitem

\bibitem{erdos2017dynamical}
\begin{bbook}[author]
\bauthor{\bsnm{Erd\H{o}s},~\bfnm{L\'{a}szl\'{o}}\binits{L.}} \AND
  \bauthor{\bsnm{Yau},~\bfnm{Horng-Tzer}\binits{H.-T.}}
(\byear{2017}).
\btitle{{A dynamical approach to random matrix theory}}.
\bseries{Courant Lecture Notes in Mathematics}
\bvolume{28}.
\bpublisher{Courant Institute of Mathematical Sciences, New York; American
  Mathematical Society, Providence, RI}.
\end{bbook}
\endbibitem

\bibitem{ErdosYY2012b}
\begin{barticle}[author]
\bauthor{\bsnm{Erd\H{o}s},~\bfnm{L\'{a}azl\'{o}}\binits{L.}},
  \bauthor{\bsnm{Yau},~\bfnm{Horng-Tzer}\binits{H.-T.}} \AND
  \bauthor{\bsnm{Yin},~\bfnm{Jun}\binits{J.}}
(\byear{2012}).
\btitle{{{Bulk universality for generalized Wigner matrices}}}.
\bjournal{Probability Theory and Related Fields}
\bvolume{154}
\bpages{341--407}.
\end{barticle}
\endbibitem

\bibitem{etemadi1985some}
\begin{barticle}[author]
\bauthor{\bsnm{Etemadi},~\bfnm{Nasrollah}\binits{N.}}
(\byear{1985}).
\btitle{{On some classical results in probability theory}}.
\bjournal{Sankhy{\=a}: The Indian Journal of Statistics, Series A}
\bpages{215--221}.
\end{barticle}
\endbibitem

\bibitem{forrester2010loggases}
\begin{bbook}[author]
\bauthor{\bsnm{Forrester},~\bfnm{P.~J.}\binits{P.~J.}}
(\byear{2010}).
\btitle{{Log-Gases and Random Matrices}}.
\bpublisher{{P}rinceton {U}niversity {P}ress}.
\end{bbook}
\endbibitem

\bibitem{forrester2001inter}
\begin{bincollection}[author]
\bauthor{\bsnm{Forrester},~\bfnm{P.~J.}\binits{P.~J.}} \AND
  \bauthor{\bsnm{Rains},~\bfnm{E.~M.}\binits{E.~M.}}
(\byear{2001}).
\btitle{{Inter-relationships between orthogonal, unitary and symplectic matrix
  ensembles}}.
In \bbooktitle{Random Matrix Models and Their Applications. Mathematical
  Sciences Research Institute publications},
\bvolume{40}
\bpages{171--207}.
\bpublisher{Cambridge University Press}.
\end{bincollection}
\endbibitem

\bibitem{fyodorov2016fractional}
\begin{barticle}[author]
\bauthor{\bsnm{Fyodorov},~\bfnm{Y.~V.}\binits{Y.~V.}},
  \bauthor{\bsnm{Khoruzhenko},~\bfnm{B.~A.}\binits{B.~A.}},
  \bauthor{\bsnm{Simm},~\bfnm{N.~J.}\binits{N.~J.}} \betal{et~al.}
(\byear{2016}).
\btitle{{Fractional Brownian motion with Hurst index $ H= 0$ and the {G}aussian
  {U}nitary {E}nsemble}}.
\bjournal{The Annals of Probability}
\bvolume{44}
\bpages{2980--3031}.
\end{barticle}
\endbibitem

\bibitem{Gotze2005}
\begin{barticle}[author]
\bauthor{\bsnm{G{\"{o}}tze},~\bfnm{F.}\binits{F.}} \AND
  \bauthor{\bsnm{Tikhomirov},~\bfnm{A.}\binits{A.}}
(\byear{2005}).
\btitle{{{The rate of convergence for spectra of GUE and LUE matrix
  ensembles}}}.
\bjournal{Central European Journal of Mathematics}
\bvolume{3}
\bpages{666--704}.
\bdoi{10.2478/bf02475626}
\end{barticle}
\endbibitem

\bibitem{gradshteyn2007}
\begin{bbook}[author]
\bauthor{\bsnm{Gradshteyn},~\bfnm{I.~S.}\binits{I.~S.}} \AND
  \bauthor{\bsnm{Ryzhik},~\bfnm{I.~M.}\binits{I.~M.}}
(\byear{2007}).
\btitle{{Table of Integrals, Series, and Products}},
\bedition{seventh} ed.
\bpublisher{Elsevier Inc.}
\end{bbook}
\endbibitem

\bibitem{james1954}
\begin{barticle}[author]
\bauthor{\bsnm{James},~\bfnm{A.~T.}\binits{A.~T.}}
(\byear{1954}).
\btitle{Normal multivariate analysis and the orthogonal group}.
\bjournal{Ann. Math. Statistics}
\bvolume{25}
\bpages{40--75}.
\bmrnumber{0060779 (15,726b)}
\end{barticle}
\endbibitem

\bibitem{johansson1998fluctuations}
\begin{barticle}[author]
\bauthor{\bsnm{Johansson},~\bfnm{Kurt}\binits{K.}}
(\byear{1998}).
\btitle{{On fluctuations of eigenvalues of random Hermitian matrices}}.
\bjournal{Duke Mathematical Journal}
\bvolume{91}
\bpages{151--204}.
\end{barticle}
\endbibitem

\bibitem{johnstone2021}
\begin{barticle}[author]
\bauthor{\bsnm{Johnstone},~\bfnm{Iain~M.}\binits{I.~M.}},
  \bauthor{\bsnm{Klochkov},~\bfnm{Yegor}\binits{Y.}},
  \bauthor{\bsnm{Onatski},~\bfnm{Alexei}\binits{A.}} \AND
  \bauthor{\bsnm{Pavlyshyn},~\bfnm{Damian}\binits{D.}}
(\byear{2021}).
\btitle{{Spin glass to paramagnetic transition in spherical
  {S}herington-{K}irkpatrick model with ferromagnetic interaction}}.
\bjournal{Manuscript in preparation}.
\end{barticle}
\endbibitem

\bibitem{johnstone2012fast}
\begin{barticle}[author]
\bauthor{\bsnm{Johnstone},~\bfnm{Iain~M.}\binits{I.~M.}} \AND
  \bauthor{\bsnm{Ma},~\bfnm{Zongming}\binits{Z.}}
(\byear{2012}).
\btitle{{{Fast approach to the Tracy-Widom law at the edge of GOE and GUE}}}.
\bjournal{The Annals of Applied Probability}
\bvolume{22}
\bpages{1962--1988}.
\end{barticle}
\endbibitem

\bibitem{johnstone2020testing}
\begin{barticle}[author]
\bauthor{\bsnm{Johnstone},~\bfnm{Iain~M.}\binits{I.~M.}} \AND
  \bauthor{\bsnm{Onatski},~\bfnm{Alexei}\binits{A.}}
(\byear{2020}).
\btitle{{Testing in high-dimensional spiked models}}.
\bjournal{Annals of Statistics}
\bvolume{48}
\bpages{1231--1254}.
\end{barticle}
\endbibitem

\bibitem{knowles2013isotropic}
\begin{barticle}[author]
\bauthor{\bsnm{Knowles},~\bfnm{Antti}\binits{A.}} \AND
  \bauthor{\bsnm{Yin},~\bfnm{Jun}\binits{J.}}
(\byear{2013}).
\btitle{{The Isotropic Semicircle Law and Deformation of Wigner Matrices}}.
\bjournal{Communications on Pure and Applied Mathematics}
\bvolume{LXVI}
\bpages{1663--1749}.
\end{barticle}
\endbibitem

\bibitem{Kosterlitz1976}
\begin{barticle}[author]
\bauthor{\bsnm{Kosterlitz},~\bfnm{J.~M.}\binits{J.~M.}},
  \bauthor{\bsnm{Thouless},~\bfnm{D.~J.}\binits{D.~J.}} \AND
  \bauthor{\bsnm{Jones},~\bfnm{Raymund~C.}\binits{R.~C.}}
(\byear{1976}).
\btitle{{{Spherical model of a spin-glass}}}.
\bjournal{Physical Review Letters}
\bvolume{36}
\bpages{1217--1220}.
\bdoi{10.1103/PhysRevLett.36.1217}
\end{barticle}
\endbibitem

\bibitem{krasovsky2007correlations}
\begin{barticle}[author]
\bauthor{\bsnm{Krasovsky},~\bfnm{I.~V.}\binits{I.~V.}}
(\byear{2007}).
\btitle{{Correlations of the characteristic polynomials in the {G}aussian
  {U}nitary {E}nsemble or a singular {H}ankel determinant}}.
\bjournal{Duke Mathematical Journal}
\bvolume{139}
\bpages{581--619}.
\end{barticle}
\endbibitem

\bibitem{lambert2020a}
\begin{barticle}[author]
\bauthor{\bsnm{Lambert},~\bfnm{Gaultier}\binits{G.}} \AND
  \bauthor{\bsnm{Paquette},~\bfnm{Elliot}\binits{E.}}
(\byear{2020}).
\btitle{{Strong approximation of {G}aussian $\beta$-ensemble characteristic
  polynomials: the hyperbolic regime}}.
\bjournal{arXiv preprint arXiv:2001.09042}.
\end{barticle}
\endbibitem

\bibitem{lambert2020b}
\begin{barticle}[author]
\bauthor{\bsnm{Lambert},~\bfnm{Gaultier}\binits{G.}} \AND
  \bauthor{\bsnm{Paquette},~\bfnm{Elliot}\binits{E.}}
(\byear{2020}).
\btitle{{Strong approximation of {G}aussian $\beta$-ensemble characteristic
  polynomials: the edge regime and the stochastic {A}iry function}}.
\bjournal{arXiv preprint arXiv:2009.05003}.
\end{barticle}
\endbibitem

\bibitem{LandonSosoeAAP}
\begin{barticle}[author]
\bauthor{\bsnm{Landon},~\bfnm{Benjamin}\binits{B.}} \AND
  \bauthor{\bsnm{Sosoe},~\bfnm{Philippe}\binits{P.}}
(\byear{2020}).
\btitle{Applications of mesoscopic {CLT}s in random matrix theory}.
\bjournal{Ann. Appl. Probab.}
\bvolume{30}
\bpages{2769--2795}.
\bdoi{10.1214/20-AAP1572}
\bmrnumber{4187127}
\end{barticle}
\endbibitem

\bibitem{LiSchnelliXu}
\begin{barticle}[author]
\bauthor{\bsnm{Li},~\bfnm{Yiting}\binits{Y.}},
  \bauthor{\bsnm{Schnelli},~\bfnm{Kevin}\binits{K.}} \AND
  \bauthor{\bsnm{Xu},~\bfnm{Yuanyuan}\binits{Y.}}
(\byear{2017}).
\btitle{{{Central limit theorem for mesoscopic eigenvalue statistics of
  deformed Wigner matrices and sample covariance matrices}}}.
\bjournal{arXiv:1909.12821v2}
\bpages{1--36}.
\end{barticle}
\endbibitem

\bibitem{maid07}
\begin{barticle}[author]
\bauthor{\bsnm{Ma\"{i}da},~\bfnm{Myl\`ene}\binits{M.}}
(\byear{2007}).
\btitle{{Large deviations for the largest eigenvalue of rank one deformations
  of {G}aussian ensembles}}.
\bjournal{Electron. J. Probab.}
\bvolume{12}
\bpages{1131--1150}.
\bdoi{10.1214/EJP.v12-438}
\bmrnumber{2336602}
\end{barticle}
\endbibitem

\bibitem{meht91}
\begin{bbook}[author]
\bauthor{\bsnm{Mehta},~\bfnm{Madan~Lal}\binits{M.~L.}}
(\byear{1991}).
\btitle{{Random matrices}},
\bedition{Second} ed.
\bpublisher{Academic Press, Inc., Boston, MA}.
\bmrnumber{1083764}
\end{bbook}
\endbibitem

\bibitem{meht04}
\begin{bbook}[author]
\bauthor{\bsnm{Mehta},~\bfnm{Madan~Lal}\binits{M.~L.}}
(\byear{2004}).
\btitle{{Random matrices}},
\bedition{third} ed.
\bseries{Pure and Applied Mathematics (Amsterdam)}
\bvolume{142}.
\bpublisher{Elsevier/Academic Press, Amsterdam}.
\bmrnumber{2129906}
\end{bbook}
\endbibitem

\bibitem{muir82}
\begin{bbook}[author]
\bauthor{\bsnm{Muirhead},~\bfnm{Robb~J.}\binits{R.~J.}}
(\byear{1982}).
\btitle{{Aspects of multivariate statistical theory}}.
\bseries{Wiley Series in Probability and Mathematical Statistics}.
\bpublisher{John Wiley \& Sons, Inc., New York}.
\bmrnumber{652932}
\end{bbook}
\endbibitem

\bibitem{olve74}
\begin{bbook}[author]
\bauthor{\bsnm{Olver},~\bfnm{F.~W.~J.}\binits{F.~W.~J.}}
(\byear{1974}).
\btitle{Asymptotics and Special Functions}.
\bpublisher{Academic Press}.
\end{bbook}
\endbibitem

\bibitem{Ona2013}
\begin{barticle}[author]
\bauthor{\bsnm{Onatski},~\bfnm{Alexei}\binits{A.}},
  \bauthor{\bsnm{Moreira},~\bfnm{Marcelo~J}\binits{M.~J.}} \AND
  \bauthor{\bsnm{Hallin},~\bfnm{Marc}\binits{M.}}
(\byear{2013}).
\btitle{{{Asymptotic power of sphericity tests for high-dimensional data}}}.
\bjournal{Annals of Statistics}
\bvolume{41}
\bpages{1204--1231}.
\bdoi{10.1214/13-AOS1100}
\end{barticle}
\endbibitem

\bibitem{Patterson2006}
\begin{barticle}[author]
\bauthor{\bsnm{Patterson},~\bfnm{Nick}\binits{N.}},
  \bauthor{\bsnm{Price},~\bfnm{Alkes~L}\binits{A.~L.}} \AND
  \bauthor{\bsnm{Reich},~\bfnm{David}\binits{D.}}
(\byear{2006}).
\btitle{Population Structure and Eigenanalysis}.
\bjournal{PLoS Genet}
\bvolume{2}
\bpages{e190}.
\bdoi{10.1371/journal.pgen.0020190}
\end{barticle}
\endbibitem

\bibitem{Peche2006a}
\begin{barticle}[author]
\bauthor{\bsnm{P{\'{e}}ch{\'{e}}},~\bfnm{S.}\binits{S.}}
(\byear{2006}).
\btitle{{{The largest eigenvalue of small rank perturbations of Hermitian
  random matrices}}}.
\bjournal{Probability Theory and Related Fields}
\bvolume{134}
\bpages{127--173}.
\bdoi{10.1007/s00440-005-0466-z}
\end{barticle}
\endbibitem

\bibitem{rio2009moment}
\begin{barticle}[author]
\bauthor{\bsnm{Rio},~\bfnm{Emmanuel}\binits{E.}}
(\byear{2009}).
\btitle{{Moment inequalities for sums of dependent random variables under
  projective conditions}}.
\bjournal{Journal of Theoretical Probability}
\bvolume{22}
\bpages{146--163}.
\end{barticle}
\endbibitem

\bibitem{tao2012central}
\begin{barticle}[author]
\bauthor{\bsnm{Tao},~\bfnm{Terence}\binits{T.}} \AND
  \bauthor{\bsnm{Vu},~\bfnm{Van}\binits{V.}}
(\byear{2012}).
\btitle{{{A central limit theorem for the determinant of a {W}igner matrix}}}.
\bjournal{Advances in Mathematics}
\bvolume{231}
\bpages{74--101}.
\end{barticle}
\endbibitem

\bibitem{trwi96}
\begin{barticle}[author]
\bauthor{\bsnm{Tracy},~\bfnm{Craig~A.}\binits{C.~A.}} \AND
  \bauthor{\bsnm{Widom},~\bfnm{Harold}\binits{H.}}
(\byear{1996}).
\btitle{On orthogonal and symplectic matrix ensembles}.
\bjournal{Communications in Mathematical Physics}
\bvolume{177}
\bpages{727-754}.
\end{barticle}
\endbibitem

\bibitem{tracy1998correlation}
\begin{barticle}[author]
\bauthor{\bsnm{Tracy},~\bfnm{Craig~A.}\binits{C.~A.}} \AND
  \bauthor{\bsnm{Widom},~\bfnm{Harold}\binits{H.}}
(\byear{1998}).
\btitle{{Correlation functions, cluster functions, and spacing distributions
  for random matrices}}.
\bjournal{Journal of Statistical Physics}
\bvolume{92}
\bpages{809--835}.
\end{barticle}
\endbibitem

\bibitem{trotter1984eigenvalue}
\begin{barticle}[author]
\bauthor{\bsnm{Trotter},~\bfnm{Hale~F.}\binits{H.~F.}}
(\byear{1984}).
\btitle{{{Eigenvalue distributions of large {H}ermitian matrices; {W}igner's
  semi-circle law and a theorem of {K}ac, {M}urdock, and {S}zeg{\"o}}}}.
\bjournal{Advances in mathematics}
\bvolume{54}
\bpages{67--82}.
\end{barticle}
\endbibitem

\bibitem{wigner1965distribution}
\begin{barticle}[author]
\bauthor{\bsnm{Wigner},~\bfnm{Eugene~P.}\binits{E.~P.}}
(\byear{1965}).
\btitle{{Distribution laws for the roots of a random {H}ermitian matrix}}.
\bjournal{Statistical Theories of Spectra: Fluctuations}
\bpages{446--461}.
\end{barticle}
\endbibitem

\bibitem{Yao2005}
\begin{barticle}[author]
\bauthor{\bsnm{Yao},~\bfnm{Jianfeng}\binits{J.}} \AND
  \bauthor{\bsnm{Bai},~\bfnm{Zhidong}\binits{Z.}}
(\byear{2005}).
\btitle{{{On the convergence of the spectral empirical process of Wigner
  matrices}}}.
\bjournal{Bernoulli}
\bvolume{11}
\bpages{1059--1092}.
\end{barticle}
\endbibitem

\end{thebibliography}
%\bibliography{mybib}       % Bibliography file (usually '*.bib')

%\printbibliography

%% or include bibliography directly:
% \begin{thebibliography}{}
% \bibitem{b1}
% \end{thebibliography}

\newpage
\begin{centering}
\textbf{Supplementary Material for ``An edge CLT for the log determinant of Wigner matrices'' by I.M. Johnstone, Y. Klochkov, A. Onatski, and D. Pavlyshyn}
\vspace{1cm}

    \textbf{Abstract: }\small{The Supplementary Material contains relatively more technical proofs of the main paper. It consists of Appendix A, which corresponds to the proofs for Gaussian ensembles, and Appendix B, which corresponds to  proofs for Wigner extension. To help the reader to navigate the Supplement, we start from a Table of Contents that encompasses the main text. The Table's references to the Supplement's sections contain short descriptions of the content of the sections.}
\end{centering}
\tableofcontents

%%%%%%%%%%%%%%%%%%%%%%%%%%%%%%%%%%%%%%%%%%%%%%
%% Single Appendix:                         %%
%%%%%%%%%%%%%%%%%%%%%%%%%%%%%%%%%%%%%%%%%%%%%%
%\begin{appendix}
%\section*{???}%% if no title is needed, leave empty \section*{}.
%\end{appendix}
%%%%%%%%%%%%%%%%%%%%%%%%%%%%%%%%%%%%%%%%%%%%%%
%% Multiple Appendixes:                     %%
%%%%%%%%%%%%%%%%%%%%%%%%%%%%%%%%%%%%%%%%%%%%%%
\begin{appendix}

%\texttt{What we (for now) call the Appendices will become the Supplementary
%Material (but let's keep in a single document while we write).
%Our initial target length for the main paper + references (but excluding
%Appendices) might be the ``average length'' for AAP papers, said to
%about 32 pages (2016-18 sample). Thus we may need to iteratively move
%material to the Appendices to stay within our target.}

\section{Proofs for the Gaussian ensembles}
\label{sec:proofs-gauss-ensembl-1}

\subsection{Proof of Lemma \ref{edge spacing lemma} (about relative location of $\lambda_i$ and $E$)}
\label{sec:proof-lemma-4}

Let $\lambda_{1} \geq \cdots \geq \lambda_{N}$ be
eigenvalues of $W = W_{N}$ and define the eigenvalue counting
function $N_W(I) = \# \{ j: \lambda_{j} \in I \}$.
The first bound in \eqref{eq:first-pair} follows from
  Tracy-Widom convergence of $\lambda_1$.

For the second bound, set $I_\gamma^+ = [2-\gamma
N^{-2/3}, \infty)$ 
and note that
\begin{equation*}
  \Pr(\lambda_k > E)
  \leq \Pr(\lambda_k > 2 - \gamma N^{-2/3})
  = \Pr( N_W(I_\gamma^+) > k-1) ).
\end{equation*}
Using the tail bound \eqref{eq:onept-bds}
% \eqref{johnstone_ma1}
for the one-point function,
\begin{align*}
  \E N_W(I_\gamma^+)
  & = \int_{2-\gamma N^{-2/3}}^\infty N \rho_{N}(x) \diff x
    = \int_{-\gamma}^\infty N^{1/3} \rho_{N}(2+sN^{-2/3}) \diff
    s \\
  & \leq C_1(\gamma) \int_{-\gamma}^\infty e^{-2s} \diff s
    \leq C_2(\gamma).
\end{align*}
The last two displays and Markov's inequality yield
$\Pr(\lambda_k > E) \leq C_2(\gamma)/(k-1) < \varepsilon$ so long as
$k = k(\varepsilon,\gamma)$ is sufficiently large.

% The second bound is immediate from Lemma \ref{lem:anticoncGOE}(ii), on
% choosing $c_1 < \epsilon/C$.
% For the second bound let $I_c = [E-cN^{-2/3},E+cN^{-2/3}]$ and again
% use \eqref{johnstone_ma1}. With $W = W_{N,\alpha}$,
% \begin{align*}
%   \Pr \big( \min_{1 \leq j \leq N} |\lambda_j-E| \leq c_1 N^{-2/3} \big)
%   & = \Pr \big( N_W(I_{c_1}) \geq 1 \big)  \\
%   &  \leq \E N_W(I_{c_1}) 
%    = \int_{I_{c_1}} N \rho_{N,\alpha}(x) \diff x \\
%   &  \leq N \rho_{N,\alpha}(2-N^{-2/3}(\gamma+c_1)) |I_{c_1}|
%     \leq 2 c_1 C(\gamma+c_1) < \varepsilon
% \end{align*}
% if $c_1 = c_1(\varepsilon,\gamma)$ is chosen sufficiently small.

For the first bound in \eqref{eq:second-pair},
set $x_{jN} = N^{2/3}(\lambda_j-2)$.
Using~\eqref{eq:first-pair},
% the first bound,
choose $k = k(\epsilon/2, \gamma+1)$ so that 
the event $\mathcal{E}_{kN} = \{ x_{kN} < -\gamma-1 \}$ has
probability at least $1 - \epsilon/2$.
On $\mathcal{E}_{kN}$, for $j \geq k$ we have
$x_{jN} - \sigma_N \leq x_{kN} + \gamma < -1$.

Consider now $j < k$. Since the $j$th Tracy-Widom law $F_j$ (of type 2 or 1 for the GUE or GOE case, respectively) has a
continuous distribution function, weak convergence of $x_{jN}$ also
implies that $\Pr\{ x_{jN} \leq x \} \to F_j(x)$ uniformly in
$x$. Since $F_j$ is uniformly continuous, we can choose $c_1$ small so
that for large $N$ and each $j < k$,
\begin{equation*}
  \Pr \{ x_{jN} \in [\sigma_N-c_1, \sigma_N + c_1]\} \leq
  \epsilon/2k.
\end{equation*}
Let $\mathcal{E}_{jN}^c$ be the corresponding event.
On the event $\cap_{j \leq k} \mathcal{E}_{jN}$, which has probability
at least $1 - \epsilon$, we have
$\min_{1 \leq j \leq N} N^{2/3}|\lambda_j-E| = \min_j |x_{jN} -
\sigma_N| \geq c_1$,  and so
the first bound in \eqref{eq:second-pair}
%the second bound
is proved.

For the second bound in \eqref{eq:second-pair},
%third bound,
% set $x_j = N^{2/3}(\lambda_j-2)$.
since $N^{2/3}|\lambda_j-E| \leq |x_{jN}| + |\sigma_N|$, we have
\begin{equation*}
  \Pr \big( \max_{j\leq k} N^{2/3}|\lambda_j-E| > C_1 + |\sigma_N|
  \big)
  \leq \Pr \big( \max_{j\leq k} |x_{jN}| > C_1 \big)
  \leq \Pr(x_{1N} > C_1) + \Pr(x_{kN} < -C_1).
\end{equation*}
Again using Tracy-Widom convergence of
$x_{1N}$ and $x_{kN}$, the right side can be made less than
$\epsilon$ for large $N$ by choosing $C_1$ large.

%\textcolor{red}{[Lemma 4 does not make statements about spiked $W_N$. Why do we need the following here?]}

%Consider now the subcritical spiked case, with eigenvalues
%$W_{h,N}=\tilde{\lambda}_j$ of $W_{N} + h \mathbf{v}\mathbf{v}^*$.
%By eigenvalue interlacing, $\lambda_j < \tilde{\lambda}_j < \lambda_{j-1}$, and so
%$\Pr \{\tilde{\lambda}_{k+1} > E\} \leq \Pr \{\lambda_k > E\}$
%Choose therefore $k = k(\epsilon/2, \gamma+1) + 1$ to recover the
%first bound in the spiked case.

%The remainder of the argument goes through with $x_{jN}$ replaced by
%$\tilde{x}_{jN} = N^{2/3} (\tilde{\lambda}_j-2)$
%since, as noted, the latter still converges weakly to $F_j$ for $h < 1$. \textcolor{red}{[Don't we need super-critical case too for the tri-diagonal based proof of spiked CLT?]}

\subsection{Proofs of technical lemmas from Section \ref{sec:proofs-gauss-ensembl}}
%\ref{section away from edge}}
\label{appendix away from edge}

%\subsubsection{A moment bound for sub-gamma random variables} 
%\label{sec:moment-bound-sub}

\subsubsection{Proof of Lemma \ref{lemma moments of ksi} (bounds on even moments of $\xi_i$)}\label{proof of lemma moments}

The existence of $C_{q}$ follows from the fact that $\xi _{i}$ are sub-gamma
random variables (see Definition~\eqref{definitions SG} and equation (\ref%
{alpha plus beta SG})) and from equation (2.7) of \cite{boucheron2013concentration}.
The lower bound on $\E\xi _{i}^{2}$ follows from
$m_ir_i = (i-1)/N \theta_N^2$: 
\begin{equation}
\E\xi _{i}^{2}=\frac{\alpha}{N\theta _{N}^{2}r_{i}^{2}}+\frac{\alpha
  m_{i}r_{i}}{
  N\theta _{N}^{2}r_{i}^{2}r_{i-1}^{2}}
  = \frac{\alpha}{N\theta _{N}^{2}r_{i}^{2}} \Big[ 1 + \frac{i-1}{N
     \theta_N^2 r_{i-1}^2} \Big].
\label{variance of ksi}
\end{equation}

\subsubsection{Proof of Lemma \ref{lemma si} (bounds on $g_i$, used to verify the Lyapunov condition)}\label{proof of lemma si}
Since $\gamma _{i}$ is increasing with $i,$ we have%
\[
g_{i}>1+\gamma _{i}+...+\gamma _{i}^{N-i+1}=\frac{1-\gamma _{i}^{N-i+2}}{%
1-\gamma _{i}}.
\]%
On the other hand, for all sufficiently large $N$%
\[
\gamma _{i}<m_{i}\leq 1-\sqrt{1-\frac{N}{N\theta _{N}^{2}}}=1-\sqrt{1-\theta
_{N}^{-2}}<1-N^{-1/3}w_{N}^{1/2}.
\]%
Hence, for $i\leq N-N^{1/3},$ any $k>0,$ and all sufficiently large $N,$%
\[
\gamma _{i}^{N-i+2}<\left( 1-N^{-1/3}w_{N}^{1/2}\right)
^{N^{1/3}}<e^{-w_{N}^{1/2}}<e^{-k\log \log N}=\log ^{-k}N.
\]%
The lower bound (\ref{lower bound on si}) follows from this and the fact
that $1-\gamma _{i}=2\left( r_{i}-1\right) /r_{i}.$

The elementary bound \eqref{eq:g-trivial} follows from the definition
of $g_i$. It shows that the lower bound \eqref{lower bound on si}
fails for $i > N-N^\alpha$ for any $\alpha < 1/3$.

For the upper bound, we seek a value $\kappa$ for which the
inequalities $(1-\gamma_i)g_i \leq 1+\kappa$ may be established by
induction for $i =N, N-1, \ldots, 1$. 
The initial step holds for any $\kappa \geq 0$, since $g_N =
1+\gamma_N$ implies that $(1-\gamma_N)g_N = 1-\gamma_N^2$.
Assuming $(1-\gamma_i)g_i \leq 1+\kappa$ and using the recursion
$g_{i-1} = \gamma_{i-1}g_i+1$, we have
\begin{equation*}
  (1-\gamma_{i-1})g_{i-1} \leq (1-\gamma_{i-1})
    \left[\frac{\gamma_{i-1}}{1-\gamma_i}(1+\kappa) + 1\right],
\end{equation*}
and so the induction step works at least so long as
\begin{equation*}
  \frac{1-\gamma_{i-1}}{1-\gamma_i} \gamma_{i-1}(1+\kappa) -
  \gamma_{i-1} \leq \kappa
\end{equation*}
for $2 \leq i \leq N$. On rearrangement this condition becomes
\begin{equation}  \label{eq:equiv-cond}
  \gamma_{i-1}(\gamma_i-\gamma_{i-1})
    \leq \kappa(1-\gamma_i - \gamma_{i-1} + \gamma_{i-1}^2).
  \end{equation}
Both sides of the inequality are monotone in such a
way that we need only work with $i=N$.
Indeed, write
$\gamma_i = \gamma(x_i)$ and note that
 $\gamma(x)=(1-R)/(1+R)$ and $\gamma'(x) = R^{-1}(1+R)^{-2}$ are
increasing, where $R(x)=\sqrt{1-x}$. Also set $R_N = R(x_N) = \sqrt{1-x_N}$.
We then have both
\begin{align*}
  \gamma_{i-1}(\gamma_i-\gamma_{i-1})
  & \leq \gamma_N \gamma_N' \Delta_N \leq \Delta_N/[R_N(1+R_N)^2]
  \\
  \intertext{and}
  1-\gamma_i - \gamma_{i-1} + \gamma_{i-1}^2
  & = (1-\gamma_{i-1})^2-(\gamma_i-\gamma_{i-1}) \\
  &  \geq (1-\gamma_N)^2-\gamma_N'\Delta_N 
    = (4R_N^3-\Delta_N)/[R_N(1+R_N)^2].
\end{align*}
From these displays and \eqref{eq:equiv-cond}, we see that any
$\kappa$ larger than $\Delta_N/(4R_N^3-\Delta_N)$ suffices for the
induction. Noting that for large enough $N$ we have $R_N > \sqrt{w_N}
  N^{-1/3}$ and
  \begin{equation*}
    4R_N^3/\Delta_N -1 > 4w_N^{3/2}\theta_N^2-1 > 3w_N^{3/2},
  \end{equation*}
we conclude that we may certainly take $\kappa = w_N^{-3/2}$ in the
induction and in our upper bound.

\subsubsection{Proof of Lemma \ref{lemma variance} (improved bounds on $\mathbf{E}\xi_i^2$)}\label{proof of lemma variance}

Using (\ref{variance of ksi}) and the 
identities $m_{i}r_{i}-m_{i-1}r_{i-1}=1/ N\theta_{N}^{2}$ and
$m_{i-1}+r_{i-1}=2$
we have
\begin{equation*}
  \alpha^{-1} N \theta_N^2 r_i^2 \E \xi_i^2
  = 1 + \frac{m_ir_i}{r_{i-1}^2}
  = 1 + \frac{m_{i-1}}{r_{i-1}} + \frac{1}{N \theta_N^2 r_{i-1}^2}
  = \frac{2}{r_i} \Big[ 1 + \frac{r_i-r_{i-1}}{r_{i-1}} \Big]
     + \frac{1}{N \theta_N^2 r_{i-1}^2}.
\end{equation*}
Appealing to the monotonicity of $r_i \in [1,2]$, we get
\begin{equation}
\label{bound1 lemma7}
  \frac{1}{2\alpha} N \theta_N^2 r_i^3 \E \xi_i^2 - 1
    = \frac{r_i-r_{i-1}}{r_{i-1}} + \frac{r_i}{2N \theta_N^2
      r_{i-1}^2}
      \leq \frac{r_i}{2N \theta_N^2
        r_{i-1}^2} \leq \frac{1}{2N}
      %\frac{1}{N \theta_N^2(r_i-1)}.
  \end{equation}
  On the other hand,
  \begin{equation}
  \label{bound2 lemma7}
      1-\frac{1}{2\alpha} N \theta_N^2 r_i^3 \E \xi_i^2=  \frac{r_{i-1}-r_{i}}{r_{i-1}} - \frac{r_i}{2N \theta_N^2
      r_{i-1}^2}
      \leq \frac{r_{i-1}-r_{i}}{r_{i}}
      \leq \frac{1}{2N\theta_N^2(r_i-1)},
  \end{equation}
  where we used
  \begin{equation}
  \label{eq:ridiff}
  \frac{r_{i-1} - r_i}{r_i}
  = \frac{-\Delta_N r'(x_*)}{r_i}
  \leq \frac{\Delta_N}{2(r_i-1)}.
  \end{equation}

 Inequalities \eqref{bound1 lemma7} and \eqref{bound2 lemma7} establish the required bound on $|\varepsilon_i|$.
  The $N^{-2/3}$ bound follows from \eqref{eq:rim1}.
\color{black}

\subsubsection{Proof of Lemma \ref{lemma L0} ($\xi_i$ and $L_i$ are sub-gamma)}\label{proof of lemma L0}

Recall the definition \eqref{definitions SG}. Straightforward
calculations using the 
definitions of $\alpha _{i}$ and $\beta _{i}$
and the monotonicity properties of $SG(v,u)$, keeping in mind that $\beta_1=0$, show that
\begin{equation*}
\alpha _{i} \in SG\left( \frac{\alpha}{N\theta _{N}^{2}r_{i}^{2}},0\right) ,\qquad
\beta _{i} \in SG\left( \frac{\alpha m_{i}}{N\theta _{N}^{2}r_{i}r_{i-1}^{2}},%
\frac{\alpha}{N\theta _{N}^{2}r_{i}r_{i-1}}\right).
\end{equation*}%
%where $\alpha$ is the parameter introduced immediately after
%\eqref{tridiag form earlier}.
To bound the moments of $\beta_i$, we use the following lemma. 

\begin{lemma}
\label{Lemma SG Egor}Suppose $X\in SG\left( v,v\right) $ with $v\leq
1/2.$
Then for any $p>2,$%
\[
\left\Vert X\right\Vert _{p}^{2}\leq 8vp^{2}.
\]
\end{lemma}

\begin{proof}[Proof of Lemma \ref{Lemma SG Egor}]

First, closely following the proof of Theorem 2.3 of \cite{boucheron2013concentration},
we obtain the following inequality
\[
\left\Vert X\right\Vert _{p}^{p}\leq p2^{p-1}\left( \left( 2v\right)
^{p/2} \Gamma(p/2)  +\left( 2v\right)^{p} \Gamma(p) \right).
\]
Since for any $x>1,$ $\Gamma \left( x\right) \leq x^{x-1}$ (see
\cite{anderson1997monotoneity}), and $\left( 2v\right) ^{p}\leq \left(
  2v\right) ^{p/2},$ we get
 \[
 \left\Vert X\right\Vert _{p}^{p}\leq \left(
 2^{p}p^{p/2}+2^{3p/2-1}p^{p}\right) v^{p/2}\leq 2^{3p/2}v^{p/2}p^{p}. 
 \]
\end{proof}

Since $m_i/(r_i r_{i-1}^2)< 1/(r_i r_{i-1}) < 1$, we have
$\beta_i \in SG(\alpha/N \theta_N^2, \alpha/N \theta_N^2)$ and so, from \cref{Lemma SG Egor},
\begin{equation} \label{eq:beta-moments}
  \| \beta_i \|_p^2 \leq \frac{8 \alpha p^2}{N \theta_N^2}
                  \lesssim \frac{p^2}{N}.
\end{equation}
In addition,
\begin{eqnarray}
\xi_{i}=\alpha _{i}+\beta _{i} &\in &SG\left( \frac{\alpha}{N\theta _{N}^{2}r_{i}^{2}}+%
\frac{\alpha m_{i}}{N\theta _{N}^{2}r_{i}r_{i-1}^{2}},\frac{\alpha}{N\theta
_{N}^{2}r_{i}r_{i-1}}\right)  \nonumber \\
&\in &SG\left( \frac{2\alpha}{N\theta _{N}^{2}r_{i}^{3}},\frac{\alpha}{N\theta
_{N}^{2}r_{i}^{2}}\right) \equiv SG\left( v_{i},u_{i}\right).
\label{alpha plus beta SG}
\end{eqnarray}%
The latter inclusion follows from the facts that $r_{i-1}>r_{i}$ and $r_{i}+m_{i}=2$.

Next, we use the identity \eqref{L as sum} that expresses $L_i$ as a weighted sum of $\xi_j$. 
For $i=1$, the inclusion $L_i\in SG\left(v_{Li},u_{Li}\right)$ follows from the identity $L_1=\xi_1$ and the observation that $v_{Li}\geq v_i$ and $u_{Li}=u_{i}$. Further, for any $1<i\leq N$ we have from \eqref{L as sum} and (\ref{alpha plus beta SG}), $L_i\in SG(\tilde{v}_{Li},\tilde{u}_{Li})$ with 
%with $v_{i}=\frac{2\alpha}{N\theta
%_{N}^{2}r_{i}^{3}}$ and $u_{i}=\frac{\alpha}{N\theta _{N}^{2}r_{i}^{2}}.$ %
\[
\tilde{v}_{Li}\leq v_{i}+\sum\nolimits_{j=0}^{i-2}\gamma _{i}^{2}...\gamma
_{i-j}^{2}v_{i-j-1}.
\]%
Since $\gamma _{i}$ is increasing in $i,$ this yields
\begin{equation*}
  \tilde{v}_{Li} \leq v_{i}( 1+\gamma _{i}^{2}+\gamma _{i}^{4}+...)
  = \frac{2 \alpha}{N \theta_N^2}\frac{1}{r_i^3(1-\gamma_i^2)}
  < \frac{\alpha}{2N \theta_N^2} \frac{1}{r_i-1}=v_{Li},
\end{equation*}
where we used
\begin{equation}
    \label{eq:neg-gamma-2}
(1-\gamma^2)^{-1}
  = r^2/(4R)  < (r-1)^{-1}
\end{equation}
with $R = r-1$, and $r>1$.
For $\tilde{u}_{Li},$ we have $%
\tilde{u}_{Li}\leq \max_{j\leq i}u_{j}=\alpha/\left( N\theta _{N}^{2}r_{i}^{2}\right)=u_{Li} $
because $r_{i}$ is decreasing in $i$. Hence, the inclusion $L_i\in SG\left(v_{Li},u_{Li}\right)$ follows by the monotonicity of the sub-gamma family.

\subsubsection{Proof of Lemma \ref{Lemma gammas Egor} (products $\gamma_{j:i+1}$ are not too small)}\label{proof of Lemma gammas Egor}

Let $n_0 = N - \left[N^{1/3} \log^3 N\right]$.
Since $x_{n_0+1} = n_0/N\theta_N^2 \geq \theta_N^{-2} - (\log^3 N)/(N^{2/3}
\theta_N^2)$, we have $\epsilon_N = 1 - x_{n_0+1} \leq 2 (\log^3 N)/N^{2/3}$ for large $N$ and $w_N\ll\log^2N$.
Recalling that for $\epsilon \in (0,\tfrac{1}{2})$, we have both
$\log(1-\epsilon) \geq - 2\epsilon$ and
$-\log(1+\epsilon) \geq -\epsilon$ and noting that
$\log \gamma(x) = \log(1-\sqrt{1-x}) - \log(1+ \sqrt{1-x})$, we get
\begin{equation*}
  \log \gamma(x_{n_0+1})
  \geq -3 \sqrt{\epsilon_N}
  \geq -3\sqrt{2} (\log^{3/2}N)/N^{1/3}.
\end{equation*}
Since $\gamma_j$ and $x_j$ are increasing in $j$, we have for large $N$
  \begin{equation*}
    \log \gamma_{j:i+1} \geq (j-i) \log \gamma_{i+1}
    \geq (j-i) \log \gamma(x_{n_0+1})
    \geq -C \log^{-1/2}N \geq \log \tfrac{1}{2}.  \qedhere
  \end{equation*}

\subsubsection{Proof of equation \eqref{eq:Rtildesmall} (about $|\tilde{R}_i|$ being small) from Lemma \ref{Lemma Egor R}}\label{proof of eq 45}

Recall the decomposition
\[
T \bar{\varepsilon}_i = - T\delta_i + T\bar{\varepsilon}_i^{\rm m}
                    + T\bar{\varepsilon}_i^{\rm s}  + T\bar{\varepsilon}_i^{\rm
                      q}.
                      \]
We start from the term $T\bar{\varepsilon}_i^{\rm m}$. 
As explained in the proof of \cref{Lemma Egor R} in the main text, this term can be viewed as a sum $
  \sum_{j=1}^i \gamma_{i:j+1} \bar{R}^{(1)}_{j-1} \beta_j\equiv \sum_{j=1}^iX_j $
of martingale differences, and the Marcinkiewicz-Zygmund-type
inequality of \cite[Theorem 2.1]{rio2009moment} says that
$\| \sum_1^i X_j \|_p^2 \leq (p-1) \sum_1^i \|X_j\|_p^2$.

Appealing to \eqref{eq:beta-moments} we obtain
\begin{equation} \label{eq:TbetaR}
  \| \sum_{j=1}^i \gamma_{i:j+1} \bar{R}^{(1)}_{j-1} \beta_j \|_p^2
    \leq \frac{p-1}{1-\gamma_i^2} \max_{j \leq i} \| \bar{R}^{(1)}_{j-1}
    \|_p^2 \| \beta_j \|^2_p
    \lesssim \frac{\alpha p^3}{N(r_i-1)} \max_{j \leq i} \| \bar{R}^{(1)}_{j-1} \|_p^2.
\end{equation}
To obtain the latter inequality, we also used
\begin{equation}
    \label{eq:gambd}
\frac{1}{1-\gamma_i^2}\leq\frac{1}{1-\gamma_i}=\frac{r_i}{r_i-m_i}=\frac{r_i}{2(r_i-1)}<\frac{1}{r_i-1}.
\end{equation}

By construction $|\bar R_i| \leq N^{-1/3}/2$ and for $N$ sufficiently large
$(1-N^{-1/3}/2)^{-1} \leq 4/3$ so that $\bar R_i^{(1)} \leq
N^{-1/3}$. For the Winsorized ``martingale'' terms, then, from
\eqref{eq:TbetaR} 
%and \eqref{eq:rim1}
\begin{equation} \label{eq:Tebm}
  \| T \bar \varepsilon_i^{\rm m} \|_p^2
  \lesssim \frac{\alpha p^3}{N^{5/3}(r_i-1)}
  \lesssim \frac{\alpha p^3}{N^{4/3}}.
\end{equation}
For the last inequality, we used
\begin{equation}
\label{eq:rim1}
\frac{1}{r_i-1}\leq\frac{1}{r_{N+1}-1}=(1-\theta_N^{-2})^{-1/2}=\frac{N^{1/3}}{\sqrt{2w_N}}(1+o(1)).
\end{equation}
For any random variable $\Pr( |X| \geq e \|X\|_{2 \log N}) \leq e^{-2
  \log N} = N^{-2}$, so taking the union bound over $i=1, \ldots, N$,
we find that there exists $C>0$ such that for all sufficiently large
$N$, with probability at least $1 - 1/N$, 
\begin{equation*}
  \max_i |T\bar \varepsilon_i^{\rm m} |
  \leq CN^{-2/3}\log^{3/2}N.%=O(N^{-2/3} \log^{3/2}N).
\end{equation*}

Turning to $T \delta_i$, 
one finds that $1/r(x)$ is convex increasing, with
$(1/r)' = - r'/r^2$, and so for $i \geq 2$
\begin{equation}
  \label{eq:deltabds}
  0 < \delta_i = m_i \Big( \frac{1}{r_i} - \frac{1}{r_{i-1}} \Big)
  \leq  \Delta_N\frac{-\gamma r'}{r}(x_i)
%  = \frac{\Delta_N}{2}\frac{\gamma_i}{r_i(r_i-1)}
  =\frac{\Delta_N}{2}\frac{m_i}{r_i^2(r_i-1)}
  \leq \frac{\Delta_N}{2(r_i-1)}.
\end{equation}
Furthermore, $\delta_i$ is increasing, so that from \eqref{eq:T-bound} followed by
\eqref{eq:deltabds}, then \eqref{eq:gambd} and \eqref{eq:rim1},
\begin{equation}
  \label{eq:Tdelta}
  T \delta_i
  \leq \frac{\delta_i}{1-\gamma_i}
  \leq \frac{1}{2N \theta_N^2 (r_i-1)^2}
  = O(N^{-1/3} w_N^{-1}).
\end{equation}

For the ``small'' term $\bar \varepsilon_i^{\rm s}$,  we have
$|\bar \varepsilon_i^{\rm s}| \leq (1-\bar 
R_{i-1})^{-1}(\gamma_i |\bar R_{i-1}|^3 + \delta_i |\bar R_{i-1}|)
\leq \frac{4}{3}[(8N)^{-1} + \delta_i N^{-1/3}]$. 
From \eqref{eq:T-bound} and \eqref{eq:gambd}, we have
\begin{equation} \label{eq:Tebs}
  |T\bar \varepsilon_i^{\rm s} |
   \leq \frac{1}{1-\gamma_i} \Big[ \frac{1}{6N} +
   \frac{4\delta_i}{3N^{1/3}} \Big]
   \leq \frac{1}{N(r_i-1)}
\end{equation}
for large $N$, since from \eqref{eq:deltabds} and \eqref{eq:rim1}
$\delta_iN^{-1/3} \leq \Delta_N[ N^{1/3} (r_i-1)]^{-1} = O(N^{-1}
w_N^{-1/2})$.   In particular, again from \eqref{eq:rim1},
$\max_{i \leq N} |T\bar \varepsilon_i^{\rm s} | = O(N^{-2/3} w_N^{-1/2})$. 

Finally, since the ``quadratic'' term
$\bar \varepsilon_i^{\rm q} \leq \bar R_{i-1}^2$,
we have from \eqref{eq:T-bound}
\begin{equation} \label{eq:Tebq}
  |T\bar \varepsilon_i^{\rm q} |
  \leq \frac{1}{4 N^{2/3}(1-\gamma_i)}
  = O(N^{-1/3} w_N^{-1/2}).
\end{equation}

Taking into account the representation $\tilde{R}_i  = L_i + T \bar{\varepsilon}_i$ and result \eqref{eq:result of lemma} of Lemma \ref{Lemma Egor},
we see that there exists $C>0$ such that on an event of probability
at least $1-o(1),$%
\[
  \max_{i\leq N}\vert \tilde{R}_{i}\vert
  \leq C N^{-1/3}/\sqrt{\log \log N},
\]%
which establishes \eqref{eq:Rtildesmall} of Lemma \ref{Lemma Egor R}.

%\begin{remark}
%In both Lemma \ref{Lemma Egor} and Lemma \ref{Lemma Egor R}, we only
%guarantee the probability $1-o(1)$ of the form $1-O( \log ^{-k}N)$,
%which is a rather slow convergence rate. 
%\end{remark}

\subsubsection{Proof of Lemma \ref{Lemma star} (bounds on components of $\sum(R_i+R_i^2/2)$)}\label{proof of lemma star}

That $| \sum_1^N T \varepsilon_i^{\rm s} | \leq  N^{-1} \sum_1^N
(r_i-1)^{-1} =  O(1)$ on a set $\mathcal{R}_N$ of probability $1 -
o(1)$ follows already from \eqref{eq:Tebs}
and
\eqref{eq:ri-sums}. 

To show that $\sum_1^N T \varepsilon_i^{\rm m} = O_\Pr(1)$ it would in
principle suffice to show that $\| T \varepsilon_i^{\rm m} \|_1$ is
summable.
Lemma \ref{Lemma Egor R} controls $\{ R_i \}$ on $\mathcal{R}_N$, but
this is not enough for moments of $R_i$ since $\Pr( \mathcal{R}_N^c)
\approx O(\log^{-k} N)$.
Instead we define a recursive mutilation $\hat{R}_i$ for $R_i$ with
two properties: (a) that $\hat{R}_i = R_i$ for $i=1,\ldots,N$
on $\mathcal{R}_N$, and
(b) for all such $i$,
\begin{equation}
  \label{eq:Rh4}
  \| \hat{R}_i \|_4 \leq 4 \alpha^{1/2} N^{-1/2} (r_i-1)^{-1/2}.
\end{equation}

To define $\{ \hat{R}_i \}$, recall that the process $\{ R_i \}$
satisfies
\begin{align}
  R_i & = L_i - T \delta_i + T\vem_i + T\ves_i + T\veq_i.  \notag \\
\intertext{We set $\hat{R}_1 = L_1 = R_1$. For $i \geq 2,$ given
$\{ \alpha_i,\beta_i,\gamma_i,\delta_i \}_1^i$ and $\{ R_j,
\hat{R}_j\}_1^{i-1}$, set}
  \hat{R}_i & = L_i - T \delta_i + T\hvem_i
              + \phi_{\frac{2}{N(r_i-1)}} (T\ves_i)
              + T\hveq_i, \label{eq:Rhatdef}
\end{align}
with the modifications
\begin{equation*}
    \hvem = \beta_i \hat{R}_{i-1}/[1-\phi_{1/2}(\hat{R}_{i-1})], \qquad
  \hveq = \gamma_i \phi_{N^{-1/3}} (\hat{R}_{i-1})\hat{R}_{i-1}.
\end{equation*}

Property (a) is verified by chasing definitions:
let $\mathcal{H}_i = \{ \hvem_j = \vem_j, \hveq_j = \veq_j \text{ for }
j = 1, \ldots i\} \cap \{ |T\ves_i| \leq 2/N(r_i-1) \}$.
By definition event $\mathcal{H}_i$ implies $\hat{R}_j = R_j$ for $j
\leq i$.
One checks by induction that $\mathcal{H}_i$ holds on $\mathcal{R}_N$
($\mathcal{H}_1$ is always true), and so $\{ \hat{R}_j\}_1^N
= \{ R_j \}_1^N$.

We turn to property (b). 
By Lemma \ref{lemma L0} and Theorem 2.3 of \cite{boucheron2013concentration},
\[
\left\Vert L_{i}\right\Vert_{4}^{4}\leq 2\left(\frac{4\alpha}{N\theta_N^{2}(r_{i}-1)}\right)^{2}+4!\left(\frac{4\alpha}{N\theta_N^{2}r_i^{2}}\right)^{4}\leq\frac{33\alpha^2}{N^2\theta^{4}_n(r_i-1)^{2}}
\]
for all sufficiently large $N$. Hence 
\[
\left\Vert L_{i}\right\Vert_{4}\leq \frac{3\alpha^{1/2}}{\theta_N \sqrt{N(r_{i}-1)}}.
\]

To bound the norm of $T \hvem_i$, observe that
since $\beta _{i}$ and $\hat{R}_{i-1}$ are independent,
$\| \hvem_i \|_p \leq 2 \| \hat{R}_{i-1} \|_p \| \beta_i \|_p$.
Since $\beta _{i}\in SG\left( \frac{\alpha}{%
\theta _{N}^{2}N},\frac{\alpha}{\theta _{N}^{2}N}\right) ,$ we have $\left\Vert
\beta _{i}\right\Vert _{p}^{2}\leq \frac{8\alpha p^{2}}{\theta
_{N}^{2}N}$ by Lemma \ref{Lemma SG Egor}, and so
\[
\| \hvem_i \|_p
\leq \frac{
2^{5/2}\alpha^{1/2}p}{\theta _{N}\sqrt{N}}\Vert \hat{R}_{i-1}\Vert _{p}. 
\]
By the Marcinkiewicz-Zygmund type inequality (see Theorem 2.1 from Rio (2009)), and using $\left( 1-\gamma _{i}^{2}\right)
^{-1}\leq \left( r_{i}-1\right) ^{-1}=o\left( N^{1/3}\right) ,$%
\begin{equation}
  \| T \hvem_i \|_4 
  \leq 3^{1/2}\frac{2^{9/2}\alpha^{1/2}}{\theta _{N}%
\sqrt{N\left( 1-\gamma _{i}^{2}\right) }}\max_{j\leq i-1}\Vert \hat{R}%
_{j}\Vert _{4}=o\left( N^{-1/3}\right) \max_{j\leq i-1}\Vert 
\hat{R}_{j}\Vert _{4}.  \label{That14}
\end{equation}%

Next, by definition, 
\[
 \big\| \phi _{\frac{2}{N\left( r_{i}-1\right) }%
}\left( T \ves_i \right)  \big\|_{4}\leq \frac{2}{N\left( r_{i}-1\right) }.
\]
And finally, we have $\Vert \hveq_i \Vert _{4}\leq
N^{-1/3}\Vert \hat{R}_{i-1}\Vert _{4}.$ Therefore, from \eqref{eq:T-bound}
\[
\Vert T \hveq_i \Vert _{4}\leq N^{-1/3}\frac{1}{1-\gamma _{i}}%
\max_{j\leq i-1}\Vert \hat{R}_{j}\Vert _{4}=o(1)\max_{j\leq
i-1}\Vert \hat{R}_{j}\Vert _{4}.
\]%
Along with \eqref{eq:Tdelta}, all this sums up to%
\begin{eqnarray*}
\Vert \hat{R}_{i}\Vert _{4} &\leq &\frac{3\alpha^{1/2}(1+o(1))}{\sqrt{N\left(
r_{i}-1\right) }}+\frac{1}{2N\left( r_{i}-1\right)
^{2}} \\
&&+o( N^{-1/3}) \max_{j\leq i-1}\Vert \hat{R}_{j}\Vert
_{4}+\frac{2}{N\left( r_{i}-1\right) }+o(1)\max_{j\leq i-1}\Vert \hat{R}%
_{j}\Vert_{4},
\end{eqnarray*}%
which implies%
\[
\Vert \hat{R}_{i}\Vert _{4}\leq \frac{1}{\sqrt{N(
r_{i}-1) }}( 3\alpha^{1/2}+o(1)) +o(1)\max_{j\leq i-1}\Vert \hat{R}%
_{j}\Vert _{4}.
\]%
Hence, by induction (and for sufficiently large $N$), we conclude that
\eqref{eq:Rh4} and property (b) hold.

Returning to $\sum T \vem_i$, the Marcinkiewicz-Zygmund type bound
\eqref{That14} and then \eqref{eq:Rh4} show that
\begin{equation}
  \label{eq:Temh1}
  \| T \hvem_i \|_4
    \leq \frac{C \alpha^{1/2}}{\sqrt{N(1-\gamma_i)}}
    \frac{\alpha^{1/2}}{\sqrt{N(r_i-1)}}
    \leq \frac{C \alpha}{N(r_i-1)}.
\end{equation}
From \eqref{eq:ri-sums},
$\sum \| T \hvem_i \|_1    \leq \sum \| T \hvem_i \|_4
    \leq C \theta_N^2 = O(1).$
Therefore $\sum_1^N T \hvem_i = O_\Pr(1)$ and the same is then true
for $\sum_1^N T \vem_i$ as they are equal on $\mathcal{R}_N$, a set
of probability $1 - o(1)$.

To finally show that $\sum T \veqd_i = O_\Pr(1)$, we first define
$\cveqd_i = \gamma_i(\hat{R}_{i-1}^2 - L_{i-1}^2)$, noting that
$\{ \cveqd_i \}_1^N = \{ \veqd_i \}_1^N$ on $\mathcal{R}_N$.
We first use \eqref{eq:Rhatdef} to bound $\| \hat{R}_i - L_i \|_2$.
Since $\|\hveq_i\|_2\leq \|\hat{R}_{i-1}^2 \|_2= \|
\hat{R}_{i-1}\|_4^2$
  , we have from \eqref{eq:Rh4} and
  \eqref{eq:gambd} that
  \begin{equation*}
    \| T \hveq_i\|_2
    \leq  \frac{1}{1-\gamma_i}
    \max_{j\leq i} \| \hveq_i \|_2
    \leq \frac{16 \alpha}{N(r_i-1)^2}.
  \end{equation*}
Combining this with \eqref{eq:Tdelta}, \eqref{eq:Temh1}, and the
trivial bounds
$\| \phi_{2 N^{-1} (r_i-1)^{-1}}(T \ves_i) \|_2 \leq 2 N^{-1}
(r_i-1)^{-1}$ and $(r_i-1)^{-1} \leq (r_i-1)^{-2}$ we get
\begin{equation*}
  \| \hat{R}_i - L_i \|_2
   \leq C \alpha N^{-1} (r_i-1)^{-2}.
\end{equation*}
Since both $\| \hat{R}_i \|_2$ and $\|L_i\|_2$ are $O(\alpha^{1/2}
N^{-1/2} (r_i-1)^{-1/2})$, we have
\begin{equation*}
  \| \cveqd_i \|_1
    \leq \| \hat{R}_{i-1} - L_{i-1} \|_2 \| \hat{R}_{i-1} + L_{i-1}
    \|_2
    \leq C \alpha^{3/2} N^{-3/2} (r_i-1)^{-5/2},
  \end{equation*}
and from \eqref{eq:T-bound}, we get $\| T \cveqd_i \|_1   
\leq C \alpha^{3/2} N^{-3/2} (r_i-1)^{-7/2}$.
Now appealing to \eqref{eq:ri-sums} with $\beta = 7/2$,
\begin{equation*}
  \sum_1^N \| T \cveqd_i \|_1   
  \leq C \alpha^{3/2} N^{-1/2} w_N^{-3/4} N^{1/2}
  = O(w_N^{-3/4}),
\end{equation*}
which establishes that $\sum T \veqd_i = O_\Pr(1)$ and completes the
proof of \eqref{eq:msqd sum}. 

The bound for $\sum R_i^2$  follows from \eqref{eq:Rh4}. Indeed,
\begin{eqnarray*}
\sum\nolimits_{i=1}^{N}\Vert \hat{R}_{i}^{2}\Vert _{1}
&=&\sum\nolimits_{i=1}^{N}\Vert \hat{R}_{i}\Vert _{2}^{2} \ \leq \ 
\sum\nolimits_{i=1}^{N}\Vert \hat{R}_{i}\Vert _{4}^{2} \\
&\leq &\sum\nolimits_{i=1}^{N}\frac{16\alpha}{N\left( r_{i}-1\right) }=O(1).
\end{eqnarray*}%
Therefore, $\sum\nolimits_{i=1}^{N}\hat{R}_{i}^{2}=O_{\Pr}(1)$ and
so is the sum of $R_{i}^{2}$, which establishes \eqref{eq: R2 sum}.

We turn to the proof of \eqref{eq:TeqL} and \eqref{eq:TeqE}.
Using \eqref{eq:gami-to-j}, then $\gamma_{i:j+1} \gamma_j =
\gamma_{i:j}$, and \eqref{eq:gjdef}, we have
\begin{align*}
  \sum_{i=1}^N T \varepsilon_i^{\rm qL}
  & = \sum_{i=1}^N \sum_{j \leq i} \gamma_{i:j+1} \gamma_j L_{j-1}^2
    = \sum_{j=1}^{N-1} (g_{j+1}-1)L_j^2 \\
  & = \bar \xi^\T T^\T G T \bar \xi,
\end{align*}
with $\bar \xi = (\xi_1, \ldots, \xi_{N-1})^\T$ and
matrices $G = \diag(g_2-1, \ldots, g_N-1)$, and $T = T_{N-1}$. 

The rescaled vector $\mathbf{x} = (\xi_1/\sigma_1,...,\xi_{N-1}/\sigma_{N-1})^\T, \sigma_i^2 = \E
\xi_i^2$ has independent components with common variance.
We use a variance bound for quadratic forms in such variables
\cite[Lemma B.26]{basi09}, namely
$\Var \, \mathbf{x}^T A \mathbf{x} \leq C \nu_4 \| A \|_{\rm HS}^2$,
where in our case $\nu_4 = \max \E x_i^4 \leq C_2/c_2^2$ by Lemma
\ref{lemma moments of ksi} and $A = DT^\T G T D$,
with $D = \diag(\sigma_i)$.
We have
\begin{equation*}
  \Var [\sum_i T(\varepsilon_i^{\rm qL} - \varepsilon_i^{\rm qE})]
  = \Var (\bar \xi^\T T^\T G T \bar \xi)
  = \Var (\mathbf{x}^\T A \mathbf{x})
  \lesssim \| A \|_{\rm HS}^2.
\end{equation*}
Again by Lemma \ref{lemma moments of ksi}, $\max \sigma_i^2 \leq C_1
\alpha/N$, and so
\begin{equation*}
  \| A \|_{\rm HS} \leq C_1 \alpha N^{-1} \| T \|_{\rm op} \| GT
  \|_{\rm HS}.
\end{equation*}

Decomposing $T$ into a sum of sub-diagonal matrices, we have by the
triangle inequality 
\[
    \| T \|_{\rm op} \leq 1 + \gamma_{N-1} + \dots + \gamma_{N-1} \dots \gamma_{2} \leq \frac{1}{1 - \gamma_{N-1}} = O(N^{1/3} w_N^{-1/2})\, .
\]
We also have,
\[
  \| G T \|_{\mathrm{HS}}^{2}
  = \sum_{i = 1}^{N-1} (g_{i+1} - 1)^2 (1 + \gamma_{i}^2 + \dots
  + \gamma_{i}^{2} \dots \gamma_{2}^{2})
  \leq \sum_{i = 1}^{N-1} (g_{i+1} - 1)^{2} \frac{1}{1 - \gamma_{i}^2} \, .
\]
By Lemma~\ref{lemma si}, for large enough \(N\), each \(g_i \geq 1 \)
satisfies $g_i-1 < (r_i-1)^{-1}$. 
Then from
\eqref{eq:gambd} and \eqref{eq:ri-sums}, 
\begin{align*}
  \frac{1}{N} \| G T \|_{\mathrm{HS}}^{2}
  < \sum_{i = 2}^{N}\frac{1}{N(r_i-1)^2(1 -  \gamma_{i}^2)} 
  < \sum_{i = 2}^{N} \frac{1}{N(r_{i} - 1)^3}
  = O(N^{1/3} w_N^{-1/2}) \, .
\end{align*}
Summing up, we get
\begin{equation*}
  \|A\|_{\rm HS}
  \leq C_1 \alpha N^{-1} O(N^{1/3} w_N^{-1/2}) O(N^{2/3} w_N^{-1/4})
  = O(w_N^{-3/4}),
\end{equation*}
which suffices to establish \eqref{eq:TeqL}.

It remains to show \eqref{eq:TeqE}.

\begin{lemma}
\label{lemma delta adjusted}There exists $C>0$ such that for all $1\leq
i\leq N,$%
\[
  \left\vert \gamma _{i} \E L_{i-1}^{2} - \alpha \delta _{i} \right\vert
  <\frac{C \alpha}{N^{2}\left( r_{i}-1\right) ^{4}} \,. 
\]
\end{lemma}
\noindent The lemma implies \eqref{eq:TeqE}, for using \eqref{eq:T-bound}, then
\eqref{eq:gambd} and \eqref{eq:ri-sums},  we have
\begin{equation*}
  \sum_i T \varepsilon_i^{qE} - \alpha T \delta_i
    \leq \sum_{i=1}^N \frac{1}{1-\gamma_i}\frac{C \alpha}{N^2
      (r_i-1)^4}
    \lesssim \frac{1}{N^2} \sum_{i=1}^N \frac{1}{(r_i-1)^5}
    = O(w_N^{-3/2}).
\end{equation*}

%\subsubsection{Proof of Lemma \ref{lemma delta adjusted}}\label{proof
%  of lemma delta adjusted}
\textit{Proof of Lemma~\ref{lemma delta adjusted}.}
Recall that $L_i = \xi_i + \gamma_i L_{i-1}$. Since $\xi_i$ and
$L_{i-1}$ are independent, $\E L_i^2$ satisfies the recursion
\begin{equation}  \label{eq:L2-recursion}
  \E L_i^2 = \E \xi_i^2 + \gamma_i^2 \E L_{i-1}^2, \qquad \qquad i \geq 1.
\end{equation}
The idea is to use $\E \xi_i^2 \approx (\alpha\Delta_N) (2/r_i^3)$
to show that an approximate solution
is  $\E L_i^2 \approx s_i$, with
\begin{equation*}
  s_i = \frac{\alpha}{2} \frac{\Delta_N}{r_i(r_i-1)},
\end{equation*}
and to use the equality, which we prove at the end of this section, 
\begin{equation}
\label{eq:deltaapprox}
  \delta_i = \gamma_i s_i /\alpha+O(N^{-2}(r_i-1)^{-3}).  
\end{equation}
%recall from \eqref{eq:deltaapprox} that $\delta_i = \gamma_i s_i /\alpha+O(N^{-2}(r_i-1)^{-3})$.
The other key ingredient is the identity
\begin{equation*}
  \frac{2}{r_i^3}
    = \frac{1}{2r_i(r_i-1)} - \frac{\gamma_i^2}{2r_{i}(r_{i}-1)}.
\end{equation*}
which, recalling \eqref{eq:neg-gamma-2}, 
follows from $1 - \gamma^2 = 4R/r^2$ and $R = r-1$.

In detail, from Lemma~\ref{lemma variance} we have
\begin{align*}
  \E\xi _{i}^2
  & = (\alpha\Delta_N) (2/r_i^3) + O(N^{-2}(r_i-1)^{-1}) \\
  & = s_i - \gamma_i^2 s_i  + O(N^{-2}(r_i-1)^{-1}) \\
  & = s_i - \gamma_i^2 s_{i-1} + \eta_i,
  \qquad \qquad \qquad
  |\eta_i| = O(N^{-2} (r_i-1)^{-3}).
\end{align*}
since the function $\rho(x) = 1/r(x)(r(x)-1)$ with $r(x)=1+\sqrt{1-x}$, $x\in(0,1)$  has $|\rho'(x)| \leq (r(x)-1)^{-3}$, and so
$0 \leq s_i - s_{i-1} \leq \Delta_N^{2} (r_i-1)^{-3}$.
 
Putting this into \eqref{eq:L2-recursion}, we obtain
\begin{equation*}
  \E L_i^2 - s_i = \gamma_i^2(\E L_{i-1}^2 - s_{i-1}) + \eta_i,
\end{equation*}
whose solution, from \eqref{eq:recurr}, is $\E L_i^2 - s_i = (\tilde{T} \eta)_i$, where the linear operator $\tilde{T}$ is defined as $T$ but with $\gamma$ replaced by $\gamma^2$.
Now bound \eqref{eq:neg-gamma-2} and an analogue of bound \eqref{eq:T-bound} imply that
\begin{equation*}
  | \E L_i^2 - s_i|
  = |\tilde{T} \eta_i|
  \leq \frac{1}{1-\gamma_i^2} \max |\eta_i|
  = O( N^{-2} (r_i-1)^{-4}).
\end{equation*}
To complete the proof, simply write
\begin{equation*}
  \alpha^{-1}\gamma_i \E L_{i-1}^2 - \delta_i
  = \alpha^{-1} \gamma_i(s_{i-1} - s_i) + O( N^{-2} (r_{i-1}-1)^{-4})
  = O( N^{-2} (r_i-1)^{-4}).
\end{equation*}

It remains to establish \eqref{eq:deltaapprox}. Similarly to \eqref{eq:deltabds}, we have 
\begin{equation}
    \label{eq:deltabds1}
    \delta_i\geq m_i\Delta_N\frac{-r'}{r^2}(x_{i-1})=\frac{\Delta_N}{2}\frac{m_i}{r_{i-1}^2(r_{i-1}-1)}.
\end{equation}
Since $f(r)=r^{-2}(r-1)^{-1}$ is convex decreasing in $r\in (1,2)$ with $\diff f/\diff r=-(3r-2)/r^3(r-1)^2$, we have
\[
\frac{1}{r_i^2(r_i-1)}-\frac{1}{r_{i-1}^2(r_{i-1}-1)}\leq \frac{3r_i-2}{r_i^2(r_i-1)^2}\frac{r_{i-1}-r_i}{r_i}\leq 
\frac{2\Delta_N}{(r_i-1)^3},
\]
where we have used \eqref{eq:ridiff}.
Hence, the difference between the upper and lower bounds in \eqref{eq:deltabds} and \eqref{eq:deltabds1} is no larger than $\Delta_N^2/(r_i-1)^3$ and
\begin{equation}
    %\label{eq:deltaapprox}
    \delta_i=\frac{\Delta_N}{2}\frac{\gamma_i}{r_i(r_i-1)}+O(N^{-2}(r_i-1)^{-3}).\qquad\qedsymbol
\end{equation}

\subsubsection{Proof of Lemma~\ref{lemma Tdeltai} (asymptotics of $\sum T\delta_i$)}\label{proof lemma Tdeltai}

Resummation \eqref{eq:resummation} gives
$  \sum_1^N T \delta_i = \sum_1^N g_{i+1} \delta_i=\sum_2^Ng_{i+1}\delta_i$.
From \eqref{eq:deltabds}, \eqref{eq:deltabds1} and Lemma \ref{lemma si} we have bounds
\begin{align*}
 & d(x_{i-1})\Delta_N \leq \delta_i \leq d(x_i) \Delta_N, \qquad\text{for }2\leq i\leq N,\\
 & g_i \leq g(x_i) \big(1+O(w_N^{-3/2})\big),\qquad \text{   for }2\leq i\leq N, \\
 & g_i \geq  g(x_i) \big(1+O(w_N^{-3/2})\big),\qquad \text{   for }2\leq i\leq N-N^{1/3},
 \end{align*}
 along with the definitions
 \[
 d(x) = d = \frac{\gamma}{2r(r-1)}, \qquad
  g(x) = g = \frac{r}{2(r-1)}, \qquad r=1+\sqrt{1-x}.
\]
Combining these, and using monotonicity of $\gamma$ and $r$, we get
\begin{equation*}
  (gd)(x_{i-1}) \Delta_N \big(1+O(w_N^{-3/2})\big)
  \leq g_{i+1} \delta_i
  \leq (gd)(x_{i+1}) \Delta_N \big(1+O(w_N^{-3/2})\big)
\end{equation*}
for $2 \leq i \leq N - N^{1/3}$ in the lower bound and
$2 \leq i \leq N$ in the upper.

Since $1-\gamma = 2(r-1)/r$, we can decompose
\begin{equation*}
  gd = \frac{\gamma}{4(r-1)^2}
  = \frac{1}{4(r-1)^2} - \frac{1}{2r(r-1)}
  = f_1 - f_2,
\end{equation*}
say, and then observe that $f_2(x) \asymp (1-x)^{-1/2}$ is integrable
for $x \in [0,1]$, so that the sums $\sum f_2(x_i) \Delta_N = O(1)$
can be ignored.
On the other hand, the integral of $f_1(x)$ is
\begin{equation*}
  I(x_a,x_b) = \frac{1}{4} \int_{x_a}^{x_b} \frac{\diff x}{1-x}
   = \frac{1}{4} \log \Big(\frac{1-x_a}{1-x_b} \Big).
\end{equation*}

Let $x^{(0)} = a_0 \Delta_N, x^{(1)} = (N-[N^{1/3}]+a_1)\Delta_N$, and
$x^{(2)} = (N+a_2)\Delta_N$ -- the right choices of fixed small integers $a_i$
legitimize bounds like \eqref{eq:riemann}, but make negligible contributions. 
Indeed, $x^{(0)} = O(N^{-1})$, and from
\eqref{eq:thetam2},
\begin{equation*}
  1-x^{(1)} = (2w_N+1)N^{-2/3} \big(1+O(N^{-1/3})\big), \qquad
  1-x^{(2)} = 2w_N N^{-2/3} \big(1+O(N^{-1/3})\big),
\end{equation*}
so that
\begin{equation*}
  I(x^{(0)},x^{(1)}) = \frac{1}{6} \log N + O(\log w_N).  
\end{equation*}
The remaining contribution is negligible:
$I(x^{(1)},x^{(2)}) =O(w_N^{-1}) +O(N^{-1/3})$,
and so to finish
\begin{equation*}
  \sum T \delta_i
  = \sum g_{i+1} \delta_i
  = [I(x^{(0)},x^{(1)}) +O(1)][1+O(w_N^{-3/2})]
  = \frac{1}{6} \log N + O(\log \log N).  
\end{equation*}
\color{black}

\subsubsection{Proof of Lemma~\ref{lem:variance-bound} (a variance bound on linear spectral statistics of GUE)}
\label{proof:lem:variance-bound}
  With determinantal structure (such as GUE), we have [\cite{tracy1998correlation}, (1.2)]
\[
R_k(x_1,...,x_k)=\det (K_N(x_i,x_j))_{i,j=1,...,k},
\]
where the kernel
$K_N(x,y)=\sum\nolimits_{k=0}^{N-1}\phi_k(x)\phi_k(y)$
% has the representation 
% \[
% K_N(x,y)=\sum\nolimits_{k=0}^{N-1}\phi_k(x)\phi_k(y)
% \]
with $\left\{ \phi_k(x)\right\}$ obtained by orthonormalizing $\big\{
  x^{k} e^{-Nx^2/4}\big\}$.
% Note that for unscaled GUE, one would orthonormalize functions
% $\big\{ x^{k} e^{-x^2/4}\big\}$. 
%  Let $s_N(x)=K_N(x,x)$.
We then have in particular
\begin{equation}
    \label{Iain2}
    R_1(x)=K_N(x,x), \qquad R_2(x,y)=R_1(x)R_2(y)-K_N^{2}(x,y).
\end{equation}
Furthermore, let $J=\E\left[N^{-1}\sum\nolimits_{i=1}^N f(l_i)\right]^{2}$. Expanding, we have
\begin{align*}
    J=N^{-2}\E\left[\sum\nolimits_{i=1}^N f^{2}(l_i)\right]+N^{-2}\E\left[\sum\nolimits_{i \neq j} f(l_i)f(l_j)\right]\\
    =N^{-1}\E f^{2}(l_1)+N^{-2}\times N(N-1)\E\left[f(l_1)f(l_2)\right].
\end{align*}
Now apply \eqref{Iain1} and then \eqref{Iain2} to get
\begin{align*}
    J&=N^{-2}\int f^{2}(x)R_1(x)\mathrm{d}x+N^{-2}\iint f(x)f(y)R_2(x,y)\mathrm{d}x\mathrm{d}y\\
    &=N^{-1}\int f^{2}(x)\rho_N(x)\mathrm{d}x+\left[\int f(x)\rho_N(x)\mathrm{d}x\right]^{2}-N^{-2}\iint f(x)f(y)K_N^{2}(x,y)\mathrm{d}x\mathrm{d}y.
\end{align*}
Dropping the last term and recalling \eqref{eq:linear-stat}, we obtain
\eqref{variance bound} which establishes \cref{lem:variance-bound}. 

\subsection{Discussion of edge bounds for one-point functions \eqref{eq:onept-bds}}
\label{sec:edge-bounds-one}

Let $ H_N(x) = e^{x^2} \left(-\frac{\diff}{\diff x}\right)^{N}
e^{-x^2} $ be the Hermite polynomials. The corresponding orthonormal
Hermite functions,
\[
\varphi_N(x) = c_N e^{-x^2/2} H_N(x), \qquad    c_N = (2^NN!
\sqrt{\pi})^{-1/2},
\]
are even/odd as $N$ is even/odd. Consequently 
$I_N = \int \varphi_N(x) \diff x$ vanishes for $N$ odd, and from a
calculation with generating functions, or \cite{gradshteyn2007} 7.373.2,
\begin{equation}
  \label{eq:I2m}
  I_{2m} = \int \varphi_{2m}
         = c_{2m} \sqrt{2\pi} \frac{(2m)!}{m!}
         \sim \frac{\sqrt{2}}{ m^{1/4}}
\end{equation}
as $m \to \infty$, from Stirling's formula.

Let $\varepsilon(x) = \frac{1}{2} \operatorname{sgn}(x)
$ and
$(\varepsilon \varphi)(x) = \int \varepsilon(x-y)\varphi(y) \diff y
= \tfrac{1}{2} \int_{-\infty}^x
\varphi - \tfrac{1}{2} \int_x^\infty \varphi$.
Let $\chi_N^{\rm e} = 1$ if $N$ is even and $0$ otherwise.
The one-point functions, scaled to bulk supported on $[-\sqrt{2N},
\sqrt{2N}]$,  are then given by
\cite[][(5.2.16),(6.3.2),(6.3.5),(6.4.3)]{meht91} or
\cite[][(6.2.10), (7.2.22), (7.2.27--28), (7.2.32)]{meht04})
\begin{align*}
  \sigma_{N,1}(x)
  & = \sum_{j=0}^{N-1} \varphi_j^2(x),  \\
  \sigma_{N,2}(x)
  & =  \sigma_{N,2}(x) + (N/2)^{1/2} \varphi_{N-1}(x) (\varepsilon
    \varphi_N)(x) + I_{N-1}^{-1} \varphi_{N-1}(x) \chi_{N-1}^{\rm e}.
\end{align*}
Note these forms have total mass $N$. To recover the forms of interest
to us, on scale $[-2,2]$ with total mass $1$, we use
\begin{equation}
  \label{eq:scaling}
  \rho_{N,\alpha}(y) = \frac{1}{\sqrt{2N}} \sigma_{N,\alpha}(\sqrt{N/2} \,  y).
%  \rho_{N,j}(y) = \frac{1}{\sqrt{2N}} \sigma_{N,j}(\sqrt{\frac{N}{2}} y).
\end{equation}

Following \cite{trwi96}, introduce
$  \varphi = (N/2)^{1/4} \varphi_N$ and 
$ \psi = (N/2)^{1/4} \varphi_{N-1}. $ We have a useful integral representation (see e.g. equation (57) in \cite{trwi96})
\begin{equation}
\label{eq:integral rep}
\sigma_{N,1}(x)=2\int_{0}^{\infty}\varphi(x+z)\psi(x+z)dz.
\end{equation}
Further,
% \begin{equation*}
%   \varphi = (N/2)^{1/4} \varphi_N, \qquad
%   \psi = (N/2)^{1/4} \varphi_{N-1},
% \end{equation*}
observe from \eqref{eq:I2m} that both
$I_\varphi = \int \varphi$ and $I_\psi = \int \psi$ converge to
$\sqrt{2}$ for large even and odd values of $N$ respectively.
We get
% Since $(\varepsilon \varphi)(x) = \tfrac{1}{2} \int_{-\infty}^x
% \varphi - \tfrac{1}{2} \int_x^\infty \varphi$, we get
\begin{equation}  \label{eq:sigdiff}
  \sigma_{N,2}(x) - \sigma_{N,1}(x)
  = \psi(x) \left[ \tfrac{1}{2} I_\varphi % \chi_N^{\rm e}
    - \int_x^\infty \varphi
    + I_\psi^{-1} \chi_{N-1}^{\rm e} \right].
\end{equation}

Now we turn to bounds for scaled Hermite functions near the bulk
edge. Set $\tau_N =  N^{-1/6}/\sqrt{2}$ and
% with $x = \sqrt{2N} + s\tau_N$,
define
\begin{equation*}
  \varphi_\tau (s) = \varphi(x), \qquad
  \psi_\tau(s) = \psi(x), \qquad
  x = \sqrt{2N} + s \tau_N
\end{equation*}
The following bounds are essentially established in
\cite[p. 403]{olve74} and \cite{johnstone2012fast}

\begin{proposition} \label{prop:Hfunc-bds}
  Fix $0<\varepsilon<2/3$. Then for large $N$, uniformly in the indicated ranges
  \begin{equation*}
    \tau_N \varphi_\tau(s),  \ \tau_N \psi_\tau(s)
    =
    \begin{cases}
      O(e^{-s})   &  s \geq 0 \\
      O \big( (1+|s|)^{-1/4} \big)  & -N^{2/3-\epsilon} < s \leq 0.
    \end{cases}
  \end{equation*}
\end{proposition}
We discuss the proof below. Taking it as given for now, observe then
that
\[
\sigma_{N,1}(x)=2\tau_N\int_{s}^{\infty}\varphi_\tau(y)\psi_\tau(y)dy.
\]
Combining this with Proposition \ref{prop:Hfunc-bds}, we obtain
\[
\sigma_{N,1}(x)
   =
     \begin{cases}
       O(\tau_N^{-1} e^{-2s})  & s \geq 0 \\
       O \big(\tau_N^{-1} (1 + |s|)^{1/2} \big) &
            -N^{2/3-\epsilon} < s \leq 0.
     \end{cases}
\]
From \eqref{eq:scaling}, we have $\rho_{N,\alpha}(2+s N^{-2/3})
= (2N)^{-1/2} \sigma_{N,\alpha}(\sqrt{2N} + s \tau_N)$ and
since $(2N)^{-1/2} \tau_N^{-1} = N^{-1/3}$, the claim
\eqref{eq:onept-bds} for $\alpha=1$ follows from this.

For $\alpha=2$, observe that
\begin{equation*}
  \int_x^\infty \varphi(x') \diff x'
  = \int_s^\infty \tau_N \varphi_\tau(s') \diff s'
  =
  \begin{cases}
    O( e^{-s})  & s \geq 0  \\
    O \big((1+|s|)^{3/4} \big)  & -N^{2/3-\epsilon} < s \leq 0.
  \end{cases}
\end{equation*}
Combining this with Proposition \ref{prop:Hfunc-bds} applied to
$\psi_\tau(s)$ in \eqref{eq:sigdiff}, we arrive at
\begin{equation*}
(\sigma_{N,2} - \sigma_{N,1})(x) 
   = - \psi_\tau(s) \left[ \int_s^\infty \tau_N \varphi_\tau + O(1)
     \right] 
   =
     \begin{cases}
       O(\tau_N^{-1} e^{-2s})  & s \geq 0 \\
       O \big(\tau_N^{-1} (1 + |s|)^{1/2} \big) &
            -N^{2/3-\epsilon} < s \leq 0.
     \end{cases}
\end{equation*}
% \begin{align*}
% (\sigma_{N,1} - \sigma_{N,2})(x) 
%    & = - \psi_\tau(s) \left[ \int_s^\infty \tau_N \varphi_\tau + O(1)
%      \right]  \\
%    & =
%      \begin{cases}
%        O(\tau_N^{-1})  & s \geq 0 \\
%        O \big(\tau_N^{-1} (1 + |s|)^{1/2} \big) &
%             -N^{2/3-\epsilon} < s \leq 0.
%      \end{cases}
% \end{align*}
This implies that \eqref{eq:onept-bds} holds for $\alpha=2$ as well.

\bigskip

\textit{Discussion of proof of Proposition \ref{prop:Hfunc-bds}.} \
Proposition \ref{prop:Hfunc-bds} is based on analysis of the second
order differential 
equation satisfied by Hermite functions $\varphi_N$ using the
Liouville-Green transform around the turning point at the upper edge.
This is detailed in \cite{erde60} (attributed to Skovgaard), and given
as an example of Theorem 11.3.1. in \cite[][Ex 4.2, 4.3 p 403]{olve74}.
This example was also worked out in detail (for another purpose) in
\cite{johnstone2012fast}, JM12 below. Although the focus there was on $s  > -c$, we
indicate how the analysis also extends to much larger ranges of
negative $s$. We focus here on the bound for $\varphi_N$; for
$\varphi_{N-1}$ it is essentially the same, see JM12.

Rescaling the $x$-axis via $x = \sqrt{2N+1} \xi$,
and setting $w_N(\xi) = \varphi_N(x)$,
the Liouville-Green transform introduces
new independent and dependent variables $\zeta$ and
%$W$ via the equations
$W = \dot{\xi}^{1/2} w_N$.
% \begin{equation*}
%   \zeta \left( \frac{\diff \zeta}{\diff \xi} \right)^2
%   = \xi^2 - 1,  \qquad
%   W = \left( \frac{\diff \zeta}{\diff \xi} \right)^{1/2} w_N.
% \end{equation*}
% More precisely
% where for $\xi \geq 1$ and $\xi \leq 1$ respectively, we set
% \begin{equation*}
%   \frac{2}{3} \zeta^{3/2}(\xi)
%   = \int_1^\xi (t^2-1)^{1/2} \diff t, \qquad \qquad
%   \frac{2}{3} (-\zeta)^{3/2}(\xi)
%   = \int_\xi^1 (1-t^2)^{1/2} \diff t.
% \end{equation*}
The transform $W$ approximately satisfies the Airy equation
$W''(\zeta) = \kappa_N^2 \zeta W(\zeta)$ with $\kappa_N = 2N+1$, and
it is shown that $\varphi_N$ is approximated by the (recessive)
solution $\Ai (\kappa^{2/3} \zeta)$ with explicit error bounds.
Indeed, cf \cite[][(71)]{johnstone2012fast}, with
$r(\xi) = [\dot \zeta(\xi)/\dot \zeta(1)]^{-1/2}$,
\begin{equation}
  \label{eq:first}
  \tau_N \varphi(x)
  = (N/2)^{1/4} \tau_N \varphi_N(x)
    = \sqrt{2} r(\xi) \{ \Ai(\kappa_N^{2/3} \xi) +
    \epsilon_2(\xi,\kappa_N) \}.
\end{equation}

The function $\zeta(\xi)$ is increasing and $C^2$ on $(0,\infty)$
\cite[][p 391]{olve74}, with $\ddot{\zeta}(\xi)$ non-negative and
bounded. The arguments leading to (78) and (85) in JM12 show that for
$|s| \geq N^{2/3 - \epsilon}$,
\begin{equation}
  \label{eq:rkappa}
  r(\xi) \leq 1 + O(N^{-\epsilon}), \qquad \qquad
  \kappa^{2/3} \zeta = s(1+O(N^{-\epsilon})).
\end{equation}
To describe error bound even in the oscillatory region of $\Ai$, \cite{olve74}
introduces continuous and positive functions $E \geq 1$ and $M \leq 1$
such that $|\Ai(x)| \leq M(x)/E(x)$ and satisfying 
\begin{equation}  \label{eq:EMbd} 
  E(x) \sim \sqrt{2} e^{(2/3)x^{3/2}}, \qquad \qquad 
  M(x) \sim \pi^{-1/2} (1+|x|)^{-1/4},
\end{equation}
% \begin{align}
%   E(x) & \sim \sqrt{2} e^{(2/3)x^{3/2}}, \qquad \text{as } x \to +
%          \infty   \label{eq:Ebd} \\
%   M(x) & \sim \pi^{-1/2} (1+|x|)^{-1/4}, \qquad \text{as } |x| \to
%          \infty.   \label{eq:Mbd}
           %   \end{align}
the former as $x \to + \infty$, the latter as $|x| \to \infty$.
In the Hermite case it follows from \cite[p. 403]{olve74} that
$|\epsilon_2(\xi,\kappa_N)| \leq N^{-1} (M/E)(\kappa_N^{2/3} \zeta)$.
From \eqref{eq:first} and boundedness of $r(\xi)$, it follows that
\begin{equation*}
  \tau_N|\varphi(x)|
  \leq C (M/E)(\kappa_N^{2/3} \zeta)
  \leq
  \begin{cases}
    C E^{-1}(\kappa_N^{2/3} \zeta)  \, \leq C e^{-2s}  & \qquad \qquad
    \quad s > 0 \\
    C M(\kappa_N^{2/3} \zeta)  \quad \leq C(1+|s|)^{-1/4} \quad  &
    -N^{2/3-\epsilon} < s \leq 0,
  \end{cases}
\end{equation*}
where the first bound follows from \eqref{eq:EMbd} and JM12, Lemma 2
and the second from \eqref{eq:rkappa} and 
\eqref{eq:EMbd}.

\subsection{Proof of Theorem \ref{theorem Damian} (linking functions of GUE and GOE)}\label{section gue to goe}
The main engine of this result is an identity stated in \cite{forrester2001inter}, which relates the eigenvalues of a GUE to the eigenvalues of two independent GOEs. In particular, we use it in the following lemma.

\begin{lemma}
\label{lem: Damian}
Let $M_{N}^{\mathbb{C}}$ be an $N\times N$ GUE, and let $f$ be a function of bounded variation with total variation $\mathrm{TV}(f)$. If $M_{N}^{\mathbb{R}},\tilde{M}_{N}^{\mathbb{R}}$ are two independent GOEs, then
\begin{equation}
    \label{Damian 1.3}
    f\left(M_{N}^{\mathbb{C}}\right)\overset{d}{=}\frac{1}{2}\left( f\left(M_{N}^{\mathbb{R}}\right)+ f\left(\tilde{M}_{N}^{\mathbb{R}}\right)\right)+X_N,
\end{equation}
where $|X_N|\leq\mathrm{TV}(f)$, and $\overset{d}{=}$ denotes equality in distribution.
\end{lemma}
\begin{proof}
Let $M_{N}^{\mathbb{R}},\tilde{M}_{N+1}^{\mathbb{R}}$ be independent $N\times N$ and $(N+1)\times (N+1)$ GOEs. Call the eigenvalues of $M_{N}^{\mathbb{R}}$ and $\tilde{M}_{N+1}^{\mathbb{R}}$ $\{\lambda_i\}_{i=1}^{N}$ and $\{\tilde{\lambda}_i\}_{i=1}^{N+1}$, respectively. Further, denote the combined set of eigenvalues $\{\lambda_i\}_{i=1}^{N} \cup \{\tilde{\lambda}_i\}_{i=1}^{N+1}$ as $\Lambda^{+}$, and enumerate its elements in decreasing order
\[
\Lambda^{+}=\{\lambda_{1}^{+}\geq...\geq\lambda_{2N+1}^{+}\}.
\]
Theorem 5.2 of \cite{forrester2001inter} implies that the even elements of this set are equal in distribution to the eigenvalues of an $N\times N$ GUE.

Thus, if $M_N^{\mathbb{C}}$ is an $N\times N$ GUE, we have
\begin{align*}
    f(M_N^{\mathbb{C}})&\overset{d}{=}\sum_{i=1}^{N}f(\lambda_{2i}^{+}) \\
    &=\frac{1}{2}\left(\sum_{j=1}^{2N+1}f(\lambda_{j}^{+})+\sum_{i=1}^{N}\left[f(\lambda_{2i}^{+})-f(\lambda_{2i-1}^{+})\right]-f(\lambda_{2N+1}^{+}) \right)\\
    &=\frac{1}{2}\left(f(W_N^{\mathbb{R}})+f(\tilde{W}_{N+1}^{\mathbb{R}})-f(\lambda_{2N+1}^{+}) +\sum_{i=1}^{N}\left[f(\lambda_{2i}^{+})-f(\lambda_{2i-1}^{+})\right] \right).
\end{align*}
Notice that, since $\lambda_{j}^{+}$ are ordered, we have
\[
\left\vert\sum_{i=1}^{N}\left[f(\lambda_{2i}^{+})-f(\lambda_{2i-1}^{+})\right]\right\vert\leq\mathrm{TV}(f).
\]
Further, let $\tilde{M}_{N}^{\mathbb{R}}$ be the principal submatrix of $\tilde{M}_{N+1}^{\mathbb{R}}$, which is thus independent and equal in distribution to $M_N^{\mathbb{R}}$. If we let $\tilde{\mu}_1,...,\tilde{\mu}_N$ be the eigenvalues of $\tilde{M}_{N}^{\mathbb{R}}$, then Cauchy's interlacing theorem yields
\[
\tilde{\lambda}_1\geq\tilde{\mu}_1\geq\tilde{\lambda}_2\geq...\geq\tilde{\lambda}_N\geq\tilde{\mu}_N\geq\tilde{\lambda}_{N+1},
\]
and so we have
\begin{align*}
\left\vert f(\tilde{M}_{N+1}^{\mathbb{R}})-f(\lambda_{2N+1}^{+})-f(\tilde{M}_N^{\mathbb{R}})\right\vert
&=\left\vert\sum_{i=1}^{N}f(\tilde{\lambda}_i)-\sum_{i=1}^{N}f(\tilde{\mu}_i)+(f(\tilde{\lambda}_{N+1})-f(\lambda_{2N+1}^{+}))\right\vert\\
&\leq \sum_{i=1}^{N}|f(\tilde{\lambda}_i)-f(\tilde{\mu}_i)|+|f(\tilde{\lambda}_{N+1})-f(\lambda_{2N+1}^{+})|\\
&\leq \mathrm{TV}(f).
\end{align*}
We conclude that \eqref{Damian 1.3} holds.
\end{proof}

An immediate useful corollary is as follows.

\begin{corollary}
\label{cor:mean-var}
Under assumptions of Lemma \ref{lem: Damian},
\begin{align*}
\E f(M_{N}^{\mathbb{R}})
&=\E f(M_{N}^{\mathbb{C}})+O(\mathrm{TV}(f)),\\
\Var f(M_{N}^{\mathbb{R}})&\leq 2\Var f(M_{N}^{\mathbb{C}})+2\mathrm{TV}^{2}(f).
\end{align*}
\end{corollary}

\begin{remark}
  \label{rem:mean-var}
  Notice that corollary \ref{cor:mean-var} also holds for scaled Gaussian matrices \(W_N^{\mathbb{R}/\mathbb{C}}=M_N^{\mathbb{R}/\mathbb{C}}/\sqrt{N}\), since \(f(W_N^{\mathbb{R}/\mathbb{C}}) = g(M_N^{\mathbb{R}/\mathbb{C}})\) for \(g(\lambda) = f(\lambda/\sqrt{N})\), which satisfy \(\TV(f) = \TV(g)\).
\end{remark}

However, to finish proving Theorem \ref{theorem Damian} in its generality, we require the following technical lemma about tightness.

\begin{lemma}
\label{lemma Damian technical}
Let $X_N,Y_N$ be iid sequences of random variables such that $X_N+Y_N$ is tight. Then $X_N$ (and thus also $Y_N$) is tight. 
\end{lemma}
\begin{proof}
For any constant $K$, we have
\[
\Pr(X_N>K)=\Pr(X_N>K,Y_N>K)^{1/2}\leq\Pr(|X_N+Y_N|>K)^{1/2},
\]
and similarly,
\[
\Pr(X_N<-K)\leq\Pr(|X_N+Y_N|>K)^{1/2},
\]
which yield
\[
\sup_{N}\Pr(|X_N|>K)\leq 2\sup_{N}\Pr(|X_N+Y_N|>K)^{1/2}.
\]
The right hand side of the latter inequality can be made arbitrarily small, by the tightness of $X_N+Y_N$.
\end{proof}

With all these results in hand, we are ready to complete the proof of Theorem~\ref{theorem Damian}. We have
\begin{align}
\label{Damian 1.20}
    \left\vert \frac{f_N(W_N^{\mathbb{R}})-a_N}{b_N+\mathrm{TV}(f_N)}+\frac{f_N(\tilde{W}_N^{\mathbb{R}})-a_N}{b_N+\mathrm{TV}(f_N)} \right\vert
    &=2\left\vert\frac{(f_N(W_N^{\mathbb{R}})+f_N(\tilde{W}_N^{\mathbb{R}}))/2-a_N}{b_N+\mathrm{TV}(f_N)} \right\vert\\
    &\leq 2\left\vert\frac{(f_N(W_N^{\mathbb{R}})+f_N(\tilde{W}_N^{\mathbb{R}}))/2+X_N-a_N}{b_N} \right\vert+2\left\vert\frac{X_N}{\mathrm{TV}(f_N)}\right\vert \notag \\
    &\overset{d}{=}2\left\vert\frac{f_N(W_N^{\mathbb{C}})-a_N}{b_N}\right\vert+2\left\vert\frac{X_N}{\mathrm{TV}(f_N)}\right\vert.\notag
\end{align}
The first term in the latter sum is tight by assumption, whereas the second term is no larger than $2$. But since the two terms on the left hand side of  \eqref{Damian 1.20} are iid, Lemma \ref{lemma Damian technical} yields that they must be tight, and so
\[
f_N(W_N^{\mathbb{R}})=a_N+O_{\Pr}(b_N+\mathrm{TV}(f_N)).
\]

%\texttt{E.g. Gaussian paper Secs A.6, A.7, B.3.}

%\newpage

\section{Proofs for Wigner extension and Section 6}
\label{sec: appendix}

%\subsection{Proof of Proposition~\ref{prop:joint-convergence}} 
 \subsection{Proof of Proposition~\ref{prop:joint-convergence} %( $\mathbf{\xi_N}\overset{d}{\rightarrow}\mathbf{\xi}\implies\mathbf{\xi'_N}\overset{d}{\rightarrow}\mathbf{\xi}$)
 (about convergence of $\xi_N$ implying convergence of $\xi_N'$)}
\label{sec:proof-proposition-21}
 Let us first show that $ \xi_{Nj}(W_N') \overset{d}{\rightarrow} \xi_j $ marginally for each $j$. Fix some $\epsilon > 0$. Then, for large enough $N$,
    \[
        \mP(\xi_{Nj}(W_N') \leq s) \leq \mP(\xi_{Nj}(W_N') \leq s + \varepsilon - \eta_N) \leq \E Q_j(W_N', s + \varepsilon) + O(N^{-A}),
    \]
    Similarly, for $N$ large
    \[
        \E Q_j(W_N, s + \varepsilon) \leq \mP(\xi_{Nj}(W_N) \leq s + 2\varepsilon) + O(N^{-A}) \leq \mP(\xi_j \leq s + 2\varepsilon) + o_{s, \varepsilon}(1),
    \]
    where the last inequality follows from the convergence $ \xi_{Nj}(W_N) \overset{d}\rightarrow \xi_j $. 
    %and $o_{s, \varepsilon}(1) \rightarrow 0$ as $ N \rightarrow \infty $ and $s, \varepsilon$ are fixed.
    Since $Q_j(\cdot, s + \varepsilon)$ satisfies condition $F(\delta_{j, N})$, we have by Proposition~\ref{prop:multi-matching},
    \[
        \E Q_j(W_N', s + \varepsilon) \leq \E Q_j(W_N, s + \varepsilon) + o_{s,\varepsilon}(1).
    \]
    %where $o_{\varepsilon}(1)$ tends to zero as $\varepsilon$ is fixed and $N$ goes to infinity. 
    We therefore obtain for $N$ large,
    \[
        \mP(\xi_{Nj}(W_N') \leq s) \leq \mP(\xi_j \leq s + 2 \varepsilon) + o_{s, \varepsilon}(N) \, .
    \]
    %Taking lim-sup of both sides, we obtain the bound that holds for all $ \varepsilon > 0 $,
   % \[
  %      \lim \sup_N \mP(\xi_{Nj}(W_N') \leq s) \leq \mP(\xi_j \leq s + 2 \varepsilon) \, .
    %\]
    %Since the right-hand side does not depend on $\varepsilon$ and $\xi_j$ has continuous distribution, we obtain that
    %\[
   %     \lim \sup_{N} \mP(\xi_{Nj}(W_N') \leq s) \leq \mP(\xi_j \leq s) \, .
  %  \]
    Similarly, we can 
    %get a lower bound for the lim-inf. It follows that $ \xi_{Nj}(W_N') \overset{d}{\rightarrow} \xi_{j} $.
    obtain a lower bound
    \[
     \mP(\xi_{Nj}(W_N') \leq s) \geq \mP(\xi_j \leq s - 2 \varepsilon) + o_{s, \varepsilon}(N) \, .
    \]
    Since $\varepsilon$ can be chosen arbitrarily small and $\xi_j$ has continuous distribution, it follows that $ \xi_{Nj}(W_N') \overset{d}{\rightarrow} \xi_{j} $.
    
  Now let $\delta_N = \max_j \delta_{j,N}$.
  It suffices to show that for each $\mathbf{s} = (s_j)$ 
  \begin{align}
    \label{eq:sandwich2}
    \mP( \bxi_N \leq \mathbf{s}-\eta_N)  &\leq \mP( \bxi_N' \leq \mathbf{s} + \eta_N) + O(\delta_N)
    \;\; \text{ and}\\
    \mP( \bxi_N' \leq \mathbf{s} - \eta_N)
    &\leq \mP( \bxi_N \leq \mathbf{s}+\eta_N) + O(\delta_N).
    \label{eq:sandwich2a}
  \end{align}
Indeed, we then have
\begin{equation*}
  |\mP( \bxi_N' \leq \mathbf{s}) - \mP( \bxi_N \leq \mathbf{s}) |
  \leq 
  %\sum_j \mP ( |\xi_{nj}-s_j| \leq \eta_N )
  {\sum_j \mP ( |\xi_{Nj}-s_j| \leq \eta_N )} + \sum_j \mP ( |\xi_{Nj}'-s_j| \leq \eta_N )
  + O(\delta_N)
  \to 0
\end{equation*}
because each $\xi_{Nj}, \xi_{Nj}'$ has a continuous limiting distribution
function.

We verify inequality \eqref{eq:sandwich2}.
For each $A > 0$ large, we have from \eqref{eq:Qj} for $W_N$, then
Proposition~\ref{prop:multi-matching}
and then \eqref{eq:Qj} again, now for $W_N'$,
that
\begin{align*}
  \mP( \bxi_N \leq \mathbf{s} - \eta_N)
  & \leq \E \prod_j Q_j(W_N, s_j) + O(N^{-A}) \\
  & \leq \E \prod_j Q_j(W_N', s_j) + O(\delta_N) 
  \leq \mP( \bxi_N' \leq \mathbf{s}+\eta_N) +O(\delta_N).
\end{align*}
Inequality \eqref{eq:sandwich2a} follows similarly.

\subsection{Proof of Lemma~\ref{lem:neg-mom-bd} (bounds on inverse power sums for Wigner matrices)}
\label{sec:proof-lemma-28}

%\begin{proof}
  Let $\delta = N^{-2/3-\epsilon}$ and 
  $A_N = \{ \min_j |E - \lambda_j| > \delta \}$:
  by Proposition \ref{prop:wign-non-conc-2} this event has probability
  at least $1 - N^{-\epsilon/2}$.
  We will work on event $A_N$, and show that there the claims hold w.o.p.
  On $A_N$ the interval $I_0 = [E-\delta,E+\delta]$ contains no
  eigenvalues.
  Consider the `coronae' defined by $I_k = \{ x \in \mR ~:~ 2^{k-1} \delta <
    |x-E| \leq 2^k \delta \}$
    for $1 \leq k \leq k' = \min \{ k: E-2^{k}\delta \leq 1 \}$, and add
  two half-infinite intervals $I_{-1}$ and $I_{k'+1}$ to obtain a
  disjoint cover of $\mR$. We may then bound (on event $A_N$)
  \begin{equation}  \label{eq:SrE}
    S_r(E)
     \leq  \sum_{k=1}^{k'} \frac{\mathcal{N}_{W_N}(I_k)}{(2^{k-1} \delta)^r}
    + \frac{N}{(2^{k'}\delta)^r}.
  \end{equation}
The semicircle density is bounded by $\sqrt{2-x} \mathbf{1}_{x \leq
  2}$ and so $\rho_{\rm sc}([2-a,2-b]) \leq a^{3/2}$.
The lower endpoint of $I_k$ is
$E - 2^k\delta \geq 2 - 2^k\delta- \check{\sigma}_N N^{-2/3}$.
Since $\check{\sigma}_N^{3/2} \leq N^{\epsilon}$ for large $N$,
\[
\rho_{\rm sc}(I_k) \leq \sqrt{2}\left( (2^{k} \delta)^{3/2} + 
N^{\epsilon-1}\right).
%\check{\sigma}_N^{3/2}N^{-1}\right).
\]
  Proposition \ref{prop:P14analog} (ii) says that, with overwhelming
probability, simultaneously for all $k \leq k' = O(\log N)$, we have
% $ \mathcal{N}_{W_N}(I_k) \leq N \rho_{\rm sc}(I_k) + O(N^\epsilon)$. 
% Now $2^{-k'} \asymp \delta$, so we
% obtain w.o.p.
% \begin{equation*}
%   S_r(E)
%      \leq   2^r N \sum_{k=1}^{k'} \frac{\rho_{\rm sc}(I_k)}{(2^k \delta)^r}
%      + C N^\epsilon\fix{\delta^{-r}} + CN.
% \end{equation*}
% Using $\check{\sigma}_N \leq N^{2 \epsilon/3}$ for all sufficiently
% large $N$, we obtain
\begin{equation*}
  \mathcal{N}_{W_N}(I_k)
    \leq N \rho_{\rm sc}(I_k) + O(N^\epsilon)
    \leq \sqrt{2} N (2^k \delta)^{3/2} + C N^\epsilon.
\end{equation*}
Putting this into \eqref{eq:SrE} and noting that $2^{k'} \delta \in
[\frac{1}{2},3]$ we obtain w.o.p.
\begin{equation*}
  S_r(E)
  \leq 2^{r+1/2} N \sum_1^{k'} (2^k \delta)^{3/2-r} + C
  N^\epsilon\delta^{-r} + 2^r N.
\end{equation*}
The sum may be bounded using
\begin{equation*}
  N \delta^{3/2-r} \sum_1^{k'} 2^{(3/2-r)k}
   \leq
   \begin{cases}
     3^{1/2} N  & \text{if } r = 1 \\
     4 N^{-3\epsilon/2} \delta^{-r}   & \text{if } r \geq 2.
   \end{cases}
\end{equation*}
Observe that $N^\epsilon \delta^{-r} = N^{(2/3+\epsilon)r + \epsilon}$.
For $r=1$, this is $o(N^{-1})$ and so
$S_1(E) \leq C_1 N$ on $A_N$ w.o.p.
For $r \geq 2$ this is the dominant term, so that
$S_r(E) \leq C_r N^{2r/3+(r+1)\epsilon}$.

The bounds also hold for $S_r(E')$ uniformly in $|E'-E| \leq
\delta/2$, by increasing $C_r$ to $2^r C_r$. Indeed, for
such $E'$, on event $A_N$ we have $|\lambda_j - E'| \geq
\frac{1}{2}|\lambda_j - E|$ for all $j$. 
% and setting $s = \frac{3}{2} - r$, we find
% \begin{align*}
%   N \sum_{k=1}^{k'} \frac{\rho_{\rm sc}(I_k)}{(2^k \delta)^r}
%    \leq \sqrt{2}N \sum_{k=1}^{k'} (2^k \delta)^{\frac{3}{2}-r}
%                      + \sqrt{2}\frac{\check{\sigma}^{3/2}}{\delta^r}
%     \sum_{k=1}^{k'} 2^{-kr}
%   \leq \sqrt{2}N \delta^s\sum_{k=1}^{k'} 2^{ks} + \sqrt{2}N^\epsilon \delta^{-r}
% \end{align*}
% for all sufficiently large $N$. Since $2^{k'}\delta\leq 2$, we have 
% \[
% \sum_{k=1}^{k'}2^{ks}=\frac{2^{(k'+1)s}-2^s}{2^s-1}\leq\left\{
% \begin{array}{cc}
%      \delta^{-s}\frac{2^{2s}}{2^s-1}\leq 6\delta^{-s}&\text{if }s=1/2  \\
%      \frac{2^s}{1-2^s}\leq 3&\text{if }s\leq-1/2 .
% \end{array}
% \right.
% \]
% Thus, if $r=1$ (so that $s=1/2$) we find that w.o.p.
% \[
%   S_r(E)
%   \leq 12\sqrt{2}N+2\sqrt{2}N^{\varepsilon}\delta^{-1}+CN^{\varepsilon}\fix{\delta^{-1}}+CN
%   \leq C_1 N,
% \]
% while if $r\geq 2$ (so that $s \leq -1/2$) we find that w.o.p.
% \[
%   S_r(E)
%   \leq 3\sqrt{2}\cdot 2^r N\delta^{3/2-r}+ \sqrt{2}\cdot 2^r N^{\varepsilon}\delta^{-r}+CN^{\varepsilon}\fix{\delta^{-r}}+CN
%   \leq C_r N^{(3/2 + \zeta)r + \varepsilon}.
% \]
%\end{proof}

\subsection{Proof of \cref{thm: isotropic} (isotropic local law and delocalization)}\label{sec:B3} \ We will modify the
proofs of Theorems 2.2 and 2.5 of \cite{knowles2013isotropic} to reach two
main goals: replace their sub-exponential assumption on the entries of
Wigner matrices by the uniform moment bound \eqref{eq:W3}, and allow for an arbitrary variance
profile along the Wigner diagonal. We are going to reformulate
$\zeta$-high probability bounds used in \cite{knowles2013isotropic} in terms
of weaker polynomial bounds. For this, it will be convenient to use
the concept of stochastic domination, as defined in
\cite[def. 2.5]{benaych2018}: 
\smallskip

\textit{Definition of stochastic domination.}\ Let
\[
X=\left(X^{(N)}
(u):N\in \mathbb{N},u\in U^{(N)}\right),\qquad Y=\left(Y^{(N)}
(u):N\in \mathbb{N},u\in U^{(N)}\right) 
\]
be two families of non-negative random variables, where $U^{(N)}$ is a possibly $N$-dependent parameter set. We say that $X$ is stochastically dominated by $Y$, uniformly in $u$, and write $X\prec Y$, if for all (small) $\epsilon>0$ and (large) $A>0$ we have
\[
\sup_{u\in U^{(N)}}\Pr\left[X^{(N)}(u)>N^{\epsilon}Y^{(N)}(u)\right]\leq N^{-A}
\]
for large enough $N\geq N_0(\epsilon,A)$. If for some complex family $X$ we have $|X|\prec Y$, we also write $X=O_\prec(Y)$.
\smallskip

With this definition, we are ready to point out necessary changes to the proofs of Theorems 2.2 and 2.5 of \cite{knowles2013isotropic}.  We use mostly the same notation, and
refer the reader to \cite{knowles2013isotropic} for definitions. For example, the resolvent matrix will be denoted as $G(z)$ instead of $R(z)$ as in the main body of our paper.
We use numbering (KY3.xx) for formula (3.xx) in \cite{knowles2013isotropic},
and (3.xxa) for a formula here which is a stochastic dominance analog
of (KY3.xx).  

\textit{Section KY3.}\ To accommodate our setting, proposition KY3.1 should be reformulated as follows.
\begin{proposition}
\label{prop: moments}
  Fix any $\tau>0$, $0<\epsilon<\tau$, and $n>0$. Then under assumptions of \cref{thm: isotropic}, for all deterministic normalized $\mathbf{v,w}\in \mathbb{C}^N$ and all $z\in\mathbf{S}(\tau)$,
  \begin{equation}
      \label{eq:prop3.1 ky}
      \mathbf{E}|G_{\mathbf{vw}}(z)-s_{sc}(z)\mathbf{v}^\ast\mathbf{w}|^n\leq(N^{\epsilon}\Psi(z))^n
  \end{equation}
  for all sufficiently large $N\geq N_0(\tau,\epsilon,n)$.
\end{proposition}
The above proposition implies \cref{isotropic law}. Indeed, let $A$ be large and take $n=A/\epsilon$. Then by Markov's inequality 
\[
\Pr(|G_{\mathbf{vw}}(z)-s_{sc}(z)\mathbf{v}^\ast\mathbf{w}|\geq N^{2\epsilon} \Psi(z))\leq \frac{\mathbf{E}|G_{\mathbf{vw}}(z)-s_{sc}(z)\mathbf{v}^\ast\mathbf{w}|^n}{N^{\epsilon n}(N^{\epsilon} \Psi(z))^n}\leq N^{-\epsilon n}= N^{-A}
\]
for sufficiently large $N$, which in our notations means that 
\[
\mathbf{v}^\ast R(z)\mathbf{w}=s_{sc}(z)\mathbf{v}^\ast\mathbf{w}+O(N^\epsilon \Psi(z))
\]
w.o.p., as required.
\smallskip

\textit{Subsection KY3.1.} Replace Lemma KY3.5 by the following one.
\begin{lemma}
\label{lem:replacement of KY lem3.5}
Let $a_1,...,a_N$ be independent random variables with zero mean and all moments bounded uniformly in $N$. Then for any deterministic complex numbers $A_i$, we have
\[
\left\vert \sum_{i=1}^NA_ia_i\right\vert\prec\left(\sum_{i=1}^N|A_i|^2\right)^{1/2}.
\]
\end{lemma}
A proof of this lemma is almost identical to the proof of Lemma 8.2 in \cite{ErdosYY2012b}, so we omit it. Further, replace Theorems KY3.6 and KY3.7 by Theorems 2.6 and 2.9 (respectively) from \cite{benaych2018}.

\textit{Subsection KY3.2}. We need to reformulate the statements formulated in terms of $\zeta$-high probability using the notion of stochastic domination. In particular, lemma KY3.8 should be reformulated as follows.
\begin{lemma}
\label{lem: 3.8 of ky}
Fix $\tau>0$. Then,
\[
|\mathcal{G}_{\mathbf{v}i}(z)|+|\mathcal{G}_{i\mathbf{v}}(z)|+|G_{\mathbf{v}i}(z)|+|G_{i\mathbf{v}}(z)|\prec \sqrt{\frac{\Im G_{\mathbf{vv}}(z)}{N\eta}}+|v_i|
\]
for all $z\in\mathbf{S}(\tau)$.

\end{lemma}
In \cite{knowles2013isotropic}'s proof of the lemma, replace the first display by
\[
|\mathcal{G}_{\mathbf{v}i}(z)|\prec \left(\frac{1}{N}\sum_k^{(i)}|G_{\mathbf{v}k}^{(i)}|^2\right)^{1/2},
\]
which holds by \cref{lem:replacement of KY lem3.5}. Further, change inequality (KY3.19) to $|G_{ii}|\prec 1$ (which follows from Theorem 2.6 of \cite{benaych2018}, and then, in all of the remaining displayed inequalities change $\leq$ to $\prec$. 
Note that so far, we have not used any information about the diagonal variance profile. This information will be used in the next subsection.

\textit{Subsections KY3.3-KY3.4} \ To manage the modification of the proof to
cover stochastic dominance and weaker conditions on diagonal moments, it is
convenient to somewhat reorganize the material in Sections KY3.3 and 3.4,
along with ideas from Case 1 of KY4.1.  At cost of some duplication of
text from KY, we thus write out this part of the proof in relatively
self-contained form. 

Given a resolvent matrix $G(z)$ as in KY Th 2.2 we consider three
cases of linear functionals $L_\bv G(z)$ and corresponding control
functions $\Pi(z)$:
\begin{equation*}
  L_\bv G =
  \begin{cases}
    \Re G_{\bv \bv} - \Re s_{sc} \\
    \Im G_{\bv \bv} - \Im s_{sc} \\
    \Im G_{\bv \bv}
  \end{cases}
  \qquad
  \Pi =
  \begin{cases}
    \Psi \\
    \Psi \\
    \Phi.
  \end{cases}
\end{equation*}
Fix $\tau > 0, 0 < \epsilon < \tau$ and even $n \geq 2$.
We seek to prove inequalities\footnote{For $L_\bv G = \Im G_{\bv
    \bv}$, Lemma KY3.9 uses $\Phi$, but the bound with $\Psi$ is
  better and allows a more uniform treatment. }
\begin{equation*}
  \tag{3.20a}
  \E (L_\bv G(z) )^n \leq ( N^{\epsilon/2} \Pi(z))^n,
%  \E (L_\bv G(z) )^n \leq (\Pi(z))^n,
\end{equation*}
for all $z \in \mathbf{S}(\tau)$.
The first two choices for $L_\bv G$ together yield
% the isotropic law Th. 2.2,
Proposition \ref{prop: moments},
while the third is needed for an intermediate step in the
proof.

We verify that (3.20a) holds when $H_0$ (matrix $W_N$ in the notations of previous sections) is a GOE/GUE matrix.
In that case $\E (L_\bv G(z))^n = \E (L_{\mathbf{e}_1} G(z))^n $ by
unitary invariance.
From the entrywise local law and (KY3.4), we have for $z \in
\mathbf{S}(\tau)$ that
\begin{equation*}
  \Im G_{11}(z) \prec \Phi(z), \qquad
  |G_{11}(z) - s_{sc}(z)| \prec \Psi(z),
\end{equation*}
so that $L_\bv G(z) \prec \Pi(z)$.
Since $\Psi(z) \gtrsim N^{-1/2}$ and $ \E (L_\bv G(z))^n  \leq N^{p(n)}$
from the rough bound $|G_{11}(z)| \leq \eta^{-1} \leq N$,
\cite[][Lemma 7.1]{benaych2018} implies that
$\E (L_\bv G(z) )^n \prec \Pi(z)^n$, which yields (3.20a) for GOE/GUE.

``From now on we work on the product space generated by the Wigner
matrix $H = (N^{-1/2} W_{ij})_{i,j}$ and the GOE/GUE matrix $(N^{-1/2}
  V_{ij})_{i,j}$. We fix a bijective
ordering map on the index set of the independent matrix elements,
\begin{equation*}
\tag{KY3.21}
  \varphi : \{(i,j) : 1 \leq i \leq j \leq N \} \to \{1, \ldots,
  \gamma_{\rm max} \} \quad
  \text{where } \quad \gamma_{\rm max} := \frac{N(N+1)}{2},
\end{equation*}
and denote by $H_\gamma = (h_{ij}^\gamma), \gamma = 0, \ldots,
\gamma_{\rm max}$, the Wigner matrix with upper-triangular
entries  defined by
\begin{equation*}
  h_{ij}^\gamma =
  \begin{cases}
    N^{-1/2} W_{ij} & \text{if } \varphi(i,j) \leq \gamma, \\
    N^{-1/2} V_{ij} & \text{otherwise}.
  \end{cases}
\end{equation*}
In particular, $H_0$ is a GOE/GUE matrix and $H_{\gamma_{\rm max}} =
H$.

Let $E^{(ij)}$ denote the matrix whose matrix elements are given by
$E^{(ij)}_{kl} = \delta_{ik} \delta_{jl}$.
Fix $\gamma \geq 1$ and let $(a,b)$ be determined by $\varphi(a,b) =
\gamma$. We shall compare $H_{\gamma - 1}$
with $H_\gamma$ for each $\gamma$ and then sum up the
differences. Note that the matrices $H_{\gamma - 1}$ and $H_\gamma$
differ only in the entries $(a,b)$ and $(b,a)$ and they can be written
as
\begin{equation*}
  \tag{KY3.22}
  H_{\gamma-1} = Q + N^{-1/2}V \quad \text{where} \quad
  V := V_{ab} E^{(ab)} + \mathbf{1}(a \neq b) V_{ba} E^{(ba)},
\end{equation*}
and
\begin{equation*}
  H_{\gamma} = Q + N^{-1/2}W \quad \text{where} \quad
  W := W_{ab} E^{(ab)} + \mathbf{1}(a \neq b) W_{ba} E^{(ba)},
\end{equation*}
here the matrix $Q$ satisfies $Q_{ab} = Q_{ba} = 0.$

Next, we introduce the Green functions
\begin{equation*}
  \tag{KY3.23}
  R:= \frac{1}{Q-z}, \quad 
  S:= \frac{1}{H_{\gamma-1}-z}, \quad 
  T:= \frac{1}{H_\gamma-z}, \quad 
\end{equation*}
which are well-defined for $\eta > 0$ since $Q$ and $H_\gamma$ are
self-adjoint. Using the notation $G^\gamma = (H_\gamma - z)^{-1}$,
we have the telescopic sum''
\begin{equation*}
  \tag{3.24a}
  \E (L_\bv G^{\gamma_{\rm max}} )^n  - \E (L_\bv G^0 )^n
     = \sum_{\gamma =1}^{\gamma_{\rm max}} (X_\gamma - X_{\gamma -1}),
\end{equation*}
where, since $n$ is even, $X_\gamma = \E (L_{\bv} G^\gamma)^n \geq 0$.
% \begin{equation*}
%   X_\gamma = \E (L_{bv} G^\gamma)^n \geq 0.
% \end{equation*}
Note that in the $R,S, T$ notation, $X_\gamma = \E (L_{\bv}T)^n$ and
$X_{\gamma -1} = \E (L_{\bv}S)^n$.

% We start by expanding the difference
% \begin{equation*}
%   (L_{\bv}S)^n - (L_{\bv} R)^n
%     = \sum_{m=1}^n \binom{n}{m} (L_{\bv}S - L_{\bv}R)^m (L_{\bv} R)^{n-m}.
% \end{equation*}
For any $K \in \mathbb{N}$ we have the resolvent expansions
\begin{align*}
  \tag{KY3.25}
  S & = R + \sum_{k=1}^{K-1} N^{-k/2} (-RV)^kR + N^{-K/2}(-RV)^K S, \\
  \tag{KY3.26}
  R & = S + \sum_{k=1}^{K-1} N^{-k/2} (SV)^kS + N^{-K/2}(SV)^K R.    
\end{align*}

With $K = 4$ in (KY3.26),  using the entrywise local law for the Wigner matrix
$S$, and the rough bound $\| R \| \leq \eta^{-1}$ to estimate the
remainder term in (KY3.26), and recalling \eqref{eq:W3} instead of (KY2.1), we find
\begin{equation*}
  \tag{3.27a}
  |R_{ij} - \delta_{ij} s_{sc}|
  \prec |S_{ij} - \delta_{ij} s_{sc}| + N^{-1/2}
  \prec \Psi.
\end{equation*}
There are trivial changes in (KY3.28) - (KY3.30), for later use we
record
\begin{align*}
  |R_{\bv a}|
  & \prec \sqrt{\frac{\Im S_{\bv \bv}}{N \eta}} + \Psi + |v_a|
    \tag{3.29a} \\
  |S_{\bv \bv} - R_{\bv \bv}|
  & \prec N^{-1/2} \bigg( \frac{\Im S_{\bv \bv}}{N \eta} + |v_a|^2 +
    |v_b|^2 \bigg). \tag{3.30a}
\end{align*}

We now apply (KY3.25) with $K = 4$ and introduce the notation
$S - R = \sum_{k=1}^4 Y_k$, whereby $Y_k$ has $k$ factors $V$.
Then $L_{\bv}S - L_{\bv} R = \sum_{k=1}^4 \dot{L}_{\bv} Y_k$,
where
$\dot L_{\bv} Y$ is either $\Re Y_{\bv \bv}$ or $\Im Y_{\bv \bv}$,
since the terms involving $s_{sc}(z)$, if present, cancel.
We expand the difference
\begin{align}
  (L_{\bv}S)^n - (L_{\bv} R)^n
     = \sum_{m=1}^n \binom{n}{m} (L_{\bv}S - L_{\bv}R)^m (L_{\bv}
      R)^{n-m} 
    = \sum_{m=1}^n  (L_{\bv} R)^{n-m} \sum_{k=m}^{4m}
      A_{m,k}  \label{eq:Lvdiff}
\end{align}
% \begin{equation*}
%   \binom{n}{m} (L_{\bv}S - L_{\bv}R)^m
%     = \sum_{k=m}^{4m} A_{m,k},
% \end{equation*}
where
\begin{equation*}
  A_{m,k} = \binom{n}{m} \sum_{k_1, \ldots, k_m = 1}^4 \mathbf{1}(k_1
  + \cdots k_m = k) \prod_{i=1}^m \dot L_{\bv} Y_{k_i}.
\end{equation*}
Thus $A_{m,k}$ collects all terms with $k$ factors $V_{ab}$ or
$V_{ba}$, and hence is of order $N^{-k/2}$.
% in the expansion of $(L_{\bv}S - L_{\bv}R)^m$.
% Thus we have
% \begin{equation*}
%   \tag{3.32}
%   (L_{\bv}S)^n - (L_{\bv} R)^n
%   = \sum_{m=1}^n  (L_{\bv} R)^{n-m} \sum_{k=m}^{4m} A_{m,k}
%   = \mathcal{A}_\gamma + \mathcal{A}_\gamma^\prime.
% \end{equation*}
Break the sum in \eqref{eq:Lvdiff} in two so that $k \leq k_\gamma -1$
and $k \geq k_\gamma$.
The value $k_\gamma$ is chosen so that the first $k_\gamma - 1$
moments of $V_{ab}$ and $W_{ab}$ are the same, meaning that for any $t_1,t_2\in\mathbb{N}$ s.t. $t_1+t_2<k_\gamma$, we have $\E\left(V_{ab}^{t_1}\overline{V_{ab}}^{t_2}\right)=\E\left(W_{ab}^{t_1}\overline{W_{ab}}^{t_2}\right)$.
Thus we take $k_\gamma = 4$ for $a \neq b$ and $k_\gamma = 2$ when $a
= b$.
We obtain
\begin{equation*}   \tag{KY3.32}
  (L_{\bv}S)^n - (L_{\bv} R)^n
  = \sum_{m=1}^n \mathcal{A}_{\gamma m}
  + \sum_{m=1}^n \mathcal{A}_{\gamma m}'
  = \mathcal{A}_\gamma + \mathcal{A}_\gamma',
\end{equation*}
where, for example,
\begin{equation*}
  \mathcal{A}_{\gamma m}'
    =  (L_{\bv} R)^{n-m} \sum_{k= k_\gamma \vee m}^{4m}  A_{m,k}.
\end{equation*}
% The terms $\mathcal{A}_\gamma$ and  $\mathcal{A}_\gamma^\prime$ are
% defined by breaking the sum on $k$ so that $k \leq k_\gamma -1$ in
% $\mathcal{A}_\gamma$ and $k \geq k_\gamma$ in
% $\mathcal{A}_\gamma^\prime$.
Thus $\E \mathcal{A}_\gamma$ depends on the randomness only through
$Q$ and the first $k_\gamma - 1$ moments of $V_{ab}$.
% The value $k_\gamma$ is chosen so that the first $k_\gamma - 1$
% moments of $V_{ab}$ and $W_{ab}$ are the same.
% Thus we take $k_\gamma = 4$ for $a \neq b$ and $k_\gamma = 2$ when $a
% = b$.
Consequently $\E \mathcal{A}_\gamma$ equals the corresponding
term in the expansion (KY3.32) of $\E (L_{\bv}T)^n - \E (L_{\bv} R)^n$.

For the higher order terms, the analog of the key inequality proved by
KY has the form 
\begin{equation*}
  \tag{3.33a}
  | \E \mathcal{A}_\gamma^\prime |
    \leq \frac{\mathcal{E}_{ab}}{N^{\epsilon/2}} \left[ \E (L_{\bv} S)^n +
      ( N^{\epsilon} \Pi)^n \right].
\end{equation*}
The factor $\mathcal{E}_{ab} = \mathcal{E}(v_a,v_b,N)$ will be
detailed below; for now we simply need that
$\varepsilon_\gamma := N^{-\epsilon/2} \mathcal{E}_{ab} \geq 0$ satisfies
$\sum_\gamma \varepsilon_\gamma \leq \frac{1}{2}$.
% \begin{equation*}
%   0 \leq \varepsilon_\gamma \leq \tfrac{1}{2}, \qquad
%   \sum_\gamma \varepsilon_\gamma \leq 1.
% \end{equation*}

Before proving (3.33a), we show how it implies (3.20a).
Repeating the derivation of (KY3.32) for $T$ instead of $S$, using that
the first $k_\gamma - 1$ moments of $V_{ab}$ and $W_{ab}$ are the
same, and using (3.33a) and its analog with $S$ replaced by $T$, we
find \footnote{KY omit the factor 2, but since there are $O(N^2)$
  inequalities, it should perhaps be tracked explicitly.} 
\begin{equation*}
  X_\gamma - X_{\gamma-1} \leq \varepsilon_\gamma(X_\gamma +
  X_{\gamma-1} + 2 \Pi_n),
\end{equation*}
for $1 \leq \gamma \leq \gamma_{\rm max} = N(N+1)/2$
and with $\Pi_n = (N^\epsilon \Pi)^n$.
% $\sum_\gamma \varepsilon_\gamma \leq \frac{1}{2}$. Then
% \begin{equation*}
%   X_{\gamma_{\rm max}} \leq 2 e^2 \Pi^n.
% \end{equation*}
% \end{lemma}
Rewriting this and making the
  % \begin{equation*}
  %   (1-\varepsilon_\gamma) X_\gamma
  %     \leq (1+\varepsilon_\gamma) X_{\gamma-1} + 2 \varepsilon_\gamma
  %     \Pi^n. 
  % \end{equation*}
abbreviation $r_\gamma =
(1-\varepsilon_\gamma)^{-1}(1+\varepsilon_\gamma) \geq 1$, we
therefore find that
\begin{equation*}
  X_\gamma \leq r_\gamma(X_{\gamma-1} + 2 \epsilon_\gamma \Pi_n).
\end{equation*}
Since (3.20a) holds for GOE/GUE, we have 
the initial estimate $X_0 \leq \Pi_n$, and find on iteration that
\begin{equation*}
  X_\gamma \leq \Big( \prod_{j=1}^\gamma  r_j \Big) \Big( 1 + 2
  \sum_{j=1}^\gamma \varepsilon_j \Big) \Pi_n
\end{equation*}
Using $\log(1 + x) \leq x$ and $\sum_\gamma \varepsilon_\gamma \leq
\frac{1}{2}$ we find that $\log( \prod r_\gamma) \leq 4 \sum
\varepsilon_\gamma \leq 2$ which implies that
\begin{equation*}
   \E (L_{\bv} G)^n = X_{\gamma_{\max}} \leq C \Pi_n,
\end{equation*}
which is (3.20a). [We may take $C = 2 e^2$.]

\bigskip
\bigskip

We turn to the proof of (3.33a). The key step is a high-probability
bound for $\prod_{i=1}^m \dot L_{\bv} Y_{k_i}$, see \eqref{eq:mainbd}
and \eqref{eq:prod-bound} below.
Consider first a term $ (Y_k)_{\bv \bv}$, which up to sign is given by
$N^{-k/2} [(RV)^k \check{S}]_{\bv \bv}$, where
$\check{S}$ is short for $R$ if $k \leq 3$ and for $S$ if $k=4$.
Multiplying out gives
\begin{equation*}
  [ (RV)^k \check{S}]_{\bv \bv}
    = \sum R_{\bv a_1} V_{a_1 b_1} R_{b_1 a_2} \cdots V_{a_k b_k}
    \check{S}_{b_k \bv},
\end{equation*}
where the sum is over all choices of $a_1,b_1, \ldots, a_k, b_k$ such
that $\{ a_\ell, b_\ell\} = \{a, b \}$ for each $\ell$.
The number of terms in the sum is bounded: $2^k \leq 2^{4n}$ if $a
\neq b$ and $1$ if $a = b$.
% When $a \neq b$, there are $2^k$ terms in the sum, as there are two
% choices for each $a_\ell$ (and then each $b_\ell$ is fixed).
Let
\begin{equation*}
  s_i = \# \text{matrix elements } R_{\bv a}, R_{a \bv}, S_{\bv a},
  S_{a \bv},
\end{equation*}
in that term when $k=k_i$,
and write $t_i$ for the corresponding number when $a$ is replaced with $b$.
We have $s_i, t_i \in \{0,1,2\}$ and $s_i + t_i = 2$.
% Let $v(k,s)$ denote the
% number of terms in the sum corresponding to a given value of
% $s$.\footnote{There is an explicit expression for $v(k,s) = 2^{k-2}
%   \binom{2}{s}$ if $k \geq 2$ and $= 2 \mathbf{1}\{s=1\}$ for $k=1$. We won't
%   need it, however.}

Now abbreviate some terms used by KY:
\begin{equation*}
  \mathcal{R}_{\bv,a} = |R_{\bv a}|+ |R_{a \bv}|+|S_{\bv a}|
  +|S_{a \bv}|,
\end{equation*}
ditto for $\mathcal{R}_{\bv,b}$.
Since $|V_{ab}| \prec 1$ from \eqref{eq:W3} and  $|R_{ij}| \prec 1$ from (3.27a) and
the entrywise local law, we obtain
\begin{equation*}
|(Y_{k_i})_{\bv \bv}| =
  N^{-k_i/2}| [(RV)^{k_i} \check{S}_i]_{\bv \bv}|
    \prec N^{-k_i/2} \sum_{s_i=0}^2 
    \mathcal{R}_{\bv,a}^{s_i} \mathcal{R}_{\bv,b}^{t_i},
  \end{equation*}
with the sum reducing to $\mathcal{R}_{\bv,a}^2$ if $a = b$. 
Since $|\dot L_\bv Y_{k_i}| \leq |(Y_{k_i})_{\bv \bv}|$, with
% if $\mathbf{k} = (k_1, \ldots, k_m)$ and
% $k_1 + \cdots + k_m = k$,
$ k = \sum_1^m k_i$ and similarly for $s$ and $t$,
we then get
\begin{equation}  \label{eq:mainbd}
  \Big| \prod_{i=1}^m \dot L_\bv Y_{k_i} \Big|
    \prec N^{-k/2} \sum_{s=0}^{2m} 
    \mathcal{R}_{\bv,a}^{s} \mathcal{R}_{\bv,b}^{t},
\end{equation}
for $a \neq b$, with the conditions
\begin{equation*}
  \tag{KY3.35}  s+t = 2m, \qquad k \geq \max \{s,t \},
\end{equation*}
while if $a=b$ in \eqref{eq:mainbd} the sum is replaced simply by
$\mathcal{R}_{\bv,a}^{2m}$. 

% Define an intermediate control term
% \begin{equation*}
%   \Omega_\bv = \Omega_\bv(z) = \frac{\Im S_{\bv \bv}}{N\eta} + \Psi.
% \end{equation*}
% As usual, the constant $C$ may change at each appearance,
% e.g. for $n \leq \phi^\zeta$, 
% $(\Omega_\bv(C_\zeta) + \Omega_\bv(C_\zeta'))^n \leq \Omega_\bv(C_\zeta'').$
% $\Omega = \Omega(\bv,z) = \phi^{C_\zeta} \frac{\Im S_{\bv \bv}}{N\eta}
% + \phi^{C_\zeta} \Psi^2 + CN^{-1/2}$.

% We now show that 
% \begin{equation}
%   \label{eq:Rs-bd}
%   \mathcal{R}_{\bv,a}^s \prec \Omega_\bv^{s/2}(1+N^{s/4} |v_a|^s).
% \end{equation}
% Indeed,
Using \cref{lem: 3.8 of ky} and (3.29a) we get the bound 
\begin{equation*}
  \mathcal{R}_{\bv,a} \prec \sqrt{\frac{\Im S_{\bv \bv}}{N\eta}}
  + \Psi + |v_a|
  =: \sqrt{x} + |v_a|.
\end{equation*}
Using KY Lemma 3.10,
% $(C|v_a| + \sqrt x)^s \leq (C|v_a|)^s + (s \sqrt x)^s$ and so
\begin{align*}
  (\sqrt x + |v_a|)^s
   \prec ( \sqrt x)^s + |v_a|^s
%  \leq (\phi^\zeta \sqrt x)^s + (C|v_a|)^s \\
   \leq (x + N^{-1/2})^{s/2} (1+N^{s/4}|v_a|^s).
\end{align*}
On $\mathbf{S}(\tau)$, we have $\Psi \lesssim N^{-\tau/2}$ and 
(KY3.5) implies that 
$x +  N^{-1/2} \lesssim \Omega_\bv$, with the intermediate control term
\begin{equation*}
  \Omega_\bv = \Omega_\bv(z) = \frac{\Im S_{\bv \bv}}{N\eta} + \Psi.
\end{equation*}
% \begin{equation*}
%   x +  N^{-1/2} \lesssim \frac{\Im S_{\bv \bv}}{N\eta} + \Psi.
% \end{equation*}
We arrive at the bound
\begin{equation}
  \label{eq:Rs-bd}
  \mathcal{R}_{\bv,a}^s \prec \Omega_\bv^{s/2}(1+N^{s/4} |v_a|^s).
\end{equation}

%\eqref{eq:Rs-bd}.
% Since $\phi^{2\zeta}x \leq 2 \phi^{2\zeta}(\phi^{2 C_\zeta} \frac{\Im
%   S_{\bv \bv}}{N\eta} + \phi^{2C_\zeta} \Psi^2)$, bound
% \eqref{eq:Rs-bd} follows after increasing $C_\zeta$.

Now recall that $k \geq k_\gamma$ in the higher order term
$\mathcal{A}_\gamma^\prime$. In the off-diagonal case, this means that
both $k \geq \max\{4,m\}$ and $k \geq \max\{s, t\} \geq (s+t)/2$
from (KY3.35) and so
\begin{equation} \label{eq:R-prodbd}
  N^{-k/2} \mathcal{R}_{\bv,a}^{s} \mathcal{R}_{\bv,b}^{t}
  \leq [N^{-\max\{1,s/4\}} \mathcal{R}_{\bv,a}^{s}]
       [N^{-\max\{1,t/4\}} \mathcal{R}_{\bv,a}^{t}].
\end{equation}
In the diagonal case, $k \geq \max\{2,m\}$, and we have just
$  N^{-k/2} \mathcal{R}_{\bv,a}^{2m}
    \leq N^{-\max\{1,m/2\}} \mathcal{R}_{\bv,a}^{2m}$.
% \begin{equation*}
%   N^{-k/2} \mathcal{R}_{\bv,a}^{2m}
%     \leq N^{-\max\{1,m/2\}} \mathcal{R}_{\bv,a}^{2m}.
% \end{equation*}

Checking cases, one sees that for $s \geq 0$,
\begin{equation*}
  N^{-\max\{1,s/4\}} (1+N^{s/4}|v_a|^s)
    \leq 2 N^{-1} + N^{-1/2}|v_a| + |v_a|^2 =: S_a.
\end{equation*}
From \eqref{eq:R-prodbd}, \eqref{eq:Rs-bd}, and recalling that $s+t =
2m$, we obtain 
\begin{align*}
  N^{-k/2}  \mathcal{R}_{\bv,a}^{s} \mathcal{R}_{\bv,b}^{t}
   \prec \Omega_\bv^m S_a S_b  \qquad \text{for } a \neq b 
\end{align*}
and $  N^{-k/2} \mathcal{R}_{\bv,a}^{2m}
 \prec \Omega^m_\bv S_a$ for  $a = b$.
% \begin{align*}
%   N^{-k/2}  \mathcal{R}_{\bv,a}^{s} \mathcal{R}_{\bv,b}^{t}
%   & \leq \Omega^m S_a S_b  \qquad \text{for } a \neq b \\
%   N^{-k/2} \mathcal{R}_{\bv,a}^{2m}
%   & \leq \Omega^m S_a  \qquad \quad \text{for } a = b.
% \end{align*}
 As these bounds are uniform in the relevant $s,t$,
% and recalling that  $\sum_{s=0}^{2m} v(\mathbf{k},s) = 2^k$, 
\eqref{eq:mainbd} implies that 
\begin{equation}   \label{eq:prod-bound}
  \Big| \prod_{i=1}^m \dot L_\bv Y_{k_i} \Big|
    \prec \mathcal{E}_{ab}  \Omega^m_\bv,
\end{equation}
with $\mathcal{E}_{ab} = S_aS_b$ for $a \neq b$ and $=S_a$ if $a = b$.

% We now argue that for each $m$,
% \begin{equation}
%   \label{eq:m-bd}
%    \mathcal{A}_{\gamma m}'
%     =  |L_\bv R|^{n-m} \sum_{k=\max\{k_\gamma,m\}}^{4m}
%       |A_{m,k}| 
%       \prec \mathcal{E}_{ab} N^{-m \epsilon/4} \left[ (L_\bv S)^n
%          + \left( \frac{\Im S_{\bv \bv}}{N\eta}
%            \right)^n
%            +  \Psi^n \right].
%   % \E \mathcal{A}_{\gamma m}'
%   %   = \E \bigg[ |L_\bv R|^{n-m} \sum_{k=\max\{k_\gamma,m\}}^{4m}
%   %     |A_{m,k}| \bigg]
%   %     \leq \mathcal{E}_{ab} \phi^{-m} \E \left[ (L_\bv S)^n
%   %        + \left( \phi^{C_\zeta} \frac{\Im S_{\bv \bv}}{N\eta}
%   %          \right)^n
%   %          + (\phi^{C_\zeta} \Psi)^n \right].
% \end{equation}
The number of terms in $A_{m,k}$ is crudely bounded by
$\binom{n}{m} 4^m \leq C_n$ and so from \eqref{eq:prod-bound} 
\begin{align*}
  \sum_{k=k_\gamma \vee m}^{4m} |A_{m,k}|
    \prec \mathcal{E}_{ab} \Omega_{\bv}^m
    = \mathcal{E}_{ab} N^{-m \epsilon} ( N^{\epsilon} \Omega_\bv)^m.
     % \mathcal{E}_{ab} (\Omega_{\bv}^m(C_\zeta)
     % = \mathcal{E}_{ab} \phi^{-m} \Omega_{\bv}^m(C_\zeta),
\end{align*}
Using $x^{n-m}y^m \leq (x+y)^n$, we have
\begin{equation*}
  \mathcal{A}_{\gamma m}'
  \leq \mathcal{E}_{ab} N^{-m \epsilon} |L_\bv R|^{n-m}
    ( N^{\epsilon} \Omega_\bv)^m
     \\
    \leq \mathcal{E}_{ab} N^{-m \epsilon}[ |L_\bv R| +
    N^{\epsilon}\Omega_\bv]^n. 
\end{equation*}
From (3.30a) and (KY3.5),
\begin{equation*}
  |L_\bv R - L_\bv S|
  \leq |R_{\bv \bv} - S_{\bv \bv}|
  \prec  \frac{\Im S_{\bv \bv}}{N\eta} + N^{-1/2}
  \lesssim  \Omega_\bv,
\end{equation*}
and so, using KY Lemma 3.10,
\begin{equation*}
  [ |L_\bv R| + N^{\epsilon} \Omega_\bv]^n
  \prec [ |L_\bv S| + N^{\epsilon} \Omega_\bv ]^n
  \prec (L_\bv S)^n + (N^{\epsilon} \Omega_\bv)^n.
%  \leq  (L_\bv S)^n +  (\Omega_\bv(C_\zeta))^n.
\end{equation*}
% Using $(a+b)^n \leq (2a)^n +(2b)^n$
% %and further increasing $C_\zeta$,
% we have
% \begin{equation*}
%   \Omega_\bv^n
%     \leq \left(\frac{\Im S_{\bv \bv}}{N\eta}
%            \right)^n  + \Psi^n,
%          \end{equation*}
But $\Omega_\bv^n \prec \left(\frac{\Im S_{\bv \bv}}{N\eta}
           \right)^n  + \Psi^n$,
which establishes
%\eqref{eq:m-bd}.
\begin{equation}
  \label{eq:m-bd}
   |\mathcal{A}_{\gamma m}'|
      \prec \mathcal{E}_{ab} N^{-m \epsilon} \left[ (L_\bv S)^n
         + \left( \frac{N^{\epsilon}\Im S_{\bv \bv}}{N\eta}
           \right)^n
           +  (N^{\epsilon} \Psi)^n \right].
\end{equation}

To turn this into a bound on expectations, we use 
\cite[][Lemma 7.1]{benaych2018}: the right side is certainly larger than
$N^{-2-n\tau - n/2}$, and so it suffices to check that
$E |\mathcal{A}_\gamma'|^2 \leq N^{C_2}$ for some constant $C_2$.
% Let $\Xi$ be the $2\zeta$-high probability event on which (0.4)
% holds. To bound $\E  [\mathcal{A}_\gamma' \mathbf{1}(\Xi^c)]$,
Note first
that \eqref{eq:W3} implies that $\E |V_{ab}|^k \leq C_k$.
%(2.1) implies that $\E |V_{ab}|^k \leq
%(C_{\vartheta}k)^{C_\vartheta k}$.
Using this and the deterministic bounds $\| R \|, \| S \| \leq N$, we
find successively by rough bounds that
\begin{align}
  \E | \prod_{i=1}^m (Y_{k_i})_{\bv \bv}|^2
  & \leq (2N)^{k+2m} C_{2k}, \notag \\
  \intertext{and that}
  \E \mathcal{A}_\gamma^{\prime 2} 
  & \leq (CN)^{Cn}.
  % \intertext{and finally that}
  % | \E \mathcal{A}_\gamma' \mathbf{1}(\Xi^c) |
  % & \leq (N \phi^\zeta)^{C \phi^\zeta} \exp (-c \phi^{2 \zeta})
  % \leq N^{-3} \Psi^n,
      \label{eq:bad-event}
\end{align}
%say, which is much smaller than the right side of (3.33). 
Thus we may take expectations to conclude that for
$N \geq N(\epsilon, \tau, n)$, and for $1 \leq m \leq n$,
\begin{equation*}
  |\E \mathcal{A}_{\gamma m}'|
    \leq \mathcal{E}_{ab} N^{\epsilon/4 -m \epsilon} \left[ \E (L_\bv S)^n
         + \E \left( \frac{N^{\epsilon} \Im S_{\bv \bv}}{N\eta}
           \right)^n
           +  (N^{\epsilon} \Psi)^n \right].
\end{equation*}

At this point we focus on the specific case $L_\bv S = \Im S_{\bv \bv}$.
On $\mathbf{S}(\tau)$ we have $N \eta \geq N^\tau > N^{\epsilon}$ and so
in this case the previous display along with $\Psi \lesssim \Phi$ implies that,
\begin{align*}
  \E \mathcal{A}_{\gamma m}'
   % & \prec  \mathcal{E}_{ab} N^{\epsilon/4 -m\epsilon} [ \E (\Im S_{\bv \bv})^n +
   % \Psi^n ] \\
   \textcolor{blue}{\leq} %\prec  
   \mathcal{E}_{ab} N^{\epsilon/4 -m\epsilon} [ \E (\Im S_{\bv \bv})^n +
   (N^\epsilon \Phi)^n ]. 
\end{align*}
Summing over $m$
% , and noting \eqref{eq:bad-event},
we obtain (3.33a)
and hence (3.20a) for $L_\bv S = \Im S_{\bv \bv}$. 

For the remaining cases of $L_\bv$, we now use the bound (3.20a) just
established for $\Im S_{\bv \bv}$ in
% \eqref{eq:m-bd}, or more precisely
$\E \mathcal{A}_{\gamma m}'  $, along with (KY3.4) to bound
\begin{equation*}
  \E \left( \frac{\Im S_{\bv \bv}}{N\eta}  \right)^n
  \prec \left( \frac{\Phi}{N\eta}  \right)^n
  \prec  \Psi^{2n}
  \prec  \Psi^n,
\end{equation*}
so that
\begin{equation*}
  |\E \mathcal{A}_{\gamma m}'|
    \textcolor{blue}{\leq} %\prec
    \mathcal{E}_{ab}  N^{\epsilon/4 -m\epsilon/4} [ \E (L_\bv S)^n +
     (N^{\epsilon} \Psi)^n ],
\end{equation*}
and hence (3.20a) and (3.33a) follow for the remaining cases of $L_\bv$
as well. 

%\bigskip
%\textit{Outside the spectrum: Proof of \eqref{eq:outside-spec}.} \
%The extension \eqref{eq:outside-spec} is proved exactly as in Case 1
%of \cite[][Theorem 10.3]{benaych2018}, simply substituting appeals to
%the isotropic law \eqref{isotropic law} in place of those to the
%entrywise law \cite[][Theorem 2.6]{benaych2018}.
%[Note also that the extra term $(\kappa+\eta)^2$ in 
%\cite[][(10.3)]{benaych2018} is an improvement to the bound that is
%only relevant when $\kappa + \eta \geq 1$, which is not relevant for us.]
%% [see also the footnote that follows (10.6) there.]
%\bigskip

It remains to prove \eqref{isotropic delocalization}. The proof is immediate and very similar to that of (KY2.14). Using \eqref{isotropic law}, we obtain
\[
\eta^{-1}|\mathbf{v}^\ast \mathbf{u}^{(j)}|^2\leq \sum_i\frac{\eta|\mathbf{v}^\ast \mathbf{u}^{(i)}|^2}{(\lambda_j-\lambda_i)^2+\eta^2}=\Im G_{\mathbf{vv}}(\lambda_j+\im \eta)\prec 1
\]
uniformly in $N^{-1+\tau}\leq \eta\leq \tau^{-1}$ for any positive $\tau$. This yields \eqref{isotropic delocalization}. In this derivation, we implicitly used Remark KY2.4 that the overwhelming probability bounds on $\mathbf{v}^\ast G(z)\mathbf{w}-s_{sc}(z)\mathbf{v}^\ast\mathbf{w}$ hold simultaneously for all $z\in\mathbf{S}(\tau)$.

\subsection{Extension of spiked CLT to the critical case for G(U/O)E}
\label{sec:proofs-section-6}

Let $W_{h,N} =W_N+ h\mathbf{v}\mathbf{v}^{\ast }$ be the spiked (scaled) GUE or GOE with spike $%
h\in \left[ 0,\infty \right) $ along the direction $\mathbf{v} \in \mathbb{C}^{N}$ for GUE and $\mathbf{v} \in \mathbb{R}^{N}$ for GOE, and $\left\Vert \mathbf{v}
\right\Vert =1$ in both cases. Since the joint distribution of the elements of $W_{N}$ is
invariant with respect to transformations $W_{N}\rightarrow UW_{N}U^{\ast },$
where $U$ is any unitary matrix for GUE and any orthogonal matrix for GOE, the joint distribution of the 
\emph{eigenvalues } of $W_{h,N}$ does not depend on the exact value
of vector $\mathbf{v} .$ Therefore, without loss of generality, it will be
convenient to set $\mathbf{v} =\left( 0,...,0,1\right)^{*}$.

For GUE and GOE, the tridiagonalization algorithm does not change the bottom right value, see proof of Proposition~{7} in \cite{tao2012central}. Hence, the analogue of \eqref{tridiag form earlier} for $W_{h,N}$ is%
\[
\sqrt{N}\widehat{W}_{h,N}=\left( 
\begin{array}{cccc}
a_{1} & b_{1} &  &  \\ 
b_{1} & \ddots  & \ddots  &  \\ 
& \ddots  & a_{N-1} & b_{N-1} \\ 
&  & b_{N-1} & a_{N}+\sqrt{N} h %
\end{array}%
\right) = \sqrt{N}\widehat{W}_N + \sqrt{N} h \mathbf{v} \mathbf{v}^{*} \, .
\]%
Recall the definition of sequence \( R_{2}, \dots, R_{N} \) for the
tridiagonal \(\widehat{W}_N\). Since $\widehat{W}_{h,N}$ differs only in the
lower right element, it will have the corresponding ratio sequence \(
R_2, \dots, R_{N-1}, R_{h,N} \),
with the only difference coming from the fact that the very last step of %the
recursion \eqref{R equation compact} now becomes%
\begin{eqnarray*}
R_{h,N} =\alpha _{N}+\frac{h}{\theta _{N}r_{N}} -\gamma
  _{N}+\frac{\gamma _{N}+\beta _{N}-\delta_{N}}{1-R_{N-1}}
= R_{N}
+\frac{h}{\theta _{N}r_{N}}.
\end{eqnarray*}%
Therefore, denoting \( D_{h,N} = \det\left( W_{h, N} - 2 \theta_N \right) \), we get that
\begin{align}
  \log \left\vert D_{h,N}\right\vert
  &=\log |D_{N}| +\log \left\vert
      \frac{1-R_{h,N}}{1-R_{N}}\right\vert                                       
    = \log \left\vert D_{N} \right\vert +
\log \left\vert 1-\frac{h/\left( \theta _{N}r_{N}\right) }{1-R_{N}}\right\vert
                                                     \, . \label{old-new}
\end{align}%

From \eqref{eq:rmdef}, \eqref{eq:thetam2} and Lemma~\ref{Lemma Egor R},
\begin{equation*}
  \theta_N = 1 + \frac{1}{2}\sigma_N N^{-2/3}, \qquad
  r_N = 1 + \sigma_N^{1/2} N^{-1/3}, \qquad
  R_N = o_\Pr(N^{-1/3}),
\end{equation*}
which implies that for fixed $h$,
\begin{equation*}
  \log \left\vert 1-\frac{h/\left( \theta _{N}r_{N}\right) }{1-R_{N}}\right\vert
 =
 \begin{cases}
   \log |1-h| + o_\Pr(\sigma_N^{1/2} N^{-1/3})   & h \neq 1 \\
   -\tfrac{1}{3} \log N + \tfrac{1}{2} \log \sigma_N + o_\Pr(1)
   & h = 1.
 \end{cases}
\end{equation*}

% If $h\neq 1$, so that the spike is either \textit{sub-} or \textit{super-critical}, we have%
%% \[
%% \log \left\vert 1-\frac{h/\left( \theta _{N}r_{N}\right) }{1-R_{N}}%
%% \right\vert =\log \left\vert 1-h\right\vert +o_{\Pr}(1), 
%% \]%
%% and therefore, Theorem \ref{theorem main} continues to hold.

%% If $h=1$, so that the spike is \emph{critical}, then%
%% \begin{eqnarray*}
%% \log \left\vert 1-\frac{h/\left( \theta _{N}r_{N}\right) }{1-R_{N}}%
%% \right\vert  &=&\log \left\vert 1-\frac{1/\left( \theta _{N}r_{N}\right) }{%
%% 1-R_{N}}\right\vert  \\
%% &=&\log \left\vert \frac{\sigma_N^{1/2}N^{-1/3}-R_{N}+o(N^{-1/3})}{1-R_N}\right\vert 
%% \end{eqnarray*}%
%% By Lemma~\ref{Lemma Egor R}, 
%% $\vert R_{N} \vert =o_{\Pr}\left( N^{-1/3}\right),$
%% so we have%
%% \[
%% \log \left\vert 1-\frac{h/\left( \theta _{N}r_{N}\right) }{1-R_{N}}%
%% \right\vert=-\frac{1}{3}\log N +\frac{1}{2}\log \sigma_N +o_{\Pr}(1),
%% \]
%% and hence it results in an extra shift \( -\frac{1}{3} \log N \).
%%and hence, by Corollary \ref{theorem subordinate} and \eqref{old-new}
%%\[
%%\frac{\log \left\vert E_{N}^{(new)}\right\vert +\frac{1}{3}\log N}{\sqrt{%
%%\frac{\alpha}{3}\log N}}\overset{d}\rightarrow \mathcal{N}\left( 0,1\right). 
%%\]

Hence there is an extra shift \( -\frac{1}{3} \log N \) when $h=1$.
Combining Theorem~\ref{theorem main} and \eqref{old-new}, we have the following theorem.

\begin{proposition}
  \label{theorem spiked away}
Let $D_{h,N}$ be the determinant of
$W_{h,N}-2\theta_{N}$, where $2\theta_N=E=2+N^{-2/3}\sigma_{N}$
with
$\left(
\log \log N\right) ^{2}\ll \sigma_{N}\ll \left( \log N\right) ^{2}.$
Then,
\begin{equation*}
  \big(\log |D_{h,N}|-\mu_N+ \mathbf{1}_{\{h=1\}}\tfrac{1}{3}\log
  N \big)/\tilde{\tau}_N 
  \overset{d}{\rightarrow} \mathcal{N}\left( 0,1\right) .
\end{equation*}
%% \[
%% \frac{\log \left\vert \mathcal{D}_{h,N}\right\vert -N/2  + \frac{\alpha-1}{6}\log N
%% -\sigma_N N^{1/3}+\frac{2}{3}\sigma_N^{3/2}+\mathbf{1}_{\{h=1\}}\frac{1}{3}\log N}{\sqrt{\alpha\log \frac{%
%% \theta_{N}+\sqrt{\theta_{N}^{2}-1}}{2\sqrt{\theta_{N}^{2}-1}}}}%
%% \overset{d}{\rightarrow} \mathcal{N}\left( 0,1\right) .
%% \]
\end{proposition}

Note that, for $h\neq 1$, this proposition is generalized by  \cref{prop:zero-diag spiked} to Wigner matrices and a wider range for the local singularity parameter $-\gamma\leq \sigma_N\ll \log^2 N$.

\end{appendix}
%\begin{appendix}
%\section{???}
%
%\section{???}
%
%\end{appendix}

\end{document}